\documentclass[]{scrartcl}

\usepackage[utf8]{inputenc}  % Only for pdflatex, and optional in newer LaTeX
\usepackage[USenglish]{babel}
\usepackage{csquotes}

\usepackage[a4paper,top=22mm,bottom=20mm,inner=25mm,outer=20mm]{geometry}

\usepackage[%
  backend=bibtex,bibencoding=ascii,
  style=numeric-comp,
  giveninits=true, uniquename=init, %abbreviate first names
  natbib=true,
  url=false,
  doi=true,
  isbn=false,
  backref=false,
  maxnames=99
  ]{biblatex}
\usepackage{amsmath}
\allowdisplaybreaks
\numberwithin{equation}{section}
\usepackage{amssymb}
\usepackage{commath}
\usepackage{tabularx}
\usepackage{mathtools}
\usepackage{bbm}
\usepackage{nicefrac}
\usepackage{subdepth}

\usepackage{algorithm}
\usepackage{algorithmicx}
\usepackage{algpseudocode}
\usepackage{xpatch}
\makeatletter
\xpatchcmd{\algorithmic}
  {\ALG@tlm\z@}{\leftmargin\z@\ALG@tlm\z@}
  {}{}
\makeatother

\usepackage{amsthm}
\usepackage{thmtools}
\usepackage{etoolbox}
\makeatletter
\patchcmd{\thmt@setheadstyle}
 {\bgroup\thmt@space}
 {\thmt@space}
 {}{}
\patchcmd{\thmt@setheadstyle}
 {\egroup\fi}
 {\fi}
 {}{}
\makeatother
\declaretheoremstyle[
  bodyfont=\normalfont\itshape,
  headformat=\NAME\ \NUMBER\NOTE,
]{myplain}
\declaretheoremstyle[
  headformat=\NAME\ \NUMBER\NOTE,
]{mydefinition}
\newcommand{\envqed}{{\lower-0.3ex\hbox{$\triangleleft$}}}
\declaretheorem[style=myplain,numberwithin=section]{theorem}

\declaretheorem[style=mydefinition,numberlike=theorem,qed=\envqed]{definition}

\usepackage[plainpages=false,pdfpagelabels,hidelinks,unicode]{hyperref}

\usepackage{color}
\usepackage{graphicx}
\usepackage[small]{caption}
\usepackage{subcaption}

\begingroup\expandafter\expandafter\expandafter\endgroup
\expandafter\ifx\csname pdfsuppresswarningpagegroup\endcsname\relax
\else
\fi

\usepackage{booktabs}
\usepackage{rotating}
\usepackage{multirow}

\usepackage{enumitem}

\usepackage[T1]{fontenc}
\usepackage{newpxtext,newpxmath}
\usepackage{accents}

\usepackage{lineno}

\usepackage{xparse}

\let\epsilon\varepsilon
\let\phi\varphi
\let\rho\varrho

\usepackage{xspace}

\newcommand{\emh}{{e - \frac{1}{2}}}
\newcommand{\eph}{{e + \frac{1}{2}}}
\newcommand{\imh}{{i - \frac{1}{2}}}
\newcommand{\iph}{{i + \frac{1}{2}}}

\newcommand{\poly}{\mathbb{P}}
\newcommand{\uu}{\boldsymbol{u}}
\newcommand{\uud}{\uu^{\delta}}

\newcommand{\uudhat}{\hat{\uu}^{\delta}}
\newcommand{\unum}{\uu^{\text{num}}}
\newcommand{\unumi}{\uu^{\text{num}(i)}}
\newcommand{\unumj}{\uu^{\text{num}(j)}}

\newcommand{\uep}{\uu_{e, p}}
\newcommand{\bss}{\boldsymbol{s}}
\newcommand{\avg}[1]{\overline{#1}}
\newcommand{\au}{\avg{\uu}}
\newcommand{\pf}{\boldsymbol{f}}
\newcommand{\fp}{\pf^+}
\newcommand{\fm}{\pf^-}
\newcommand{\fnum}{\pf^\text{num}}
\newcommand{\fnump}{\pf^{\text{num}+}}
\newcommand{\fnumm}{\pf^{\text{num}-}}
\newcommand{\fnumep}{\pf^{e,\text{num}+}}
\newcommand{\fnumem}{\pf^{e,\text{num}-}}
\newcommand{\fnumepo}{\pf^{e+1,\text{num}-}}
\newcommand{\fnumepm}{\pf^{e,\text{num}\pm}}
\newcommand{\fnumpm}{\pf^{\text{num}\pm}}

\newcommand{\fnumpj}{\pf^{(j)\text{num}+}}
\newcommand{\fnummj}{\pf^{(j)\text{num}-}}
\newcommand{\fnumpmj}{\pf^{(j)\text{num}\pm}}

\newcommand{\fnumncpmj}{\pf^{(j)\text{num}^\pm_{\text{nc}}}}

\newcommand{\fnumncp}{\pf^{\text{num}^+_{\text{nc}}}}
\newcommand{\fnumncm}{\pf^{\text{num}^-_{\text{nc}}}}
\newcommand{\fnumncpm}{\pf^{\text{num}^\pm_{\text{nc}}}}

\newcommand{\F}{\boldsymbol{F}}
\newcommand{\Ftotp}{\F^{\text{tot}+}}
\newcommand{\Ftotm}{\F^{\text{tot}-}}
\newcommand{\Ftotpm}{\F^{\text{tot}\pm}}

\newcommand{\ftotpj}{\pf^{(j)\text{tot}+}}
\newcommand{\ftotmj}{\pf^{(j)\text{tot}-}}
\newcommand{\ftotpmj}{\pf^{(j)\text{tot}\pm}}

\newcommand{\ftotp}{\pf^{\text{tot}+}}
\newcommand{\ftotm}{\pf^{\text{tot}-}}
\newcommand{\ftotpm}{\pf^{\text{tot}\pm}}
\newcommand{\Fnum}{\F^\text{num}}
\newcommand{\Fnump}{\F^{\text{num}+}}
\newcommand{\Fnumm}{\F^{\text{num}-}}
\newcommand{\Fnumpm}{\F^{\text{num}\pm}}
\newcommand{\FnumHOpm}{\F^{\text{HO}\pm}}

\newcommand{\Fnumncpm}{\F^{\text{num}^\pm_{\text{nc}}}}

\newcommand{\Fncpm}{\F^{{\text{nc}}\pm}}
\newcommand{\Bloc}{\mathcal{B}^\delta}
\newcommand{\myvector}[1]{\text{{\hspace{0.1em}}#1{\hspace{0.15em}}}}
\newcommand{\vD}{\myvector{D}}
\newcommand{\vu}{\myvector{u}}

\renewcommand{\vec}[1]{\underline{#1}}

\newcommand{\rev}[1]{{\color{black}#1}}

\newcommand{\discf}{\ensuremath{\pf_h^\delta}}
\newcommand{\discfref}{\ensuremath{\hat{\pf}_{h,e}^\delta}}
\newcommand{\discfrefe}{\ensuremath{\hat{\pf}_h^\delta}}
\newcommand{\discF}{\ensuremath{\F_h^\delta}}
\newcommand{\discFref}{\ensuremath{\hat{\F}_{h,e}^\delta}}
\newcommand{\discFrefe}{\ensuremath{\hat{\F}_h^\delta}}

\newcommand{\bzero}{\boldsymbol{0}}
\newcommand{\bx}{\boldsymbol{x}}
\newcommand{\uU}{\boldsymbol{U}}
\newcommand{\bE}{\boldsymbol{E}}
\newcommand{\bv}{\boldsymbol{v}}
\newcommand{\bp}{\boldsymbol{p}}

\newcommand{\bc}{\boldsymbol{c}}
\newcommand{\bb}{\boldsymbol{b}}
\newcommand{\bB}{\boldsymbol{B}}

\newcommand{\pg}{\boldsymbol{g}}

\newcommand{\fg}{\boldsymbol{g}}

\newcommand{\bS}{\boldsymbol{S}}

\newcommand{\pdx}{\partial_x}

\newcommand{\dfrx}{\partial_x^\text{FR}}
\newcommand{\dlocx}{\partial_x^\text{loc}}
\newcommand{\nc}{M}
\newcommand{\ad}{P}
\newcommand{\Uad}{\mathcal{U}_{\textrm{ad}}}
\newcommand{\re}{\mathbb{R}}
\newcommand{\ud}{\textrm{d}}
\newcommand{\dv}[2]{{\frac{{\ud}#1}{{\ud}#2}}}
\newcommand{\df}[2]{{\frac{{\ud}#1}{{\ud}#2}}}
\newcommand{\xep}{x^e_p}
\newcommand{\vR}{\myvector{R}}
\newcommand{\half}{\frac{1}{2}}
\newcommand{\pph}{{p + \frac{1}{2}}}
\newcommand{\pmh}{{p - \frac{1}{2}}}

\newcommand{\Nmh}{{N-\frac{1}{2}}}
\newcommand{\Nph}{{N+\frac{1}{2}}}

\newcommand{\uez}{\boldsymbol{u}_{e, 0}}

\newcommand{\ueN}{\boldsymbol{u}_{e, N}}
\newcommand{\fx}{{\pf}_{1}}
\newcommand{\fy}{{\pf}_{2}}
\newcommand{\B}{\boldsymbol{B}}

\newcommand{\Diss}{\mathcal{D}}
\newcommand{\vel}{v}

\newcommand{\cE}{E}

\newcommand{\Eonetwo}{E_{12}}
\newcommand{\Etwotwo}{E_{22}}

\newcommand{\cR}{p}

\newcommand{\SE}{E}
\newcommand{\bSE}{\boldsymbol{\SE}}
\newcommand{\SP}{\mathcal{P}}
\newcommand{\bSP}{\boldsymbol{\SP}}
\newcommand{\SR}{\mathcal{R}}
\newcommand{\bSR}{\boldsymbol{\SR}}
\newcommand{\btop}{\boldsymbol{b}}

\newcommand{\Nabla} {\vec{\nabla}}

\newcommand{\hydroEner}{\mathcal{E}}
\newcommand\threeMatrix[1]{\underline{ #1}}
\newcommand\state[1]{\mathbf{#1}}

\newcommand{\mat}[1]{\underline{\mathbf{#1}}}
\newcommand{\amsvect}{%
  \mathpalette {\overarrow@\vectfill@}}
\def\vectfill@{\arrowfill@\relbar\relbar{\raisebox{-3.81pt}[\p@][\p@]{$\mathord\mathchar"017E$}}}
\newcommand{\vecr}[1]{\boldsymbol{#1}}

\newcommand{\orcid}[1]{ORCID:~\href{https://orcid.org/#1}{#1}}
\usepackage{authblk}

\newenvironment{keywords}{\par\textbf{Key words.}}{\par}

\title{Compact Runge-Kutta flux reconstruction methods for non-conservative hyperbolic equations}

\author[1]{Arpit~Babbar\thanks{\orcid{0000-0002-9453-370X}}}
\affil[1]{Institute of Mathematics, Johannes Gutenberg University Mainz, Staudingerweg 9, 55128 Mainz, Germany}

\author[1]{Hendrik~Ranocha\thanks{\orcid{0000-0002-3456-2277}}}

\makeatletter
\hypersetup{pdfauthor={Arpit~Babbar, Hendrik~Ranocha}} % TODO: authors
\hypersetup{pdftitle={Compact Runge-Kutta flux reconstruction methods for non-conservative hyperbolic equations}} %TODO: title
\makeatother

\begin{document}

\maketitle

\begin{abstract}
\noindent
  Compact Runge-Kutta (cRK) Flux Reconstruction (FR) methods are a variant of RKFR methods for hyperbolic conservation laws with a compact stencil including only immediate neighboring finite elements.
We extend cRKFR methods to handle hyperbolic equations with stiff source terms and non-conservative products.
To handle stiff source terms, we use IMplicit EXplicit (IMEX) time integration schemes such that the implicitness is local to each solution point, and thus does not increase inter-element communication.
Although non-conservative products do not correspond to a physical flux, we formulate the scheme using \textit{numerical fluxes} at element interfaces.
We use similar numerical fluxes for a lower order finite volume scheme on subcells of each element, which is then blended with the high order cRKFR scheme to obtain a robust scheme for problems with non-smooth solutions.
Combined with a \textit{flux limiter} at the element interfaces, the subcell based blending scheme preserves the physical admissibility of the solution, e.g., positivity of density and pressure for compressible Euler equations.
The procedure thus leads to an admissibility preserving IMEX cRKFR scheme for hyperbolic equations with stiff source terms and non-conservative products.
The capability of the scheme to handle stiff terms is shown through numerical tests involving Burgers' equations, reactive Euler's equations, and the ten moment problem.
The non-conservative treatment is tested using variable advection equations, shear shallow water equations, the GLM-MHD, and the multi-ion MHD equations.

\end{abstract}

% TODO: keywords

%TODO: MSC
% \begin{AMS}
%   65M06, % NA, PDEs, IVPs, IBVPs: Finite difference methods for initial value and initial-boundary value problems involving PDEs
%   65M20, % NA, PDEs, IVPs, IBVPs: Method of lines for initial value and initial-boundary value problems involving PDEs
%   65M70  % NA, PDEs, IVPs, IBVPs: Spectral, collocation and related methods
% \end{AMS}

\section{Introduction}

Compact Runge-Kutta Discontinuous Galerkin methods (cRKDG) were introduced in~\cite{chen2024} as a class of high-order space-time discretizations for solving hyperbolic conservation laws with more compact stencils than the standard Runge-Kutta Discontinuous Galerkin (RKDG) methods.
In~\cite{babbar2025crk}, a cRK scheme in a time averaged flux reconstruction (FR) framework was introduced.
The time averaged FR framework of~\cite{babbar2025crk} enabled provable admissibility preservation by using a flux limiter.
In addition,~\cite{babbar2025crk} required only one numerical flux computation per time step, in contrast to the original cRKDG scheme of~\cite{chen2024} which required as many numerical flux computations as the number of stages in the Runge-Kutta method.
We extend the physical admissibility-preserving compact Runge-Kutta flux reconstruction (cRKFR) method of ~\cite{babbar2025crk} to hyperbolic systems of the form
\begin{equation}
\uu_t + \nabla \cdot \pf \left( \uu \right) + \bB
\left( \uu \right) \nabla \uu = \bss \left( \uu \right), \quad \pf = \{ \pf_i \}_{i: 1}^d, \quad \bB
\left( \uu \right) \nabla \uu = \{ \bB_i \partial_{x_i}\uu \}_{i: 1}^d,
\label{eq:general.hyperbolic.equation}
\end{equation}
where $\pf_i' \left( \uu \right) + \bB_i \left( \uu \right)$ is a
diagonalizable matrix with real eigenvalues for $1 \leq i \leq d$, where $d$ denotes the physical dimension.
The term $\bB \left( \uu \right) \nabla \uu$ is called a \textit{non-conservative product}, and it usually cannot be written as the divergence of a flux.
The source term $\bss$ in~\eqref{eq:general.hyperbolic.equation} can be
\textit{stiff}, i.e., it evolves at time scales that are different from the
advective part of the equation and thus requires significantly smaller time steps when solved
with a fully explicit method.

In the pioneering works of DG/FR methods~\cite{Cockburn1991,Cockburn1989a,Cockburn1989,Huynh2007}, these spectral element methods are used as the spatial discretization to obtain a semi-discrete system.
Following the commonly used method of lines (MOL), the semi-discrete system is then solved using a multistage Runge-Kutta (RK) method.
This requires applying the FR/DG method at every RK stage, necessitating inter-element communication at each stage.
The alternative to the multistage approach are single-stage solvers which specify a complete space-time discretization.
We briefly mention some well-known single-stage methods, an in-depth review can be found in~\cite{babbar2022lax,babbar2024thesis}.
ADER (Arbitrary high order DERivative) schemes are one class of single-stage evolution methods which originated as ADER Finite Volume (FV) schemes~\cite{Titarev2002,Titarev2005} and are also used as ADER-DG schemes~\cite{dumbser2008,dumbser2014}.
Another class of single-stage evolution methods are the Lax-Wendroff schemes which were originally proposed in the finite difference framework in~\cite{Qiu2003} with the WENO approximation of spatial derivatives~\cite{Shu1989} and later extended to the DG framework in~\cite{Qiu2005b,sun2017stability}.
Compact Runge-Kutta (cRK) DG/FR methods~\cite{chen2024,babbar2025crk} are a recent class of schemes that perform single-stage evolution by doing the inter-element communication only in the last Runge-Kutta stage.

The earlier works on cRK methods~\cite{chen2024,liu2024,babbar2025crk} have focused on hyperbolic equations in a conservation law form
with applications like Computational Fluid Dynamics (CFD), astrophysics and weather modeling.
However, there are also several applications involving hyperbolic equations which cannot be written in a conservative form and require~\eqref{eq:general.hyperbolic.equation} to be solved.
Non-conservative products arise for example in the shear shallow water equations~\cite{Teshukov2007,Richard2012,chandrashekar2020ssw}, magnetohydrodynamics (MHD) equations with divergence cleaning~\cite{dedner2002,derigs2018}, and multi-ion MHD equations~\cite{toth2010multi,ramirez2025}.
These are the equations which have been used to test the proposed scheme in this paper.
However, the scheme is relevant for general hyperbolic equations with non-conservative products which additionally include applications like multilayer shallow water equations~\cite{audusse2005,castro2006}, two-phase flow models~\cite{baer1986,saurel1999}, and general relativity~\cite{Reula1998,dumbser2018}.

The most common examples of hyperbolic equations with stiff source terms include chemical reactions~\cite{benartzi1989,bihari1999}, but there are other applications like multilayer shallow water, e.g., when there is a viscosity coefficient~\cite{audusse2005}, problems involving friction~\cite{marcati1990,bouchut2007}, relaxation processes~\cite{jin1996}, and relaxation numerical methods~\cite{toro2014}.

Several phenomena governed by hyperbolic equations exhibit complexity and scale, often necessitating substantial computational resources.
This makes the development of efficient numerical methods for solving these PDE an important area of research. Higher-order methods, because of their high arithmetic intensity, are  thus particularly relevant in the current state of memory-bound High Performance Computing (HPC) hardware~\cite{attig2011,subcommittee2014}.
However, low order methods still continue to be widely used in practice because of their robustness, especially for problems with discontinuities.
In this work, the blending limiter and the subcell based blending scheme of~\cite{babbar2025crk} are extended to obtain a provably admissibility-preserving scheme for the general hyperbolic equations~\eqref{eq:general.hyperbolic.equation}.

\subsection{Related works}

Discontinuous Galerkin (DG) methods were introduced in~\cite{reed1973} for neutron transport equations and developed for hyperbolic conservation laws by Cockburn and Shu and others, see~\cite{Cockburn1991,Cockburn1989a,Cockburn1989} and other references in~\cite{cockburn2000}.
The flux reconstruction (FR) method introduced by Huynh~\cite{Huynh2007} is a finite-element type high-order method which is quadrature-free. The key idea in this method is to construct a continuous flux approximation and then use collocation at solution points which leads to an efficient implementation that can exploit optimized matrix-vector operations and vectorization capabilities of modern CPUs. The continuous flux approximation requires a correction function whose choice affects the accuracy and stability of the method~\cite{Huynh2007,Vincent2011a,Vincent2015,Trojak2021}; by properly choosing the correction function and solution points, the FR method can be shown to be equivalent to some discontinuous Galerkin and spectral difference schemes~{\cite{Huynh2007,Trojak2021}}.

The above-mentioned works on FR do not apply to non-conservative products.
One novelty of this work is to show how the idea of FR can be extended to hyperbolic equations with non-conservative products.
Path conservative schemes~\cite{castro2006,pares2006} based on the theory of Dal Maso, LeFloch and Murat~\cite{maso1995} are widely used to solve hyperbolic equations with non-conservative products.
The Weighted Essentially Non-Oscillatory (WENO) schemes of~\cite{Shu1989} were used in~\cite{castro2006} to obtain high order accuracy for solving hyperbolic equations with non-conservative products in a finite volume framework.
More recently, the work of~\cite{Balsara2024} developed WENO schemes in a finite difference framework for hyperbolic equations with non-conservative products.
The application of discontinuous Galerkin methods for hyperbolic equations with non-conservative products was first performed in~\cite{dumbser2009} using the Arbitrary DERivative (ADER) approach which is a single stage method based on solving a local space-time predictor to obtain high order accuracy~\cite{dumbser2008a}.
Our treatment of non-conservative terms is not based on path conservative schemes, but rather on a formal integration by parts treatment of the non-conservative products.
This gives a scheme that is similar to the schemes of~\cite{derigs2016,bohm2020,ramirez2021,ramirez2025} where non-conservative products are treated without using path conservative schemes.
For degree $N=0$, the treatment simply involves performing a central difference approximation with a Rusanov type~\cite{rusanov1962} dissipation. Such a finite volume scheme has been used in~\cite{yadav2023}.

The general idea for treating stiffness is to use implicit time stepping schemes.
However, using implicit time stepping schemes for the entire hyperbolic equation~\eqref{eq:general.hyperbolic.equation} will require solving a global non-linear systems of equations.
\rev{This is required when the stiffness is also present in the flux terms, e.g., the Quinpi method for hyperbolic conservation laws~\cite{quinpi2023}.
However, it is not efficient when the stiffness is only present in the source term.}
A natural solution, which is also taken in this work, is to use IMplicit-EXplicit (IMEX) schemes~\cite{ascher1995,pareschi2005} where the stiff source terms are treated implicitly, and the advective terms are treated explicitly.
The implicit equation is local to each solution point, and it does not require solving a global system of equations.
For the standard multistage Runge-Kutta methods, such IMEX methods for source terms can be used non-intrusively.
The treatment of stiff source terms for single stage methods requires more ingenuity.
As mentioned in~\cite{dumbser2008}, it is not possible to use the single stage Lax-Wendroff methods to be arbitrarily high order accurate and handle stiff source terms.
Thus, the Arbitrary DERivative (ADER) methods with a locally implicit space-time predictor were introduced in~\cite{dumbser2008} to handle stiff source terms.
Since the implicit equation in ADER is local to each space-time element, there is no large global system to be solved.
In this work, by exploiting the Runge-Kutta structure of the cRKFR methods, we add implicitness to the cRK approach directly using IMEX Runge-Kutta schemes of~\cite{ascher1995,pareschi2005}.
The crucial properties of the obtained single stage IMEX scheme that differ from the ADER scheme are that, the implicitness is point-wise and only in the source term.

Several problems involving hyperbolic equations consist of discontinuities, for which high order methods require additional numerical dissipation to avoid spurious oscillations~\cite{godunov1959}.
A recent review of such dissipation mechanisms, also known as limiters, for DG/FR methods can be found in~\cite{babbar2024admissibility}.
In this work, we use the subcell based blending limiter as in~\cite{babbar2024admissibility,babbar2025crk} inspired from~\cite{hennemann2021} to control spurious oscillations.
On the subcells, the finite volume scheme is obtained by performing a central difference approximation of the non-conservative term, with Rusanov type dissipation~\cite{rusanov1962}.
The resulting scheme is similar to the approach of~\cite{yadav2023}.
The \textit{numerical flux} for this scheme is no longer continuous across the subcell interfaces because of the non-conservative products.
The same is true for the high order cRKFR scheme, i.e., the high order numerical flux is discontinuous across FR element interfaces.
In order to maintain the conservation property when the contribution of non-conservative products is zero, similar to the conservative cases in~\cite{babbar2024admissibility,babbar2025crk}, at the FR element interfaces, we compute a convex combination of the high-order and low-order numerical fluxes and use it as the numerical flux for both the high-order and low-order schemes.
This enables the usage of the flux limiter of~\cite{babbar2024admissibility} to ensure admissibility preservation of the low-order scheme on subcells.
In the conservative case, the flux limiter would ensure that the numerical scheme is admissibility-preserving in means allowing the usage of the scaling limiter of~\cite{zhang2010c} to obtain an admissibility-preserving scheme, which is also how admissibility of standard Strong Stability Preserving Runge-Kutta (SSP-RK) FR/DG schemes~\cite{Gottlieb2001} is obtained~\cite{zhang2010c}.
This does not occur in the non-conservative case because the non-conservative products imply that the element mean of the high and low order schemes are not the same.
However, since the lower order scheme on subcells is made admissible by the flux limiter, further blending can be performed with the subcell based scheme to ensure admissibility at all solution points within each element.
Another strategy to obtain admissibility preservation for non-conservative systems is the MOOD approach~\cite{dumbser2014}.
Some differences between our approach and MOOD are as follows.
The MOOD approach discards the high order evolution in troubled cells, and instead uses a lower order evolution on $2N + 1$ subcells for a degree $N$ scheme.
A least squares reconstruction is used to transform the solution back to the high order polynomial representation at the next time step, and to maintain conservation.
However, our approach takes a convex combination of the high order and low order solutions, and maintains the conservation property because of the choice of subcells and the numerical fluxes.
The MOOD approach is completely a posteriori, whereas our approach uses a posteriori limiting only for admissibility preservation.
Admissibility-preserving techniques which are a priori also exist, for example monolithic convex limiting~\cite{kuzmin2020}.

\subsection{Contributions and organization of the paper}

In summary, the contribution of this work is to develop a provably admissibility-preserving compact Runge-Kutta flux reconstruction (cRKFR) scheme for solving~\eqref{eq:general.hyperbolic.equation} by treating stiff source terms with a locally implicit IMEX approach.
\rev{The current works on cRK methods~\cite{chen2024,liu2024,babbar2025crk} are limited to hyperbolic conservation laws without stiff source terms, and thus cannot be used for the general hyperbolic equations~\eqref{eq:general.hyperbolic.equation}.}

This contribution is important not only as an extension of the cRKFR scheme of~\cite{babbar2025crk} to non-conservative hyperbolic equations with stiff source terms, but also because admissibility-preserving IMEX schemes are challenging to construct in general.
\rev{To the best of our knowledge, there are no existing works that guarantee admissibility-preservation for general IMEX schemes for hyperbolic equations with stiff source terms and non-conservative products.}
The standard technique of Zhang and Shu~\cite{zhang2010c} to obtain admissibility-preserving schemes for explicit SSP-RK schemes does not directly extend to IMEX schemes or non-conservative systems.

The rest of the paper is organized as follows.
In Section~\ref{sec:crkfr.imex}, we introduce the cRKFR scheme for hyperbolic conservation laws with stiff source terms.
In Section~\ref{sec:non.cons.products}, the scheme is extended to handle non-conservative products.
Problems with discontinuities are treated with a subcell based shock capturing approach in Section~\ref{sec:blending.scheme}, which is also used to obtain an admissibility-preserving scheme as long as a first-order scheme with the same property is available.
Numerical results are presented in Section~\ref{sec:numerical.results} to demonstrate the performance of the proposed scheme.
Finally, conclusions are drawn in Section~\ref{sec:conclusions}.

\section{Compact Runge-Kutta flux reconstruction for conservation laws with source terms} \label{sec:crkfr.imex}

In this section, we describe our proposed IMplicit-EXplicit (IMEX) compact Runge-Kutta flux reconstruction (cRKFR) scheme for 1-D hyperbolic conservation laws with source terms
\begin{equation}
\label{eq:con.law.source} \uu_t + \pf (\uu)_x = \bss(t, \uu), \qquad \uu (x, t_0) = \uu_0
  (x), \qquad x \in \Omega,
\end{equation}
where the solution $\uu \in \re^{\nc}$ is the vector of conserved quantities, $\Omega$ is the physical domain, $\pf (\uu)$ is the physical flux, $\bss(t, \uu)$ is the source term which can be stiff requiring smaller time steps when used with explicit methods (Appendix~\ref{app:stiff.source}), $\uu_0$ is the initial condition, and some boundary conditions are additionally prescribed.
In many applications involving~\eqref{eq:con.law.source}, it is essential that the solution is physically sensible, i.e., it belongs to an admissible set denoted by $\Uad$.
For example, in gas dynamics with chemical reactions governed by the reactive Euler's equations~\eqref{eq:reactive.euler}, the density, pressure and chemical reactant are positive.
In case of the ten moment equations~\eqref{eq:tmp}, the pressure tensor is a positive definite matrix.
It will be assumed that the equations of interest~(\ref{eq:con.law.source}, \ref{eq:general.hyperbolic.equation}) satisfy the admissibility preservation property of the solution, i.e.,
\begin{equation} \label{eq:hyp.adm.preservation}
\uu(\cdot, t_0) \in \Uad \qquad  \implies \qquad \uu(\cdot, t) \in \Uad,\qquad t>t_0.
\end{equation}
For the models that we study, we assume that the admissible set is a convex subset of $\re^{\nc}$, and can be written as
\begin{equation}
\label{eq:uad.form} \Uad = \{ \uu \in \re^{\nc} : \ad_k (\uu) > 0, 1 \le k \le K\}.
\end{equation}
For the reactive Euler's equations~\eqref{eq:reactive.euler}, $K = 3$ and $\ad_1, \ad_2$ and $\ad_3$ are the density, reactant and pressure constraint respectively; if the density is positive then the pressure constraint is a concave function of the conserved variables.
In case of the ten moment problem~\eqref{eq:tmp}, $K = 3$ and $\ad_1, \ad_2, \ad_3$ are the density, the trace and the determinant of the pressure tensor respectively.
The determinant $\ad_3$ is not a concave function of the conserved variables~\eqref{eq:det.concave}, but the admissibility-preserving scheme developed in this work is still applicable as discussed in Section~\ref{sec:flux.limiter}.
In~\cite{babbar2025crk}, the authors developed an admissibility-preserving compact Runge-Kutta flux reconstruction (cRKFR) discretization to solve hyperbolic conservation laws.
Following~\cite{babbar2025generalized}, an admissibility-preserving treatment of the source terms using the time averaging cRKFR schemes was already given in Appendix B of~\cite{babbar2025crk}, but it was fully explicit and thus incapable of handling stiff source terms.
In this section, we extend the time averaged framework of~\cite{babbar2025crk} to handle stiff source terms by using IMplicit-EXplicit (IMEX) Runge-Kutta methods, where the flux term is treated explicitly, and the source term is treated implicitly.
The implicit equation obtained for the source term is local to each solution point, and is thus efficient as it, e.g. does not require solving a global system of equations.

\begin{figure}
\begin{center}
\begin{tabular}{cc}
\includegraphics[width=0.45\columnwidth]{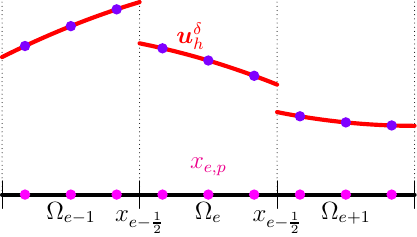} &
\includegraphics[width=0.45\columnwidth]{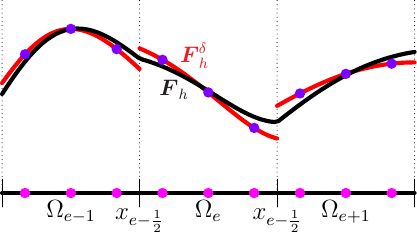}\\
(a) & (b)
\end{tabular}
\end{center}
\caption{(a) Piece-wise polynomial solution at time $t_n$, and (b)
discontinuous and continuous flux. The figure is inspired from~{\cite{babbar2022lax}}.}
\label{fig:solflux1}
\end{figure}

We now describe the notations for the finite element grid, which follow our previous works~\cite{babbar2022lax,babbar2025crk}.
The physical domain $\Omega$ is divided into disjoint elements $\{\Omega_e \}$ with
\begin{equation}
  \label{eq:reference.element} \Omega_e = [x_{\emh}, x_\eph] \qquad
  \text{and} \qquad \Delta x_e = x_\eph - x_{\emh} .
\end{equation}
Each element is mapped to a reference element, $\Omega_e \to [0, 1]$, by
\begin{equation}
x \mapsto \xi = \frac{x - x_{\emh}}{\Delta x_e} .
\label{eq:ref.map}
\end{equation}
Inside each element, we approximate the numerical solution
to~\eqref{eq:con.law.source} by $\poly_N$ functions which are degree $N \geq 0$
polynomials so that the basis for the numerical solution is
\begin{equation}
  \label{eq:fr.basis} V_h = \{v_h : v_h |_{\Omega_e} \in \poly_N \},
\end{equation}
which allows the solution to be discontinuous at element interfaces, see Figure~\ref{fig:solflux1}a.

For each basis function $v_h \in V_h$~\eqref{eq:fr.basis} and a physical element $\Omega_e$~\eqref{eq:reference.element}, we define its representative in the reference element $\hat{v}_{h, e} = \hat{v}_{h, e} (\xi)$ through the reference map~\eqref{eq:ref.map} as
\begin{equation}
  \label{eq:vh.reference} \hat{v}_{h, e} (\xi) = v_h  (x_{\emh} + \xi \Delta
  x_e) = v_h (x) .
\end{equation}
The element index $e$ will often be suppressed for brevity.
To construct
$\poly_N$ polynomials in the basis, choose $N + 1$ distinct nodes
\begin{equation}
  \label{eq:soln.points} 0 \le \xi_0 < \xi_1 < \cdots < \xi_N \le 1,
\end{equation}
which will be Gauss-Legendre (GL) or Gauss-Lobatto-Legendre (GLL) nodes, and will also be referred to as the \textit{solution points}.
There are associated quadrature weights $w_j$ such that the quadrature rule is exact for
polynomials of degree up to $2 N + 1$ for GL points and up to degree $2 N - 1$
for GLL points.
Note that the nodes and weights we use are with respect to the
interval $[0, 1]$ whereas they are usually defined for the interval $[- 1, +
1]$.
The numerical solution $\uud_h \in V_h^{\nc}$, where $\nc$ is the number of conservative variables~\eqref{eq:con.law.source}, inside an element $\Omega_e$ is given in reference coordinates as
\begin{equation}
  \label{eq:soln.poly} x \in \Omega_e : \qquad \uudhat_{h, e} (\xi, t) =
  \sum_{p = 0}^N \uu_{e, p} (t) \ell_p (\xi),
\end{equation}
where each $\ell_p$ is a Lagrange polynomial of degree $N$ given by
\begin{equation}
  \label{eq:defn.lagrange} \ell_q (\xi) = \prod_{p = 0, p \ne q}^N \frac{\xi -
  \xi_p}{\xi_q - \xi_p} \in \poly_N, \qquad \ell_q (\xi_p) = \delta_{p
   q} .
\end{equation}
The numerical methods to solve~\eqref{eq:con.law.source} require computation of
spatial derivatives.
The spatial derivatives are computed on the reference
interval $[0, 1]$ using the differentiation matrix $\vD = [D_{ij}]$ whose
entries are given by
\begin{equation}
  \label{eq:diff.matrix} D_{ij} = \ell_j' (\xi_i), \qquad 0 \le i, j \le N.
\end{equation}
Figure~\ref{fig:solflux1}a illustrates a piece-wise polynomial solution at some
time $t_n$ with discontinuities at the element boundaries.

A crucial ingredient of the cRKFR scheme and the FR scheme in general is the degree $N$ \textit{discontinuous flux} approximation $\discf$ (Figure~\ref{fig:solflux1}b), which is defined in reference coordinates for element $e$ as
\begin{equation} \label{eq:discts.flux}
\discfref(\xi) = \sum_{p=0}^N \pf(\uu_{e,p}) \ell_p(\xi).
\end{equation}
Thus, the flux reconstruction (FR) differentiation operator is defined as
\begin{equation}
\begin{gathered}
\dfrx \pf (\uud_h) = \pdx \pf_h, \quad \pf_h = \discf +
(\fnum_\eph - \fm_\eph) g_R + (\fnum_{\emh} - \fp_\emh) g_L, \\
\fm_\eph = \discfrefe (1), \quad \fp_\emh = \discfrefe (0),
\label{eq:dfrx.defn}
\end{gathered}
\end{equation}
where $g_L, g_R \in \poly_{N + 1}$ are FR correction functions~\cite{Huynh2007} satisfying
\[
g_L (0) = g_R (1) = 1, \qquad g_L (1) = g_R (0) = 0,
\]
and $\fnum_\eph$ is the numerical flux at the element interface $x_\eph$. We use
\begin{equation} \label{eq:num.flux.fr}
\fnum_\eph(\uud_h) = \half(\pf(\uu_\eph^-) + \pf(\uu_\eph^+)) - \frac {\lambda_\eph}{2} (\uu_\eph^+ - \uu_\eph^-), \qquad \uu_\eph^\pm = \uud_h(x_\eph^\pm),
\end{equation}
where $\lambda_\eph = \lambda(\uu_\eph^-,\uu_\eph^+)$ is a Rusanov/local Lax-Friedrichs~\cite{rusanov1962} estimate of the wave speed at the interface $x_\eph$ given by
\begin{equation} \label{eq:rusanov.wave.speed}
\lambda(\uu_l, \uu_r) = \max_{\uu \in \{\uu_l, \uu_r\}} \sigma(\pf'(\uu)).
\end{equation}
Here $\sigma(A)$ denotes the spectral radius of a matrix $A$.
The $\pf_h$ defined in~\eqref{eq:dfrx.defn} is called the \textit{continuous
flux approximation} in the FR literature~{\cite{Huynh2007}} as it is globally
continuous, taking the numerical flux value $\fnum_{e + 1 / 2}$ at each
element interface $e + 1 / 2$.
We choose correction functions $g_{L/R} \in
\mathbb{P}_{N + 1}$ to be $g_{\text{Radau}}$ for GL and $g_2$ for GLL nodes~{\cite{Huynh2007}}.
A semi-discretization of~\eqref{eq:con.law.source} is thus obtained as
\begin{equation}
  \label{eq:semidiscretization.dg} \partial_t \uud_h + \dfrx \pf (\uud_h) = \bss(t, \uud_h).
\end{equation}
Solving~\eqref{eq:semidiscretization.dg} using an RK method gives the original RKFR method of~{\cite{Huynh2007}}.
We now describe the novel fully discrete cRKFR IMEX scheme to solve~\eqref{eq:con.law.source}.
The $s$-stage scheme is specified by double \textit{Butcher tableaux}
\begin{equation}
  \begin{array}{c|c}
    \widetilde{\bc} & \tilde{A}\\
    \hline
    & \widetilde{\bb}^T
  \end{array}, \qquad \begin{array}{c|c}
    \bc & A\\
    \hline
    & \bb^T
  \end{array},
\label{eq:imex.butcher}
\end{equation}
where $\tilde{A} = (\tilde{a}_{ij})_{s \times s}$ is a strictly
lower triangular matrix corresponding to the explicit part.
The matrix $A =
(a_{ij})_{s \times s}$ corresponds to the implicit part.
We use diagonally implicit RK (DIRK) methods for the implicit part, i.e., the $A$ matrix is lower triangular.
The $\rev{\widetilde{\bc} = (\tilde{c}_1,
\ldots, \tilde{c}_s)}^T, \bc = (c_1, \ldots, c_s)^T$~\cite{Hairer1991, Hairer1996} are the
\textit{quadrature nodes} given by
\begin{equation}
  \tilde{c}_i = \sum_{j = 1}^{i - 1} \tilde{a}_{ij}, \qquad c_i =
  \sum_{j = 1}^i a_{ij} .
\end{equation}
The $\widetilde{\bb} = (\tilde{b}_1, \tilde{b}_2, \ldots, \tilde{b}_s)^T, \bb = (b_1, b_2, \ldots, b_s)^T$ are the \textit{weights} used to perform the final evolution with evaluations from the $s$ stages.
The choices of IMEX Runge-Kutta schemes used in this work are given in Appendix~\ref{app:imex.schemes}.
Following~\cite{chen2024,babbar2025crk}, we define a local differentiation operator $\dlocx$ where the inter-element terms from~\eqref{eq:dfrx.defn} are dropped to be
\begin{equation}
\dlocx \pf (\uud_h) (\xi_p) = \pdx \discf (\xi_p) .
\label{eq:dlocx.defn}
\end{equation}
The operator $\dlocx$ computes derivatives of the degree $N$ discontinuous flux approximation~\eqref{eq:discts.flux}.
In practice, this operation is performed with a differentiation matrix~\eqref{eq:diff.matrix}.
The compact Runge-Kutta flux reconstruction method in the time average framework applied to~\eqref{eq:con.law.source} is given by
\begin{eqnarray}
  \uu^{(i)} & = & \uu^n - \Delta t \sum_{j = 1}^{i - 1} \tilde{a}_{i j} \partial_x^{\text{loc}}  \pf_h ( \uu^{(j)} ) +
  \Delta t \sum_{j = 1}^i a_{ij}  \bss ( t_n +
  c_j \Delta t, \uu^{(j)} ), \quad i = 1, \ldots, s, \label{eq:crkfr.inner}\\
  \uu^{n + 1} & = & \uu^n
  - \Delta t \dfrx \F
  + \Delta t \bS_h^{\delta}, \label{eq:crkfr}\\
  \dfrx \F = \pdx \F_h  &=& \pdx \discF + (\Fnum_\eph - \discFrefe(1))\pdx g_R + (\Fnum_\emh - \discFrefe(0)) \pdx g_L, \label{eq:cts.F}
\end{eqnarray}
where the time averaged flux $\F$ is computed in the reference coordinates as
\begin{equation}
  \discFref = \sum_{p = 0}^N \F_{e, p} \ell_p (\xi), \qquad \F_{e, p} =
  \sum_{i = 1}^s \rev{\tilde{b}_i}  \pf (\uu_{e, p}^{(i)}), \label{eq:disc.avg.flux}
\end{equation}
and the time averaged source term $\bS$ as
\begin{equation}
\hat{\bS}_{h, e}^{\delta} = \sum_{p = 0}^N \bS_{e, p} \ell_p (\xi), \qquad \bS_{e, p} =
  \sum_{i = 1}^s b_i  \bss (t_n + c_i \Delta t, \uu^{(i)}_{e,p}) .
\label{eq:avg.source}
\end{equation}
The scheme also makes use of the time average solution to compute the numerical flux $\Fnum_\eph$~\eqref{eq:cts.F}, which is computed as
\begin{equation}
\hat{\uU}_{h,e}^\delta = \sum_{p = 0}^N \uU_{e, p} \ell_p (\xi), \qquad \uU_{e, p} =   \sum_{i = 1}^s \rev{\tilde{b}_i}  \uu_{e, p}^{(i)}. \label{eq:time.avg.soln}
\end{equation}
A formal justification for why~(\ref{eq:disc.avg.flux},~\ref{eq:avg.source},~\ref{eq:time.avg.soln}) are approximations to the time averaged flux, source and solution respectively can be found in Appendix~A of~{\cite{babbar2025crk}}.
The \textit{time averaged numerical flux} $\Fnum_\eph$ used to construct the continuous flux approximation~\eqref{eq:cts.F} is computed as
\begin{equation} \label{eq:numflux.cons}
\Fnum_\eph = \half (\F_\eph^- + \F_\eph^+) - \frac{\lambda_\eph}{2} (\uU_\eph^+ - \uU_\eph^-),
\end{equation}
where $\F_\eph^\pm, \uU_\eph^\pm$ are the approximations of the time average flux and solution~(\ref{eq:disc.avg.flux}) at the element interface $x_\eph$, and $\lambda_\eph$ is the maximum wave speed estimate at the interface $x_\eph$~\eqref{eq:rusanov.wave.speed}.
The dissipation term is chosen to be $\lambda_\eph (\uU_\eph^+ - \uU_\eph^-) / 2$ so that the scheme has the same CFL number as the cRK scheme of~\cite{chen2024} while requiring only a single numerical flux $\Fnum_\eph$ to be computed per time step, see Section 3.1 of~\cite{babbar2025crk} for details.
The idea of using the time average solution to compute the dissipation term in the numerical flux was first introduced for Lax-Wendroff schemes in~\cite{babbar2022lax} where it was shown to improve the CFL numbers.
The difference between the cRKFR scheme~\eqref{eq:crkfr} and the standard RKFR scheme lies in the fact that the inner stages that evaluate $\{\uu_h^{(i)} \}_{i = 1}^s$ in~\eqref{eq:crkfr.inner} use the local operator $\dlocx$ instead of $\dfrx$ in~\eqref{eq:semidiscretization.dg}.

We now briefly review the admissibility preserving terminologies and strategies used in~\cite{babbar2025crk,babbar2024admissibility} for conservation laws, as they can already be applied to the cRKFR IMEX scheme for solving~\eqref{eq:con.law.source} because of its time averaged flux and source term framework.
The discrete variant of the admissibility preservation property~\eqref{eq:hyp.adm.preservation} for a flux reconstruction scheme is given as follows.
\begin{definition}\label{defn:admissibility.preserving}
A flux reconstruction scheme is said to be admissibility-preserving if
\[
\uu_{e,p}^n \in \Uad \quad \forall e,p \qquad \implies \qquad \uu_{e,p}^{n+1} \in \Uad \quad \forall e,p,
\]
where $\Uad$ is the admissible set~\eqref{eq:uad.form} of the physical equation~(\ref{eq:con.law.source},~\ref{eq:general.hyperbolic.equation}).
\end{definition}
We also mention the weaker property of admissibility preservation in means which is successfully used to develop admissibility-preserving schemes for conservation laws in~\cite{zhang2010c,babbar2025crk,babbar2024admissibility}.
\begin{definition}\label{defn:admissibility.preserving.means}
A flux reconstruction scheme is said to be admissibility-preserving in means if
\begin{equation}
\uu_{e,p}^n \in \Uad \quad \forall e,p \qquad \implies \qquad \au_{e}^{n+1} \in \Uad \quad \forall e,
\end{equation}
where $\au_{e}$ is the element mean of $\uu_h$ in element $e$ computed as
\[
\au_e = \sum_{p=0}^N w_p \uu_{e,p},
\]
where $w_p$ are the quadrature weights associated with the solution points~\eqref{eq:soln.points}, and $\Uad$ is the admissible set~\eqref{eq:uad.form} of the equations of interest~(\ref{eq:general.hyperbolic.equation},~\ref{eq:con.law.source}).
\end{definition}
% Thus, it is essential to develop a different strategy to obtain admissibility preservation for non-conservative systems, which is done in Section~\ref{sec:blending.scheme}.
As we have written the cRKFR IMEX scheme in a time averaged flux and source term framework, the admissibility preservation procedure of~\cite{babbar2024admissibility,babbar2025generalized} can be applied to obtain an admissibility-preserving in means scheme (Definition~\ref{defn:admissibility.preserving.means}) as long as a first-order scheme with the same property is available.
This would be done by applying a limiter specifically on the time averaged flux and source terms.
Once the scheme preserves admissibility in means, an admissibility-preserving scheme (Definition~\ref{defn:admissibility.preserving}) can be obtained by using the scaling limiter of Zhang and Shu~\cite{zhang2010c}.
% However, for non-conservative systems~\eqref{eq:general.hyperbolic.equation}, there is no natural strategy to obtain admissibility preservation in means because of the non-conservative products.
In this work, we do not use a source term limiting strategy since it is not enough for obtaining an admissibility-preserving scheme when non-conservative products are present.
Thus, a more generally applicable strategy for admissibility preservation is described in Section~\ref{sec:flux.limiter}.

\section{Non-conservative products} \label{sec:non.cons.products}
In this section, we describe our proposed extension of the cRKFR scheme to hyperbolic equations with non-conservative products of the form
\begin{equation}
\uu_t + \pf (\uu)_x + \bB (\uu)  \uu_x
   = \bzero,
\label{eq:flux.and.non.cons}
\end{equation}
where $\pf'(\uu) + \bB(\uu)$ is a diagonalizable matrix with real eigenvalues.
The source terms are omitted for simplicity, but their inclusion is independent of the non-conservative products and is performed exactly as in Section~\ref{sec:crkfr.imex}.
For further simplicity, we begin by describing the finite volume method (FVM) for the equation without the conservative flux, i.e.,
\begin{equation}
\uu_t = - \bB (\uu)  \uu_x. \label{eq:pure.non.conservative}
\end{equation}
A natural discretization of~\eqref{eq:pure.non.conservative} is obtained by using a central difference approximation
\begin{equation}
\dv{\uu_i}{t} = - \frac{1}{\Delta x}  \bB(\uu_i) (\unum_\iph - \unum_\imh),
\qquad \unum_\iph = \unum(\uu_i, \uu_{i+1}) := \half  (\uu_i + \uu_{i + 1}).
\label{eq:semi.fv}
\end{equation}
In order to write it in a finite volume form, we introduce the following numerical fluxes
\[ \dv{\uu_i}{t} = - \frac{1}{\Delta x}  (\fnumncm_\iph - \fnumncp_\imh), \]
where $\fnumncpm_\iph = \fnumncpm(\uu_i, \uu_{i+1})$ which is defined as
\begin{equation}
\fnumncpm(\uu_-, \uu_+) = \bB({\uu_\pm})  \unum(\uu_-, \uu_+).
\label{eq:non.cons.num.flux.basic}
\end{equation}
The crucial difference between~\eqref{eq:non.cons.num.flux.basic} and the standard conservative numerical fluxes~\eqref{eq:num.flux.fr} is that each finite volume interface $\iph$ now has two numerical fluxes $\fnumncpm_\iph$, which is equivalent to saying that the numerical fluxes are no longer continuous across the interfaces.
From the common experience with conservation laws, the above scheme will need to have additional dissipation terms to be stable.
Thus, the final numerical flux for solving~\eqref{eq:pure.non.conservative} is actually chosen to be
\begin{equation}
\fnumncpm_\iph
- \frac{\lambda_\iph}{\rev{2}}  (\uu_{i + 1} - \uu_i),
\label{eq:num.flux.non.cons}
\end{equation}
where $\lambda_\iph$ is defined by the Rusanov wave speed estimate~\eqref{eq:rusanov.wave.speed} extended to non-conservative products as $\lambda(\uu_l, \uu_r) = \max_{\uu \in \{\uu_l, \uu_r\}} \sigma(\bB(\uu))$.
However, note that, unlike the central part~\eqref{eq:non.cons.num.flux.basic}, the dissipation term is actually continuous across the interface $\iph$.
This allows us to include the dissipation term of the entire system in the conservative numerical flux which is continuous across interfaces.
Thus, for the system with both conservative and non-conservative terms~\eqref{eq:flux.and.non.cons}, we define the numerical fluxes $\fnumpm_\iph$ as
\begin{equation}
\begin{split}
\fnumpm_\iph &=
\fnumpm(\uu_i, \uu_{i+1}) :=
\fnum(\uu_i, \uu_{i+1}) + \fnumncpm(\uu_i, \uu_{i+1}), \\
\fnum_\iph &= \fnum(\uu_i, \uu_{i+1}) :=
\half(\pf(\uu_i) + \pf(\uu_{i+1}))
- \frac {\lambda_\iph}{2} (\uu_{i+1} - \uu_{i}),
\end{split}
\label{eq:combined.num.flux}
\end{equation}
where the dissipation coefficient is defined as
\begin{equation}
\lambda_\iph = \max_{\uu \in \{\uu_i, \uu_{i+1}\}} \sigma(\pf'(\uu) + \bB(\uu)),
\label{eq:combined.wave.speed}
\end{equation}
where $\sigma(A)$ denotes the spectral radius of matrix $A$.
The general FR scheme to solve~\eqref{eq:flux.and.non.cons} that applies to any choice of correction functions and solution points is derived via DG methods in Appendix~\ref{app:dg.fr}.
Here we simply state it, while mentioning that it is similar to the FR scheme~\eqref{eq:semidiscretization.dg}, after the operator $\dfrx$ has been expanded using \eqref{eq:dfrx.defn}, as follows
\begin{equation}
\begin{aligned}
\dv{\uud_h}{t}
+ \dlocx \pf(\uud_h) + \bB (\uud_h) (\dlocx \uud_h)
& + [
  % \bB (\uu_\eph^-) \unum_\eph +
  \fnumm_\eph - \ftotm_\eph
  % \bB (\uu_\eph^-) \uu_\eph^- - \pf_\eph^-
  ] \partial_x g_R  + [
  % \bB (\uu_{\emh}^+) \unum_{\emh} +
  % \fnumncpm_{\emh} - \bB (\uu_{\emh}^+) \uu_{\emh}^+ - \pf_{\emh}^+
  \fnump_{\emh} - \ftotp_{\emh}
  ] \partial_x g_L = \bzero
,
\end{aligned} \label{eq:general.non.conservative.fr.one}
\end{equation}
where $\dlocx$ is the local derivative operator~\eqref{eq:dlocx.defn}
% which is used in~\eqref{eq:general.non.conservative.fr.one} to compute the derivative of the discontinuous flux~\eqref{eq:discts.flux} and the discontinuous solution~\eqref{eq:soln.poly}
, $\fnumpm_\eph$ is defined using the finite volume conservative and non-conservative numerical fluxes~(\ref{eq:combined.num.flux}, \ref{eq:non.cons.num.flux.basic}) as
\begin{equation}
\begin{split}
\fnumpm_\eph &= \fnum_\eph + \fnumncpm_\eph, \qquad \fnumncpm_\eph = \fnumncpm(\uu_\eph^-, \uu_\eph^+),\\
\fnum_\eph &= \fnum(\uu_\eph^-, \uu_\eph^+),
% \bB (\uu_\eph^{\pm})  \unum_\eph,
%=- \lambda_\eph  (\uu_{e + 1} - \uu_e)=%,
\end{split} \label{eq:combined.num.flux.fr}
\end{equation}
where $\uu_\eph^{\pm}$ are the extrapolations of the solution as in~\eqref{eq:num.flux.fr}.
The $\ftotpm_\eph$ consist of the extrapolations of the conservative and non-conservative terms at the element interface $x_\eph$ given by
\begin{equation} \label{eq:total.flux.non.cons}
\ftotpm_\eph = \bB (\uu_\eph^\pm) \uu_\eph^\pm + \pf_\eph^\pm,
\end{equation}
where, as in~\eqref{eq:dfrx.defn}, $\pf_\eph^\pm$ are extrapolations obtained from the discontinuous fluxes~\eqref{eq:discts.flux} at the element interface $x_\eph$.
The motivation for why the extrapolations of the non-conservative terms are chosen as in~\eqref{eq:total.flux.non.cons} is also given in Appendix~\ref{app:dg.fr}.
It is further shown in Appendix~\ref{app:dg.fr} that the scheme~\eqref{eq:general.non.conservative.fr.one} reduces to the first order finite volume scheme~\eqref{eq:combined.num.flux} when $N = 0$, and that for the FR scheme with Gauss-Lobatto-Legendre (GLL) solution points and $g_2$ correction function can be written in terms of the FR operators $\dfrx$ like in the conservative case~\eqref{eq:semidiscretization.dg} as
\begin{equation}
\begin{gathered}
\dv{\uud_h}{t} +  \dfrx \pf(\uud_h) + \bB (\uud_h) \dfrx  \uud_h = \bzero,\qquad  \\
\qquad \dfrx \uud_h := \partial_x  \uu_h,
\qquad \uu_h := \uud_h +
(\unum_\eph - \uu_\eph^-)  g_R +
(\unum_{\emh} - \uu_{\emh}^+)  g_L.
\end{gathered} \label{eq:fr.non.cons.gll}
\end{equation}
In order to now derive the cRKFR scheme for~\eqref{eq:flux.and.non.cons}, we show the fully discrete scheme obtained by applying an explicit RK method to~\eqref{eq:general.non.conservative.fr.one}.
Following the Butcher tableau notation~\eqref{eq:imex.butcher} and taking only the explicit part from the tableau~\eqref{eq:imex.butcher}, we obtain the fully discrete RKFR scheme for the non-conservative system~\eqref{eq:flux.and.non.cons} to be
\begin{subequations}\label{eq:rkfr.non.cons}
\begin{align}
&\begin{aligned}
\ \ \uu^{(i)}& = \uu^n -  \Delta t\sum_{j = 1}^{i-1} \tilde{a}_{ij} \Big( \dlocx \pf(\uu^{(j)}) + \bB (\uu^{(j)}) \dlocx \uu^{(j)} \\
&\qquad\qquad\qquad + \big[\fnummj_\eph - \ftotmj_\eph\big]\partial_x g_R
+ \big[\fnumpj_{\emh} - \ftotpj_{\emh}\big]\partial_x g_L \Big),
\ \ i = 1, 2, \dots, s,
\end{aligned}\label{eq:rkfr1} \\
&\begin{aligned}
\uu^{n+1} &= \uu^n - \Delta t \sum_{j = 1}^s \tilde{b}_j \Big( \dlocx \pf(\uu^{(j)} ) + \bB (\uu^{(j)}) \dlocx \uu^{(j)} \\
&\qquad\qquad\qquad + \big[\fnummj_\eph - \ftotmj_\eph\big]\partial_x g_R
+ \big[\fnumpj_{\emh} - \ftotpj_{\emh}\big]\partial_x g_L \Big),
\end{aligned}
\label{eq:rkfr2}
\end{align}\label{eq:rkfr}
\end{subequations}
where $\fnumpmj_\eph = \fnumpm(\uu_\eph^{(j)-}, \uu_\eph^{(j)+})$ are the sum of conservative and non-conservative numerical fluxes~\eqref{eq:combined.num.flux.fr} at stage $j$ evaluated as~\eqref{eq:combined.num.flux}, and $\ftotpmj_\eph$ are the extrapolations of the conservative and non-conservative terms at stage $j$ computed as in~\eqref{eq:total.flux.non.cons}.
Similar to the conservative case~\eqref{eq:cts.F} the cRKFR scheme is obtained by dropping the inter-element terms from the inner stages~\eqref{eq:rkfr1} while using FR operators on the time averages in the final update~\eqref{eq:rkfr2} as follows
\begin{eqnarray}
\uu^{(i)} & = & \uu^n -  \Delta t\sum_{j = 1}^{i-1} \tilde{a}_{ij} \Big( \dlocx \pf(\uu^{(j)}) + \bB (\uu^{(j)}) \dlocx \uu^{(j)} \Big),
\quad i = 1, \ldots, s,\\
\uu^{n + 1} & = & \uu^n
- \Delta t \pdx \discF
- \Delta t \Bloc_h
- \Delta t (\Fnumm_\eph - \Ftotm_\eph)\pdx g_R - \Delta t (\Fnump_\emh - \Ftotp_\emh) \pdx g_L,
\label{eq:crkfr.non.cons.final}
\\ \Ftotpm_\eph &=& \F_\eph^\pm  + \Fncpm_\eph,\quad \Fnumpm_\eph = \Fnum_\eph + \Fnumncpm_\eph, \label{eq:total.F.non.cons}
\end{eqnarray}
where the local time averaged flux $\discF$ in~\eqref{eq:crkfr.non.cons.final} is computed as in~\eqref{eq:disc.avg.flux}, the local time average of the non-conservative term, $\Bloc_h$, is computed as a time average of the \rev{local approximations of $\bB(\uu) \partial_x \uu$ used in~(\ref{eq:general.non.conservative.fr.one},~\ref{eq:rkfr})} to give
\begin{equation}
\Bloc_h = \sum_{j=1}^s \tilde{b}_j \bB (\uu^{(j)}) \dlocx \uu^{(j)}.
\label{eq:local.time.avg.non.cons}
\end{equation}
The conservative part of the time averaged numerical flux $\Fnum_\eph$ is computed as in~\eqref{eq:numflux.cons} with the dissipation coefficient $\lambda_\eph$ using the wave speed estimate for the full non-conservative system~\eqref{eq:flux.and.non.cons} as in~\eqref{eq:combined.wave.speed}.
The conservative part $\F_\eph^\pm$ in the extrapolation term $\Ftotpm_\eph$ in~\eqref{eq:total.F.non.cons} is computed by extrapolating the time averaged flux~\eqref{eq:disc.avg.flux} to the element interface $x_\eph$ as in~\eqref{eq:cts.F}.
The non-conservative part $\Fncpm_\eph$ of the extrapolation term~$\Ftotpm_\eph$ in~\eqref{eq:total.F.non.cons} is computed as
\[
\Fncpm_\eph = \sum_j \rev{\tilde{b}_j} \bB (\uu_\eph^{(j)\pm}) \uu_\eph^{(j)\pm}.
\]
The expression of the non-conservative time averaged numerical flux $\Fnumncpm_\eph$ that is obtained by time averaging~(\ref{eq:combined.num.flux.fr},~\ref{eq:non.cons.num.flux.basic}) directly would be
\begin{equation}
\Fnumncpm_\eph = \sum_{j=1}^s \rev{\tilde{b}_j} \fnumncpmj_\eph = \sum_{j=1}^s \rev{\tilde{b}_j} \bB (\uu_\eph^{(j), \pm})  \unumj_\eph.
\label{eq:bad.numflux.non.cons}
\end{equation}
However,~\eqref{eq:bad.numflux.non.cons} requires the inter-element communication of the inner stage values $\{\uu_\eph^{(j), \pm} \}_{j = 1}^s$ at the element interfaces between neighboring elements.
We would like to remark that this would be similar to the communication volume required for the first cRKDG scheme proposed in~\cite{chen2024}, and the ADER scheme~\cite{dumbser2008} as ADER requires communication of the local space-time predictor at the element interfaces when dealing with non-conservative products.
Thus, the increased communication volume just might be unavoidable for some practical problems.
However, in all the problems that we tried, we have reduced our inter-element communication by using
\begin{equation}
\Fnumncpm_\eph = \sum_{j=1}^s \rev{\tilde{b}_j} \bB (\uu_\eph^{n\pm})  \unumi_\eph = \bB (\uu_\eph^{n\pm}) \sum_{j=1}^s \rev{\tilde{b_j}} \unumj_\eph = \bB (\uu_\eph^{n\pm}) \uU_\eph^\text{num}.
\label{eq:good.numflux.non.cons}
\end{equation}
The computation of~\eqref{eq:good.numflux.non.cons} does not require any additional communication in comparison to the conservative case.
The expression~\eqref{eq:good.numflux.non.cons} has been used in all the numerical results of this work, and has given optimal order of accuracy for problems with smooth analytical solutions and results that are extremely close to those obtained by~\eqref{eq:bad.numflux.non.cons} for all the problems we tested.
It is possible that the expression~\eqref{eq:good.numflux.non.cons} might not work for some problems, but we have not encountered such a problem yet.
A more natural approach to approximate~\eqref{eq:bad.numflux.non.cons} without increasing the communication volume would be an interesting topic for future research.

\section{Subcell based scheme}\label{sec:blending.scheme}

By the formal high order accuracy of the RK scheme (Appendix~A of~\cite{babbar2025crk}), the cRKFR scheme~(\ref{eq:crkfr}, \ref{eq:crkfr.non.cons.final}) is a fully discrete high order accurate scheme.
Thus, as demonstrated in the numerical results (like Section~\ref{sec:mhd.convergence}), the proposed scheme gives the optimal order of accuracy for smooth solutions.
However, a direct usage of high order methods for problems with non-smooth solutions is known to produce Gibbs oscillations, as is also formalized in Godunov's order barrier theorem~{\cite{godunov1959}}.
Non-smooth solutions involving rarefaction waves, shocks and other discontinuities occur in a wide variety of practical applications involving hyperbolic equations.
Thus, a subcell based blending limiter has been developed for cRKFR methods for conservation laws~\cite{babbar2025crk}, based on~\cite{babbar2024admissibility,babbar2024multi,hennemann2021}.
The idea of the subcell based blending limiter is to take a convex combination of the high order cRKFR scheme and a lower order subcell scheme that is oscillation free and preserves the physical admissibility of the solution.
In this section, we discuss how to extend the subcell based scheme of~\cite{babbar2025crk} to solve hyperbolic equations with non-conservative products and stiff source terms, which are given in 1-D as
\begin{equation}\label{eq:1d.non.cons.source}
\uu_t + \pf (\uu)_x + \bB (\uu)  \uu_x = \bss.
\end{equation}
The description of the blending scheme begins by writing a single evolution step of the cRKFR scheme~(\ref{eq:crkfr},~\ref{eq:crkfr.non.cons.final}) as
\[
\vu^{H, n + 1}_e = \vu^n_e - \frac{\Delta t}{\Delta x_e}  \vR^H_e,
\]
where $\vu_e$ is the vector of nodal values in the element $e$.
Note that the $\vR^H_e$ consists of local implicit solves to handle the stiff source terms.
We use a lower order scheme written in the same form as
\[
\vu^{L, n + 1}_e = \vu^n_e - \frac{\Delta t}{\Delta x_e}  \vR^L_e,
\]
which will also consist of a local implicit solve for the stiff source terms.
Then the blended scheme is obtained by taking a convex combination of the two schemes as
\begin{equation}
  \vu^{n + 1}_e = (1 - \alpha_e)  \vu^{H, n + 1}_e + \alpha_e  \vu^{L, n +
  1}_e = \vu^n_e - \frac{\Delta t}{\Delta x_e}  [(1 - \alpha_e) \vR^H_e +
  \alpha_e  \vR^L_e], \label{eq:blended.scheme}
\end{equation}
where $\alpha_e \in [0, 1]$ is chosen based on a local smoothness indicator
taken from~{\cite{hennemann2021}} (also see Section 3.3
of~{\cite{babbar2024admissibility}} for a discussion on the indicator
of~{\cite{hennemann2021}}).
If $\alpha_e = 0$, the blending scheme is the high
order cRKFR scheme, while if $\alpha_e = 1$, the scheme becomes the low order
scheme that does not produce spurious oscillations and is admissibility
preserving (Definition~\ref{defn:admissibility.preserving}).
In subsequent sections, we explain the details of the lower order scheme, which is a first order finite volume method in this work.

\begin{figure}
\begin{center}
\includegraphics[width = 0.7\columnwidth]{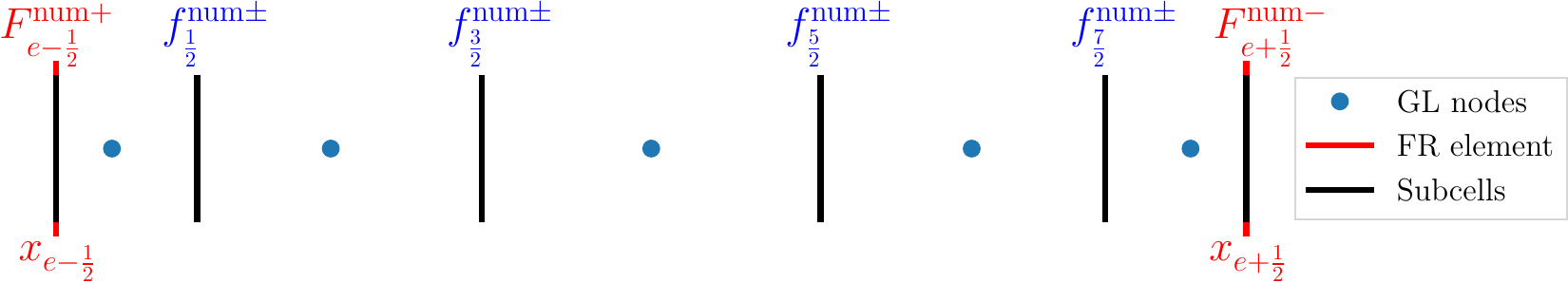} \\
(a) \\
\includegraphics[width = 0.7\columnwidth]{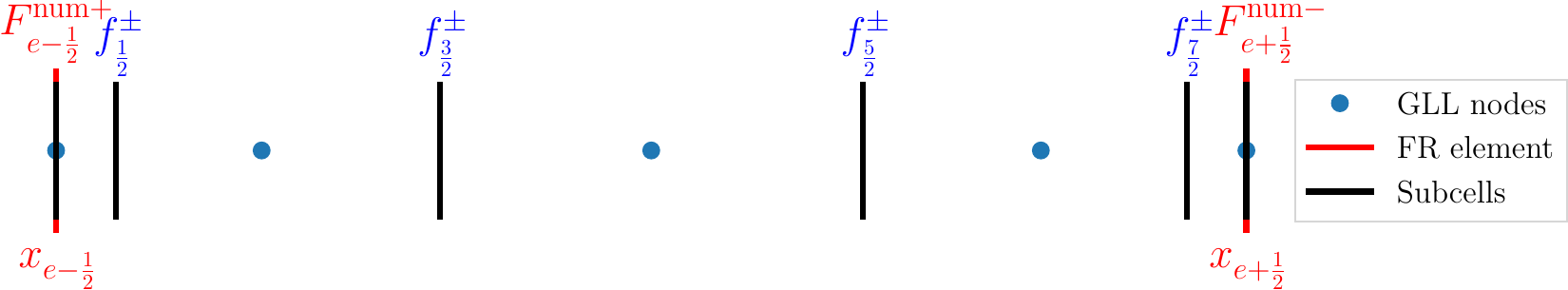} \\
(b)
\end{center}
\caption{Subcells used by lower order scheme for degree $N = 3$ using (a) Gauss-Legendre (GL) solution points (GL), (b) Gauss-Lobatto-Legendre (GLL) solution points.\label{fig:subcells}}
\end{figure}

Let us subdivide each element $\Omega_e$ into $N + 1$ subcells associated to
the solution points $\{ \xep, p = 0, 1, \ldots, N\}$ of the cRKFR scheme.
Thus, we will have $N + 2$ subfaces denoted by $\{x^e_{\pph}, p = - 1, 0,
\ldots, N\}$ with $x^e_{- \half} = x_{\emh}$ and $x^e_{\Nph} = x_\eph$.
For
maintaining a conservative scheme, the $p^{\text{th}}$ subcell is chosen so
that
\begin{equation}
  \label{eq:subcell.defn} x_{\pph}^e - x_{\pmh}^e = w_p \Delta x_e, \qquad 0
  \le p \le N,
\end{equation}
where $w_p$ is the $p^{\text{th}}$ quadrature weight associated with the
solution points.
Figure~\ref{fig:subcells} gives an illustration of the
subcells for the degree $N = 3$ case for Gauss-Legendre and
Gauss-Lobatto-Legendre solution points.
The low order scheme is obtained by updating the solution in each of the subcells by a finite volume scheme,

{\subequations{\begin{align}
  \uez^{L, n + 1} & = \uez^n - \frac{\Delta t}{w_0 \Delta x_e}
  [\fnumem_\half - \Fnump_{\emh}] + \Delta t \bss(\uu^{n+1}_{e,0}),  \label{eq:low.order.update.a}\\
  \uep^{L, n + 1} & = \uep^n - \frac{\Delta t}{w_p \Delta x_e}
  [\fnumem_\pph - \fnumep_\pmh]  + \Delta t \bss(\uu^{n+1}_{e,p}), \qquad 1 \le p \le N - 1, \\
  \ueN^{L, n + 1} & = \ueN^n - \frac{\Delta t}{w_N \Delta x_e}  [\Fnumm_\eph
  - \fnumep_\Nmh] + \Delta t \bss(\uu^{n+1}_{e,N}).
 \label{eq:low.order.update.b}
\end{align}
\label{eq:low.order.update}}}
The lower order numerical fluxes are chosen as
\begin{equation}
\fnumepm_{\pph} = \fnumpm(\uu_{e,p}, \uu_{e,p + 1}) \label{eq:numflux.args},
\end{equation}
where $\fnumpm$ is as defined in~\eqref{eq:combined.num.flux}.
% \begin{equation}
% \begin{gathered}
% \fnumpm(\uu_-, \uu_+) =  \half (\pf(\uu_-) + \pf(\uu_+) +  \bB(\uu_\pm) (\uu_- + \uu_+))
% - \lambda(\uu_-, \uu_+)  (\uu_+ - \uu_-), \\
% \lambda(\uu_-, \uu_+) = \max_{\uu \in \{\uu_-, \uu_+\}}\{ \sigma(\pf'(\uu) + \bB(\uu)) \}.
% \end{gathered} \label{eq:low.order.num.flux}
% \end{equation}
% The choice of $\lambda$ in~\eqref{eq:low.order.num.flux} is a natural extension of Rusanov's dissipation~\cite{rusanov1962} to non-conservative systems.
At the faces that are shared by both high and low order methods, the same inter-element flux $\Fnumpm_\eph$ is used in both the low and high order schemes.
This ensures that the conservative property is maintained when the contribution of non-conservative terms is zero.
To achieve high order accuracy, this flux needs to be high order accurate in smooth regions.
However, a high order flux may produce spurious oscillations near discontinuities when used in the low order scheme.
A natural choice to balance accuracy and oscillations is to take
\begin{equation}
  \Fnumpm_\eph = (1 - \alpha_\eph)  \FnumHOpm_\eph +
  \alpha_\eph  \fnumpm_\eph, \quad \alpha_\eph = \frac 12 (\alpha_e + \alpha_{e+1})\label{eq:Fblend},
\end{equation}
where $\FnumHOpm_\eph$ are the high order inter-element time-averaged numerical fluxes used in the cRKFR scheme~\eqref{eq:total.F.non.cons} and $\fnumpm_\eph = \fnumpm(\uu_{e,N}, \uu_{e+1,0})$~\eqref{eq:numflux.args} are lower order fluxes at the face $x_\eph$ shared between FR elements and subcells.
Since the blending coefficient $\alpha_\eph$ is based on a local smoothness indicator, it will bias the flux towards the lower order flux $\fnumpm_\eph$ near regions of lower solution smoothness.
However, to enforce admissibility preservation of the lower order scheme (Definition~\ref{defn:admissibility.preserving}), the inter-element fluxes have to be further limited, as we now explain.

\subsection{Admissibility preservation} \label{sec:flux.limiter}
In this section, we describe how we exploit the subcell based blending scheme to obtain admissibility preserving (Definition~\ref{defn:admissibility.preserving}) cRKFR IMEX schemes for the non-conservative equations~\eqref{eq:1d.non.cons.source}.
We begin by giving some notation and background on admissibility preservation strategies.
In this work, the admissibility constraints $\ad_j$~\eqref{eq:uad.form} need not be concave functions of the solution $\uu$.
We do assume that the admissible set $\Uad$ is convex, and the following weaker condition on the admissibility constraints for states $\uu_a, \uu_b$ (also mentioned in Section 5.2 of~\cite{babbar2024thesis}):
\[
\ad_j(\uu_a), \ad_j(\uu_b) > 0,\quad \forall j \le k\qquad \implies \qquad \ad_j(\theta \uu_a + (1 - \theta) \uu_b) > \epsilon_j(\uu_a, \uu_b),\quad \forall j \le k, \; \theta \in [0, 1].
\]
For the ten moment problem~\eqref{eq:tmp} whose admissibility constraint $\ad_3$ (the determinant of the pressure tensor) is not concave, $\epsilon_3(\uu_a, \uu_b) = \half \min(\ad_3(\uu_a), \ad_3(\uu_b))$ (see (2.9) of~\cite{meena2017}).
The admissibility preservation strategy is crucially reliant on finding, for $\uu_a$ and $\uu_b$ such that
\[
\ad_j(\uu_a), \ad_j(\uu_b) > 0, \quad \forall j < k, \qquad \ad_k(\uu_a) > 0, \ad_k(\uu_b) < 0,
\]
a $\theta \in [0, 1]$ such that
\begin{equation}
\ad_k(\theta \uu_a + (1 - \theta) \uu_b) \ge \epsilon_k(\uu_a),
\label{eq:admissibility.convex}
\end{equation}
where $\epsilon_k(\uu_a)$ is a positive number which is taken to be $\frac{1}{10} \ad_k(\uu_a)$ following~\cite{RuedaRamrez2021}.
A solution which is known to preserve high order accuracy~\cite{zhang2010c} is to find $\theta$ by solving the non-linear equation
\begin{equation}
\ad_k(\theta \uu_a + (1 - \theta) \uu_b) = \epsilon_k(\uu_a),
\label{eq:theta.eqn}
\end{equation}
which reduces to a cubic equation for the ten moment equations~\eqref{eq:tmp}~\cite{meena2017}.
In this work, when $\ad_k$ is concave, we use the simpler direct formula
\begin{equation} \label{eq:theta.formula}
\theta = \min \left( 1, \left| \frac{\epsilon_k(\uu_a) - \ad_k(\uu_a)}{\ad_k(\uu_b) - \ad_k(\uu_a)} \right| \right),
\end{equation}
which ensures~\eqref{eq:admissibility.convex} by the Jensen's inequality.
When $\ad_k$ is not concave, the formula~\eqref{eq:theta.formula} will not guarantee~\eqref{eq:admissibility.convex}, i.e., the convex combination $\theta \uu_a + (1 - \theta) \uu_b$ need not satisfy the admissibility constraint $\ad_k$ for $\theta$ computed using~\eqref{eq:theta.formula}.
Thus, in this case, we need to solve~\eqref{eq:theta.eqn} for $\theta$.
In the past works like~\cite{zhang2010c,meena2017}, the equation~\eqref{eq:theta.eqn} reduced to a polynomial equation and was thus solved using polynomial root finding methods.
In this work, we use the more general approach of using a Newton's method to solve~\eqref{eq:theta.eqn} for $\theta$.
For simplifying the description of the following algorithms, we introduce the notation
\begin{equation}
\theta_{\ad_k}(\epsilon_k, \uu_a, \uu_b) :=
\begin{cases}
\eqref{eq:theta.formula}, & \text{if } \ad_k \text{ is concave}, \\
\text{solution of~\eqref{eq:theta.eqn} using Newton's method}, & \text{otherwise}.
\end{cases}
\label{eq:theta.adk}
\end{equation}

We now describe the admissibility preservation strategy.
The first step is to ensure that the lower order blending scheme is admissibility-preserving at each solution point.
Since we assume that the first order numerical flux~\eqref{eq:combined.num.flux} gives an admissibility-preserving finite volume scheme, the only part that needs to be taken care of is the evolution using inter-element fluxes $\F_\eph^{\pm}$ where a convex combination of high and low order fluxes is used~\eqref{eq:Fblend}.
Thus, as a first step, we use a flux limiting procedure to find the inter-element fluxes $\F_\eph^\pm$ as described in Algorithm~\ref{alg:flux.limiting} to enforce admissibility preservation of the lower order method.

\begin{algorithm}
  \caption{Flux Limiting and evolution in subcells neighboring FR interface $x_\eph$} \label{alg:flux.limiting}
  \begin{algorithmic}
    \State $\Fnumpm_\eph \gets (1 - \alpha_\eph) \FnumHOpm_\eph
      + \alpha_\eph \fnumpm_\eph$ \Comment{Initial guess for the flux}

    \State $\tilde{\uu}_0^{n + 1} \gets
      \uu_{e + 1, 0}^n - \dfrac{\Delta t}{w_0 \Delta x_{e + 1}}
        ( \fnumepo_\half - \Fnump_\eph )$ \Comment{Evolutions whose admissibility is to be enforced}
    \State $\tilde{\uu}_N^{n + 1} \gets
      \uu_{e, N}^n
      - \dfrac{\Delta t}{w_N \Delta x_e}
        ( \Fnumm_\eph - \fnumep_{N - \half} )$ \Comment{Evolutions whose admissibility is to be enforced}

    \State $\tilde{\uu}_0^{\text{low}, n + 1} \gets
      \uu_{e + 1, 0}^n - \dfrac{\Delta t}{w_0 \Delta x_{e + 1}}
        ( \fnumepo_\half - \fnump_\eph )$ \Comment{Low order admissible evolutions}
    \State $\tilde{\uu}_N^{\text{low}, n + 1} \gets
      \uu_{e, N}^n
      - \dfrac{\Delta t}{w_N \Delta x_e}
        ( \fnumm_\eph - \fnumep_{N - \half} )$ \Comment{Low order admissible evolutions}

    \For{$k = 1$ to $K$}
      \State $\epsilon_0, \epsilon_N \gets
        \tfrac{1}{10} \ad_k(\tilde{\uu}_0^{\text{low}, n+1}),\;
        \tfrac{1}{10} \ad_k(\tilde{\uu}_N^{\text{low}, n+1})$
      \Comment{As in~\eqref{eq:admissibility.convex}, the factor $\frac{1}{10}$ is chosen following~\cite{RuedaRamrez2021}.}

      \State $\displaystyle
        \theta \gets \min_{p = 0, N} \left \{
          \theta_{\ad_k}\left(
            \epsilon_p,
            \tilde{\uu}_p^{\text{low}, n + 1},
            \tilde{\uu}_p^{n + 1}
          \right)\right \}$
          % \left|
          %   \frac{\epsilon_p - \ad_k(\tilde{\uu}_p^{\text{low}, n + 1})}{
          %         \ad_k(\tilde{\uu}_p^{n + 1})
          %       - \ad_k(\tilde{\uu}_p^{\text{low}, n + 1})}
          % \right|, 1
        \Comment{The $\theta_{\ad_k}$ is defined in~\eqref{eq:theta.adk}.}

      \State $\Fnumpm_\eph \gets
        \theta \Fnumpm_\eph + (1 - \theta) \fnumpm_\eph$

      \State $\tilde{\uu}_0^{n + 1} \gets
        \uu_{e + 1, 0}^n
        - \dfrac{\Delta t}{w_0 \Delta x_{e + 1}}
          ( \fnumepo_{\half} - \Fnump_\eph )$

      \State $\tilde{\uu}_N^{n + 1} \gets
        \uu_{e, N}^n - \dfrac{\Delta t}{w_N \Delta x_e} ( \Fnumm_\eph - \fnumep_{N - \half} )$
    \EndFor

    \State $\uu_0^{n + 1} \gets
      \tilde{\uu}_0^{n + 1} + \Delta t \bss(\uu_0^{n + 1})$     \Comment{Add source term contribution implicitly}
    \State ${\uu}_N^{n + 1} \gets
      \tilde{\uu}_N^{n + 1} + \Delta t \bss(\uu_N^{n + 1})$     \Comment{Add source term contribution implicitly}
  \end{algorithmic}
\end{algorithm}
In~{\cite{babbar2025crk}}, the conservative analogue of Algorithm~\ref{alg:flux.limiting} was enough to obtain an admissibility-preserving in means scheme for conservation laws (Definition~\ref{defn:admissibility.preserving.means}), following which the scaling limiter of~{\cite{zhang2010c}} was used in~\cite{babbar2025crk} to obtain an admissibility-preserving scheme.
However, for general equations of the form~\eqref{eq:1d.non.cons.source}, Algorithm~\ref{alg:flux.limiting} does not ensure admissibility preservation in means.
This is because the non-conservative products also contribute to the evolution of the cell average, and not just the inter-element fluxes.
Thus, we cannot use the scaling limiter of~\cite{zhang2010c} at this stage to obtain an admissibility-preserving scheme.
However, an alternative to the scaling limiter is to further limit using the blending coefficient $\alpha_e$~\eqref{eq:blended.scheme} so that the final solution in the element is admissible.
This was already used for conservative equations in~\cite{RuedaRamrez2021}, and mentioned as a possible alternative in Section 5 of~\cite{babbar2024admissibility}.
The same approach works for the non-conservative equations.
The rough idea is to find, for each element $e$, the blending coefficient $\alpha_e \in [0, 1]$ so that
\[
\vu^{n + 1}_e = (1 - \alpha_e)  \vu^{H, n + 1}_e + \alpha_e  \vu^{L, n + 1}_e \in \Uad .
\]
This is performed by Algorithm~\ref{alg:alpha.limiting}.
At the end of the procedure, the solution is admissible at each solution point in each element.
\begin{algorithm}
\caption{Final admissibility enforcement limiting} \label{alg:alpha.limiting}
\begin{algorithmic}
\State $\vu^{n + 1}_e \gets (1 - \alpha_e)  \vu^{H, n + 1}_e + \alpha_e  \vu^{L, n + 1}_e$. \Comment{Initial limiting from smoothness indicator}
\For{$e$ in \texttt{eachelement(grid)}}
\For{$k = 1$ to $K$}
\State $\displaystyle \epsilon \gets \tfrac{1}{10} \min_{0 \le p \le N} \ad_k(\uu_{e,p}^{L, n + 1})$
\Comment{As in~\eqref{eq:admissibility.convex}, the factor $\frac{1}{10}$ is chosen following~\cite{RuedaRamrez2021}.}
\State $\displaystyle
\theta \gets
% \min_{0 \le p \le N} \left|
% \frac{\epsilon - \ad_k(\uu_{e,p}^{L, n + 1})}{
% \ad_k(\uu_{e,p}^{n + 1})
% - \ad_k(\uu_{e,p}^{L, n + 1})}
% \right|$
\min_{0 \le p \le N} \left \{
\theta_{\ad_k}\left(
\epsilon,
\uu_{e,p}^{L, n + 1},
\uu_{e,p}^{n + 1}
\right)\right \}$
\Comment{The $\theta_{\ad_k}$ is defined in~\eqref{eq:theta.adk}.}

\State $\vu_{e}^{n + 1} \gets \theta \vu_{e}^{n + 1} + (1 - \theta) \vu_{e}^{L, n + 1}$
\EndFor
\EndFor
\end{algorithmic}
\end{algorithm}

\section{Numerical results} \label{sec:numerical.results}

We now show various numerical results to demonstrate the performance of the proposed cRKFR scheme for hyperbolic equations with non-conservative products and stiff source terms.
Here we summarize the test cases that are done in this section.
In Section~\ref{sec:scalar.tests}, we test the scheme using the Jin-Xin relaxation system~\cite{jinxin1995} which is a standard test for schemes that can handle stiff source terms, three Burgers' (and Burgers-type) equations tests with stiff source terms from~\cite{svard2011} and non-conservative variable advection equations.
In Section~\ref{sec:reactive.euler}, we test the scheme for reactive Euler equations with stiff source terms from~\cite{svard2011} and a non-stiff test with flow over a step to show our scheme's capability to run physically relevant problems.
Section~\ref{sec:ten.moment} contains a test of the ten moment equations with two rarefactions and a very stiff Gaussian type source term, leading to a near-vacuum state.
A shear shallow water (SSW) equation model is tested in Section~\ref{sec:ssw} against physical reference data to validate the treatment of non-conservative products.
Magnetohydrodynamics and multi-ion magnetohydrodynamics equations are tested in Sections~\ref{sec:mhd},~\ref{sec:multi.ion.mhd} for 2-D problems to validate the scheme with some standard test cases.
Convergence tests are also done for MHD and multi-ion MHD equations to show the benefit of using GL points over GLL points.
% A stiff source term collision test is also done for multi-ion MHD.

The most general form of hyperbolic equations with non-conservative products and stiff source terms that we consider in this section are the 2-D equations of the form
\begin{equation}
\uu_t + \pf (\uu)_x + \pg (\uu)_y + \bB_1(\uu)  \uu_x + \bB_2(\uu)  \uu_y = \bss(t, \uu).
\end{equation}
The time step size is computed as
\begin{align} \label{eq:time.step.2d}
\Delta t = C_s \min_e \left( \frac{\sigma(\pf'(\au_e) + \bB_1(\au_e))}{\Delta x_e} + \frac{\sigma(\pg'(\au_e) + \bB_2(\au_e))}{\Delta y_e} \right)^{-1} \text{CFL}(N),
\end{align}
where $e$ is the element index, $\avg{\uu}_e$ is the element mean (Definition~\ref{defn:admissibility.preserving.means}), $\Delta x_e, \Delta y_e$ are the element lengths in $x,y$ directions respectively,  $C_s \le 1$ is a safety factor, $\text{CFL}(N)$ is the optimal CFL number obtained by Fourier stability analysis.
We use the optimal CFL numbers of the LWFR scheme from Table 1 of~\cite{babbar2022lax} since the cRKFR scheme is linearly equivalent to the LWFR scheme~\cite{chen2024,babbar2025crk}.
These CFL numbers are also in Table 1.1 of~\cite{chen2024}.
The safety factor is $C_s = 0.9$ unless specified otherwise.
For the 1-D equations,~\eqref{eq:time.step.2d} is used by setting $\fg = \bzero, \bB_2 = \bzero$.
The developments have been contributed to the Julia~\cite{bezanson2017julia} package \texttt{Tenkai.jl}~\cite{tenkai} and the run files for generating the results in this paper are available at~\cite{babbar2025crknonconsrepro}.

\subsection{Scalar equations}\label{sec:scalar.tests}

We begin by testing with scalar equations involving stiff source terms.
In these tests, we find that the choice of the IMEX scheme is also relevant for handling extremely stiff source terms.
The simple second order IMEX scheme HT (1,1,2) of~\cite{Hairer1996}~\eqref{eq:ht112.butcher} is found to be insufficient for extremely stiff source terms, while the L-stable third order SSP3-IMEX(4,3,3) scheme of~\cite{pareschi2005}~\eqref{eq:ssp33.butcher} works for all problems.

\subsubsection{Jin-Xin relaxation system} \label{sec:jin.xin}

The Jin-Xin relaxation system~\cite{jinxin1995} is a simple and popular example of a hyperbolic system with relaxation source terms.
The idea is that, given a scalar conservation law
\begin{equation}
  u_t + f(u)_x = 0, \label{eq:scalar.cons.law}
\end{equation}
we introduce a new variable $v$ and obtain the $2 \times 2$ system
\begin{align}
\partial_t u + \partial_x v &= 0, \label{eq:jinxin1}\\
\partial_t v + a^2 \partial_x u &= -\frac{1}{\varepsilon} \left( v - f(u) \right),
\label{eq:jinxin2}
\end{align}
where $a > 0$ defines the characteristic speed of the system.
% It is easy to see that, as $\varepsilon \to 0$, the system converges to the scalar conservation law.
% In order to determine conditions on $a$, we do the \emph{Chapman–Enskog expansion} following~\cite{boscarino2024}.
% We can observe from~\eqref{eq:jinxin2} that
% \begin{equation}
%     v = f(u) - \varepsilon \left( v_t + a^2 u_x \right) = f(u) + \mathcal{O}(\varepsilon),
%     \label{eq:v_expansion}
% \end{equation}
% from which it follows that
% \begin{align*}
%     v_t &= f(u)_t + \mathcal{O}(\varepsilon)
%     = f'(u) u_t + \mathcal{O}(\varepsilon)
%     = - f'(u) v_x + \mathcal{O}(\varepsilon) \\
%         &= - f'(u) f(u)_x + \mathcal{O}(\varepsilon)
%     = - f'(u)^2 u_x + \mathcal{O}(\varepsilon).
% \end{align*}

\begin{figure}
\centering
\begin{tabular}{cccc}
\includegraphics[width=0.2225\columnwidth]{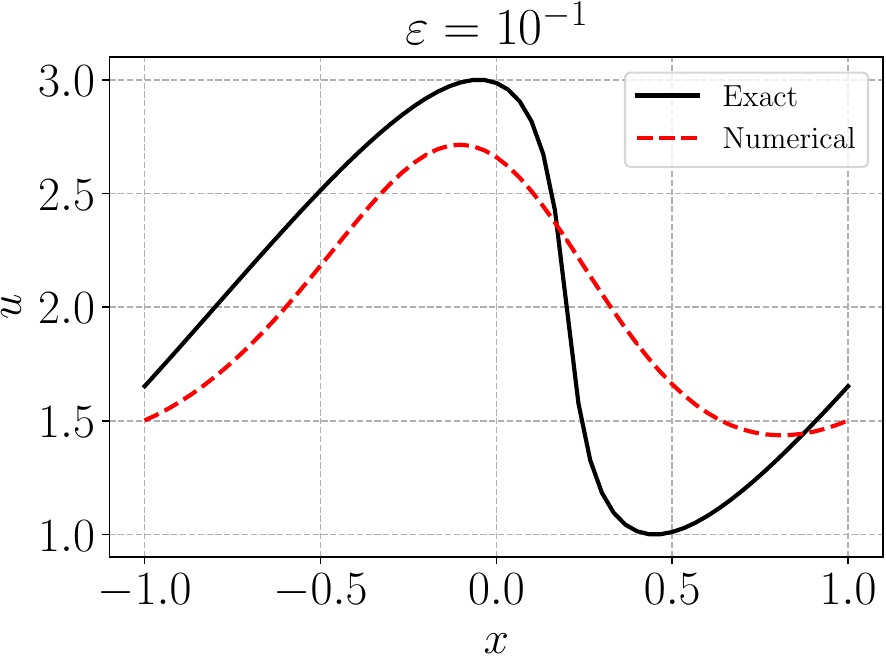} &
\includegraphics[width=0.2225\columnwidth]{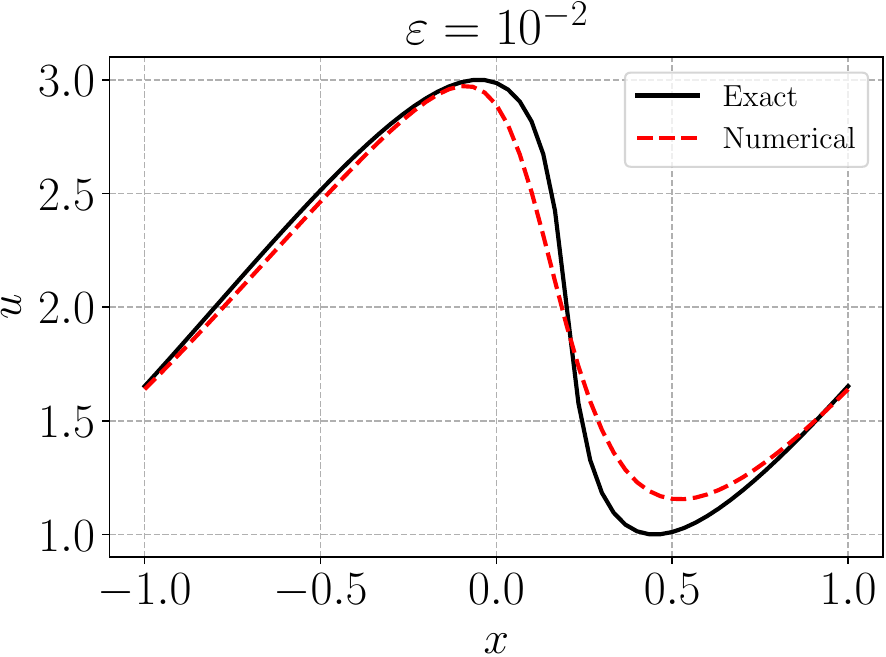} &
\includegraphics[width=0.2225\columnwidth]{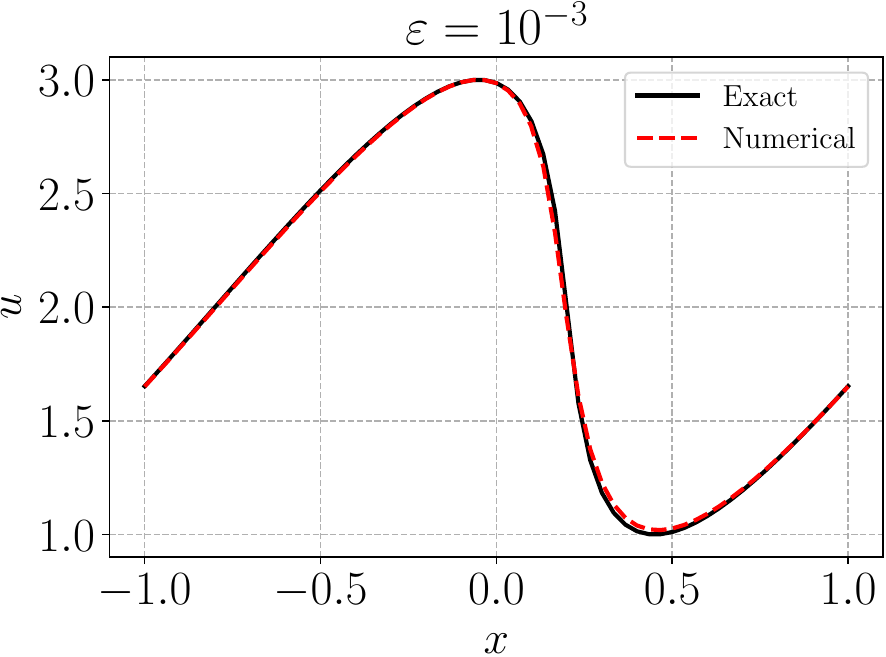} &
\includegraphics[width=0.2225\columnwidth]{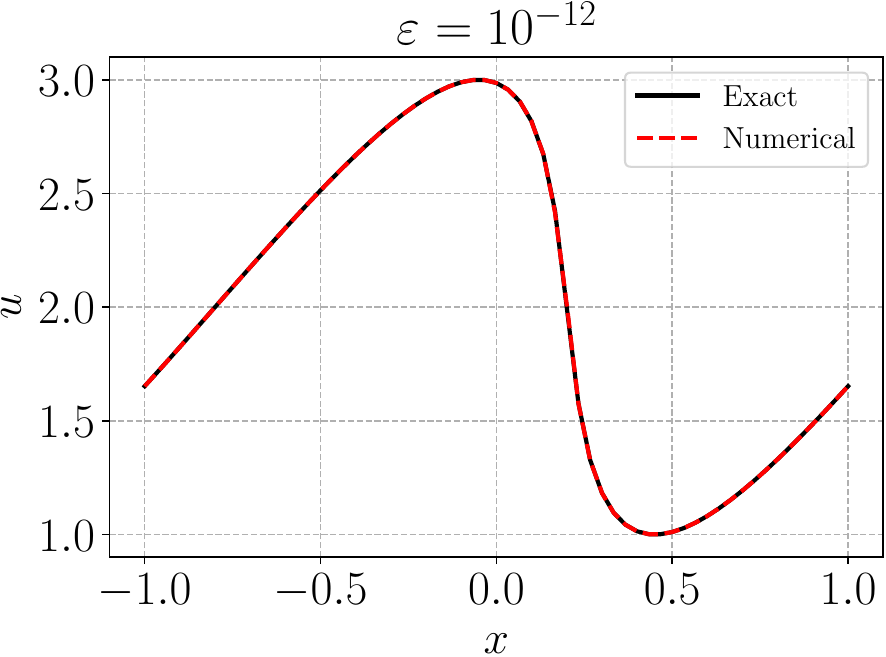} \\
(a) & (b) & (c) & (d)
\end{tabular}
\caption{Jin-Xin relaxation system for Burgers' equation with $\varepsilon = 10^{-1}$ (a), $\varepsilon = 10^{-2}$ (b), $\varepsilon = 10^{-3}$ (c), $\varepsilon = 10^{-12}$ (d) using initial condition~\eqref{eq:jin.xin.ic} at time $t=0.25$.}\label{fig:jin.xin}
\end{figure}

% Using this relation in~\eqref{eq:v_expansion}, and substituting the result into the first equation of the system~\eqref{eq:jinxin1}, we obtain
% \begin{equation}
% u_t + f(u)_x = \varepsilon \left[ \left( a^2 - f'(u)^2 \right) u_x \right]_x + \mathcal{O}(\varepsilon^2).
% \label{eq:conv_diff}
% \end{equation}
% Neglecting second-order terms in the small parameter \(\varepsilon\), the above expression represents a nonlinear convection–diffusion equation, in which the diffusion coefficient is
% \begin{equation}
% \nu = \varepsilon \left( a^2 - f'(u)^2 \right). \label{eq:nu}
% \end{equation}
For small values of \(\varepsilon\), by performing a \emph{Chapman-Enskog expansion}, it is shown in~\cite{boscarino2024} that the well-posedness of~\eqref{eq:jinxin2} requires the \emph{subcharacteristic condition}
\begin{equation}
|f'(u)| \le a.
\label{eq:subchar}
\end{equation}
We test the Jin-Xin relaxation system for the Burgers' flux function $f(u) = \frac{u^2}{2}$, with the initial condition
\begin{equation} \label{eq:jin.xin.ic}
u(x,0) = 2 + \sin(\pi (x - 0.7)), \qquad x \in [-1, 1].
\end{equation}
The subcharacteristic condition~\eqref{eq:subchar} is satisfied by choosing $a = 3$.
We run the simulation until time $t = 0.25$ with periodic boundary conditions and show the results for different values of $\varepsilon$ in Figure~\ref{fig:jin.xin}.
The SSP3-IMEX(4,3,3) scheme of~\cite{pareschi2005}~\eqref{eq:ssp33.butcher} is used.
The simplest IMEX scheme HT (1,1,2) of~\cite{Hairer1996}~\eqref{eq:ht112.butcher} could not be used for tiny $\varepsilon$ values.
This is likely because the HT (1,1,2) scheme is only A-stable, while the SSP3-IMEX(4,3,3) scheme is L-stable (Appendix~\ref{app:imex.schemes}).
By using the SSP3-IMEX(4,3,3) scheme, we were able to run tests with $\varepsilon$ as small as $10^{-12}$ without any issues.
This shows that our IMEX cRKFR scheme can handle extremely stiff source terms with the right choice of IMEX scheme.
\subsubsection{Burgers-type equations with stiff source terms}
\begin{figure}
\centering
\begin{tabular}{ccc}
\includegraphics[width=0.3\columnwidth]{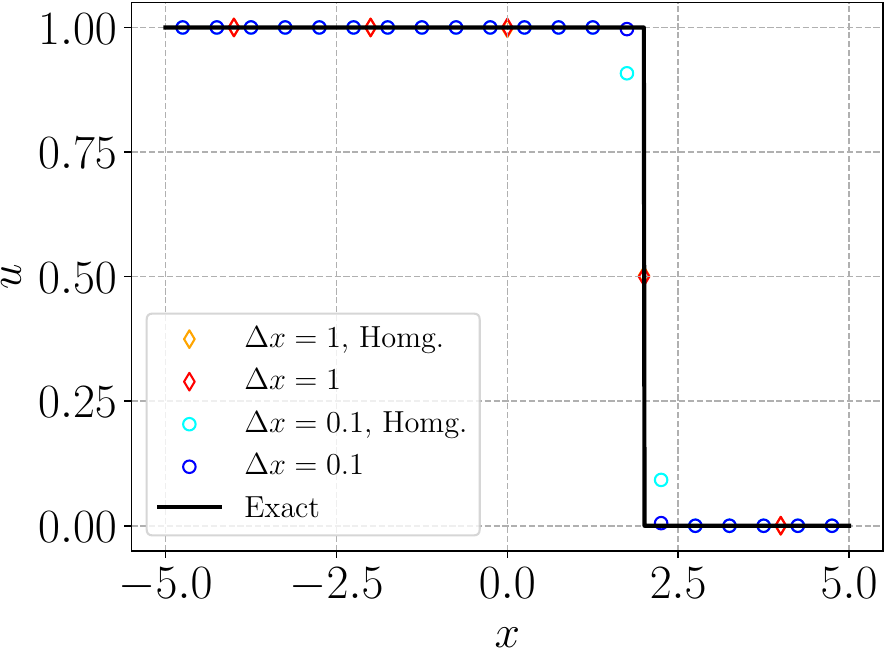} &
\includegraphics[width=0.3\columnwidth]{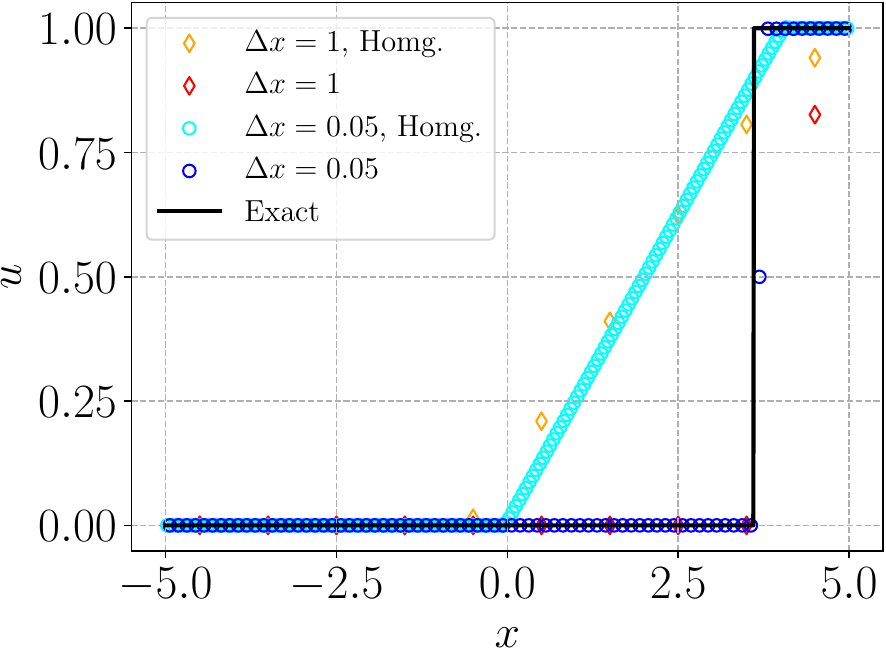} &
\includegraphics[width=0.3\columnwidth]{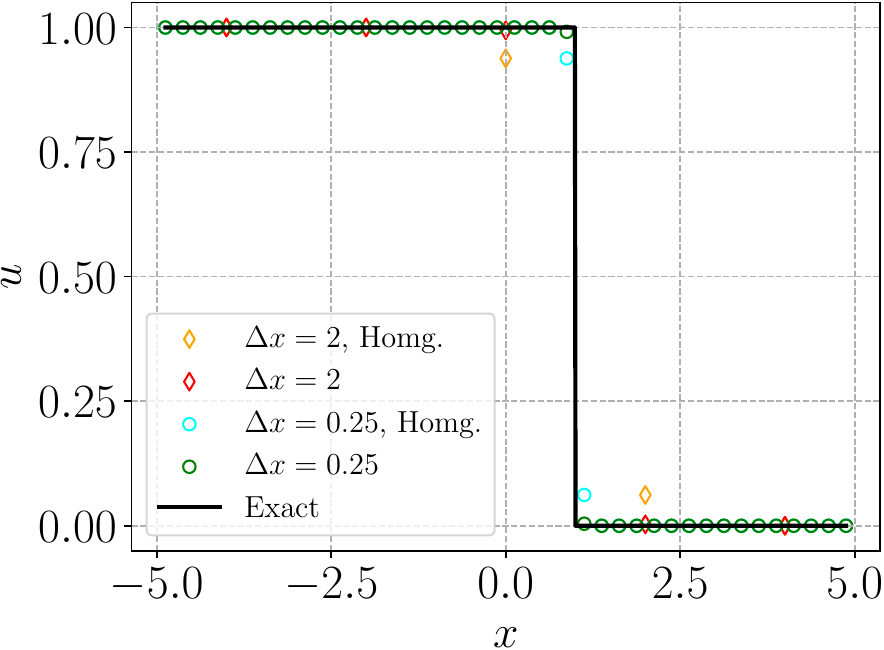} \\
(a) & (b) & (c)
\end{tabular}
\caption{Scalar conservation laws~\eqref{eq:scalar.cons.law} with shocks and stiff source term~\eqref{eq:svard.source} with $\nu = 10^6$, simulated using SSP3-IMEX(4,3,3)~\cite{pareschi2005} and polynomial degree $N=3$
(a) $f(u) = u^2 / 2$ (Burgers' equation) with
$\beta = 0.9$ at $t = 4$ with initial condition~\eqref{eq:ic.shock},
(b) $f(u) = u^2 / 2$ with $\beta = 0.9$ at $t = 4$ with initial condition~\eqref{eq:ic.rare},
(c) $f(u) = u^4 / 4$, $\beta = 0.4$ at $t=4$ with initial condition~\eqref{eq:ic.shock}.}
\label{fig:svard}
\end{figure}

The next scalar tests are from~\cite{leveque1990,svard2011} and mentioned to be good benchmarks for numerical methods that can handle stiff source terms in~\cite{svard2011}.
For the source term
\begin{equation}
s_{\nu, \beta}(u) = \nu (1 - u) (u - \beta) u,
\label{eq:svard.source}
\end{equation}
we look at the scalar conservation law~\eqref{eq:scalar.cons.law} with fluxes $f(u) = \frac{u^2}{2}$ (Burgers' equation) and $f(u) = \frac{u^4}{4}$.
The source term~\eqref{eq:svard.source} was mentioned in~\cite{svard2011} as a simple toy problem to test numerical methods for the reactive Euler equations~\eqref{eq:reactive.euler}.
The $\nu$ in~\eqref{eq:svard.source} is a large number that prescribes the stiffness of the source term.
The timescale of the source term is $\mathcal{O}(1/\nu)$.
We are able to run tests with all values of $\nu$ using the SSP3-IMEX(4,3,3) scheme.
The HT (1,1,2) scheme could also run these tests, but produced more oscillations which could be controlled by using the smoothness indicator of~\cite{basak2025bound} for the blending limiter~\eqref{eq:blended.scheme} in place of the indicator of~\cite{hennemann2021} which is used in all the numerical results shown in this work.
However, for some $\nu$ in $[10^3, 10^6)$, we needed to decrease the step sizes in the iterative solver used in solving the implicit equations for the source term.
The results are shown for $\nu = 10^6$, but the reproducibility repository~\cite{babbar2025crknonconsrepro} contains files for simulating with other values of $\nu$ including some that required smaller step sizes in the iterative solver.
The $\beta$ in~\eqref{eq:svard.source} is an unstable root of the source term.
If we look at an ODE with the source term~\eqref{eq:svard.source}, then any initial data $u_0 > \beta$ converges to $1$, while $u_0 < \beta$ goes to $0$.
Thus, choosing a good initial guess for the implicit solve is important to get the correct solution.
We use the strategy mentioned in~\cite{svard2011} which uses the solution of the homogeneous equation to choose the initial guess for the implicit solve.
Thus, the scheme requires solving the homogeneous equation first to get the solution $v_{e,p}$ at each solution point, and then the initial guess at the solution point is chosen to be
\[
v_{e,p}^* = \begin{cases}
1, \quad & \text{if } v_{e,p} > \beta,\\
0, \quad & \text{if } v_{e,p} < \beta.
\end{cases}
\]
Additionally, we found it necessary in some cases (e.g., on a much finer mesh) to compute the homogeneous solution with the purely first order scheme on subcells (i.e., setting $\alpha = 1$ in~\eqref{eq:blended.scheme} for the homogeneous solver) because even a few spurious solutions in the initial guess caused us to converge to the wrong solution.
This is in fact advantageous as the first order scheme is cheaper to compute.
Despite that, having two solvers does increase the computational cost.
This increase in computations is mentioned to be worth the cost in~\cite{svard2011} as it allows the solution of the stiff problem with an advective time step restriction~\eqref{eq:time.step.2d} instead of the source term time step restriction of $\mathcal{O}(1/\nu)$ (Appendix~\ref{app:stiff.source}).
The results are shown in Figure~\ref{fig:svard} where Figures~\ref{fig:svard}a,c use the initial condition
\begin{equation}
u(x,0) =
\begin{cases}
1, \quad &  x \ge 0,\\
0, & x < 0,
\end{cases}
\label{eq:ic.shock}
\end{equation}
with fluxes $f(u) = \half u^2$ (Burgers' equation) and $f(u) = \frac{1}{4} u^4$ respectively.
By Rankine-Hugoniot conditions, the solution is a shock traveling at speeds $0.5, 0.25$ respectively, and the source term does not activate for the exact solution.
However, as can be seen from the solution of the homogeneous solution in Figure~\ref{fig:svard}, numerical dissipation causes the solution to reach values other than $0,1$ leading to an activation of the source term.
Thus, the fact that our IMEX scheme is able to capture the correct shock speeds shows its capability to handle stiff source terms.
In Figure~\ref{fig:svard}b, we show the solution for Burgers' equation with initial condition where states from~\eqref{eq:ic.shock} are swapped
\begin{equation}
u(x,0) =
\begin{cases}
0, \quad &  x \ge 0,\\
1, & x < 0.
\end{cases}
\label{eq:ic.rare}
\end{equation}
The initial condition~\eqref{eq:ic.rare} results in a rarefaction wave in the absence of source terms, which increases the difficulty in our case because the solution of the non-homogeneous equation is now different from the homogeneous equation, where the former is a shock travelling at the speed $\beta$~\cite{svard2011} while the latter is a rarefaction wave.
In all tests in Figure~\ref{fig:svard} we show results using coarse and fine meshes always using the hyperbolic CFL condition~\eqref{eq:time.step.2d}, showing that the scheme can handle stiff source terms with the large hyperbolic time step sizes.

\subsubsection{Non-conservative variable coefficient advection equation}

\begin{figure}
\centering
\begin{tabular}{cc}
\includegraphics[width=0.43\columnwidth]{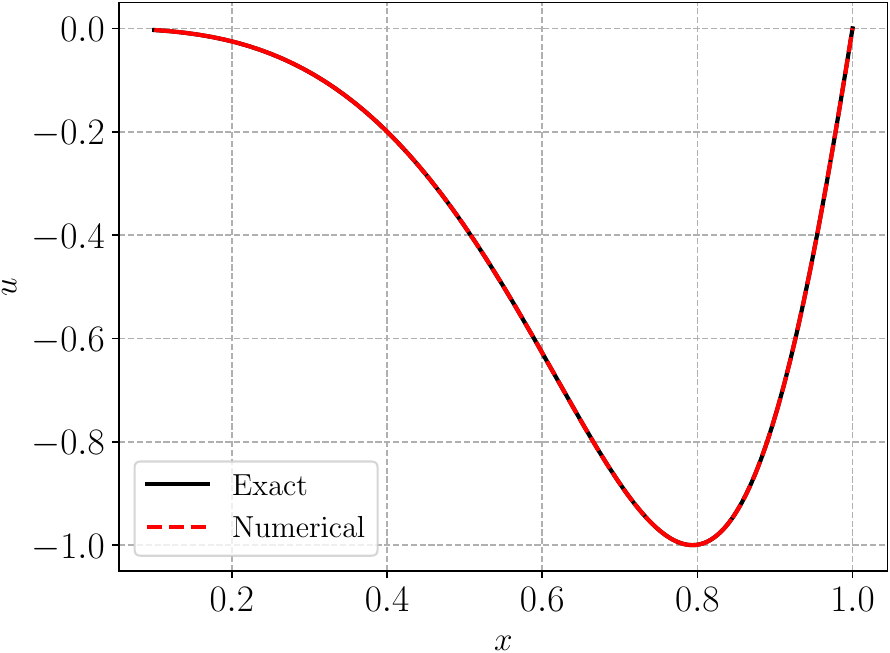} &
\includegraphics[width=0.35\columnwidth]{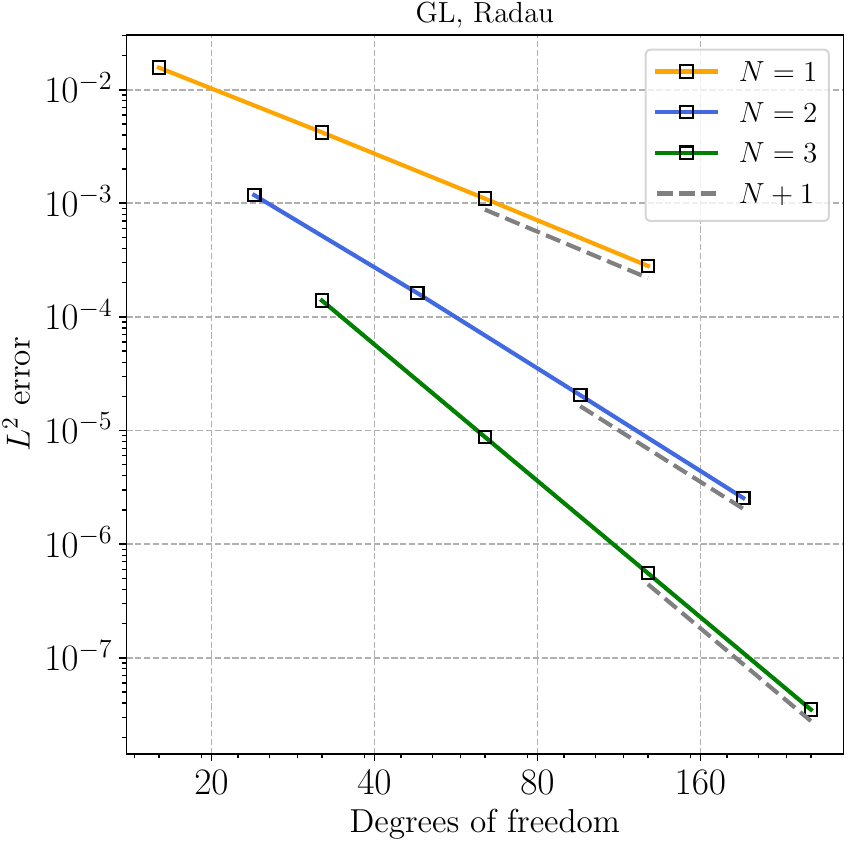} \\
(a) & (b)
\end{tabular}
\caption{Numerical solution of variable coefficient equation~\eqref{eq:var.coeff.adv} with $a(x) = \frac{1}{x^2}$ showing (a) Numerical solution at time $t=1$ with initial condition $u_0(x) = \sin(\pi x)$ using $8$ elements and polynomial degree $N=3$, (b) convergence study showing $L^2$ errors and orders of accuracy for polynomial degrees $N=1,2,3$. }
\label{fig:var.coeff.adv}
\end{figure}

Variable coefficient equations have already been used for numerical tests in single stage flux reconstruction schemes in~\cite{babbar2022lax,babbar2025crk,babbar2024multi}.
However, they were limited to the equations that could be written in the conservative form $u_t + f(u, x)_x$.
By extending the scheme to non-conservative equations~\eqref{eq:1d.non.cons.source}, we are now able to also solve the equations of the form
\begin{equation}
\begin{aligned}
\partial_t v + a (x) \partial_x v &= 0,\\
v (x, 0) &= v_0 (x) .
\end{aligned}
\label{eq:var.coeff.adv}
\end{equation}
It is possible (and computationally efficient) to extend the previously described framework to variable coefficient equations by replacing $\bB (\uu)$ in~\eqref{eq:1d.non.cons.source} with $\bB (\uu, x)$.
However, for simplicity, we instead cast the variable coefficient equation~\eqref{eq:var.coeff.adv} in form of a non-conservative system~\eqref{eq:pure.non.conservative}.
To do that, we define $\uu = (u_1, u_2) := (v, a)$ so that the above system is equivalent to
\begin{equation}
\begin{aligned}
\partial_t u_1 + u_2 \partial_x u_1 & = 0,\\
\partial_t u_2 & = 0,\\
\uu (x, 0) & = (v_0 (x), x) .
\end{aligned}
\label{eq:var.coeff.adv.non.cons}
\end{equation}
Our experience with numerical experiments shows that we are able to solve the system~\eqref{eq:var.coeff.adv.non.cons} as long as it is strictly hyperbolic, i.e., the matrix
\[ A = \left[\begin{array}{cc}
     u_2 & 0\\
     0 & 0
   \end{array}\right]
\]
has distinct real eigenvalues, which is true as long as $u_2 = a(x) \ne 0$.
We test the scheme with $a(x) = x^{-2}$ in domain $x \in [0.1,1]$ with initial condition $u_0(x) = \sin (\pi x)$.
The general exact solution of the equation is obtained using the method of characteristics to be $u(x,t) = u_0\bigl( (x^3 - 3t)^{1/3} \bigr)$.
Thus, the numerical solution is shown in Figure~\ref{fig:var.coeff.adv}a where good agreement with the exact solution is observed.
In order to validate the scheme in Figure~\ref{fig:var.coeff.adv}b, we perform a convergence study where we see optimal order of accuracy for polynomial degrees $N=1,2,3$.

\subsection{Reactive Euler's equations} \label{sec:reactive.euler}
\begin{figure}
\centering
\begin{tabular}{cc}
\includegraphics[width=0.46\columnwidth]{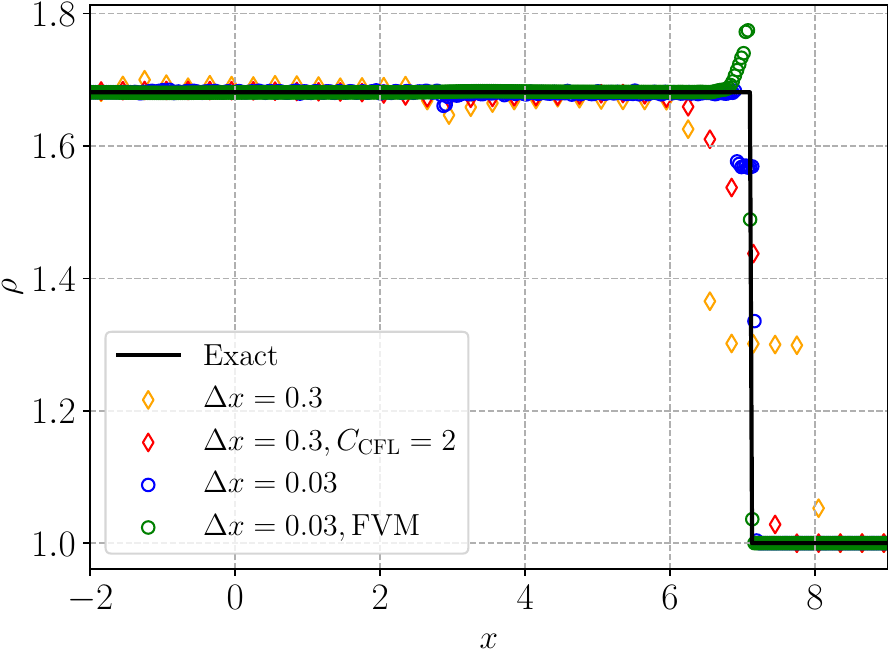} &
\includegraphics[width=0.46\columnwidth]{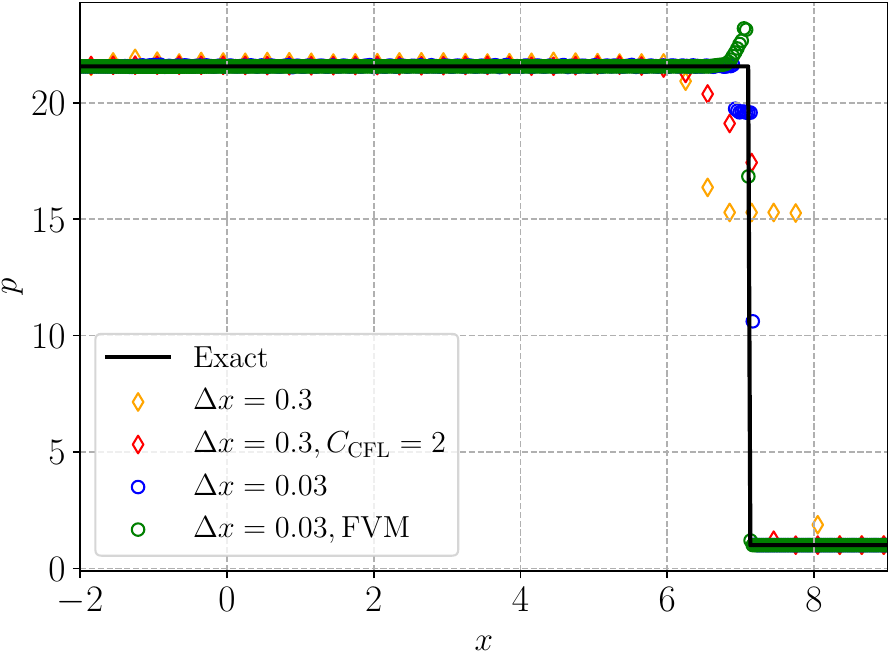} \\
(a) $\rho$ & (b) $p$
\end{tabular}
\caption{Riemann problem with reactive Euler equations~\eqref{eq:reactive.euler} with stiff source term simulated using SSP3-IMEX(4,3,3)~\cite{pareschi2005} and polynomial degree $N=3$ at time $t=1$.}
\label{fig:reactive.euler.rp1}
\end{figure}

We now test the scheme for the reactive Euler equations where the system~\eqref{eq:con.law.source} is specified by
\begin{equation}
\begin{gathered}
\uu = (\rho, \rho v, E, \rho Y),
\quad \pf = ( \rho v, p + \rho v^2, (E + p) v, \rho v Y )^T,
\quad \bss = (0, 0, 0, Q), \\
E = E (\rho, v, p) = \frac{p}{\gamma - 1} + \half \rho v^2 + q_0 \rho Y, \\
Q =  - K(T) \rho Y, \quad K(T) = A \exp\left(-\frac{T_A}{T}\right), \quad T = \frac{p}{\rho}, \quad A \gg 1.
\end{gathered} \label{eq:reactive.euler}
\end{equation}
The variables $(\rho, v, E, Y, p, T)$ are the density, velocity, total energy, mass fraction of the chemical reactant, pressure and temperature respectively.
The constant $\gamma$ is the ratio of specific heats, $q_0$ is the heat release energy, and $K(T)$ is the chemical reaction rate where $A, T_A$ are experimentally determined constants.
The chemical reaction rate $K$ is typically large which leads to stiffness in the source term $Q$, similar to the linear case (Appendix~\ref{app:stiff.source}).
It was mentioned in~\cite{svard2011} that the implicit equation that arises from the source term~\eqref{eq:reactive.euler} requires special treatment.
We describe the special treatment here for obtaining the solution at stage $i$ of the IMEX scheme~\eqref{eq:crkfr.inner}.
Defining
\[
\tilde{\uu}^{(i-1)} = \uu^n - \Delta t \sum_{j = 1}^{i - 1} \tilde{a}_{ij} \partial_x^{\text{loc}}  \pf_h ( \uu^{(j)} ) + \Delta t \sum_{j = 1}^{i - 1} a_{ij}  \bss ( t_n + c_j \Delta t, \uu^{(j)} ),
\]
% \uu^{(i)} & = & \uu^n - \Delta t \sum_{j = 1}^{i - 1} \tilde{a}_{i
%  j} \partial_x^{\text{loc}}  \pf_h ( \uu^{(j)} ) +
%   \Delta t \sum_{j = 1}^i a_{ij}  \bss ( t_n +
%   c_j \Delta t, \uu^{(j)} )
the implicit equation that we solve at each point $(e,p)$~\eqref{eq:reactive.euler.implicit} is replaced by the following equation inspired from~\cite{svard2011}
\begin{equation}
(\rho Y)_{e,p}^{(i)} + a_{ii} \Delta t K(T_{e,p}^{(i-1)}) (\rho Y)_{e, p}^{(i)} = \tilde{\uu}(\rho Y)^{(i-1)}_{e,p},
\label{eq:reactive.euler.implicit}
\end{equation}
where $\tilde{\uu}^{(i-1)}(\rho Y)$ denotes the component of $\tilde{\uu}^{(i-1)}$ corresponding to the conservative variable $\rho Y$.
Note that the difference between the above equation~\eqref{eq:reactive.euler.implicit} and the general implicit equation~\eqref{eq:crkfr.inner} is that the temperature term in $K(T)$ is taken from the previous stage $(i-1)$ instead of the current stage $i$.
This actually leads to a simple scalar linear equation for $(\rho Y)_{e,p}^{(i)}$ at each solution point $(e,p)$.
The first test with~\eqref{eq:reactive.euler} is a Riemann problem with stiff source term from~\cite{svard2011}.
The physical parameters in~\eqref{eq:reactive.euler} are $A = 164180, T_A = 25, q_0 = 25$ and $\gamma = 1.4$. The initial condition on the physical domain $[-5, 25]$ is
\[
(\rho, v, p, Y) = \begin{cases}
(1.6812, 2.8867, 21.5682, 0), \quad & x < 0,\\
(1, 0, 1, 1), & x \ge 0.
\end{cases}
\]
Seeing the linear equation with stiff source term as a prototype~\eqref{eq:advection.stiff.source} for the reactive Euler equations~\eqref{eq:implicit.evolution},~\cite{svard2011} argued from~\eqref{eq:implicit.evolution} that the stiff source term is dissipating the solution.
Thus, they recommended reducing the dissipation from the numerical scheme to avoid getting a wrong shock location.
Since the source term only acts on the variable $\rho z$, we follow~\cite{svard2011} to reduce the dissipation in the numerical flux only for this variable by replacing the Rusanov dissipation~\eqref{eq:rusanov.wave.speed} for the variable $\rho z$ at the face $\eph$ with $W_k \lambda$ where
\begin{equation}
W_k = \begin{cases}
1 - \Delta t K(\avg{T}_\eph), \quad & \text{if } \Delta t K(\avg{T}_\eph) < 1,\\
0, & \text{otherwise},
\end{cases}
\label{eq:dissipation.reduction}
\end{equation}
where $\avg{T}_\eph = T((\avg{\uu}_{e+1} + \avg{\uu}_e)/2)$ and $\avg{u}_e$ is element mean of $u$ in element $e$ (Definition~\ref{defn:admissibility.preserving.means}).
Although the dissipation reduction strategy~\eqref{eq:dissipation.reduction} is shown to work reasonably well for the first order scheme in~\cite{svard2011}, it need not be enough for the high order schemes especially when combined with a subcell based limiter where the numerical fluxes are also used in the limiting procedure.
We show our results in Figure~\ref{fig:reactive.euler.rp1}, where it is observed that the strategy gives a \textit{spurious shock} for a coarse grid with $\Delta x = 0.3$.
However, the issue is also observed in Figure~\ref{fig:reactive.euler.rp1} to be resolved by choosing \emph{larger time steps} obtained by using $C_s = 2$~\eqref{eq:time.step.2d} instead of $C_s = 0.9$.
This is likely because the larger time step changes the dissipation behavior of the scheme.
The results on a finer grid with $\Delta x = 0.03$ do still show some under-shoot near the shock.
However, the order of the under-shoot is the same as the overshoot in the original scheme of~\cite{svard2011}, also shown in Figure~\ref{fig:reactive.euler.rp1}.
Although the results are not shown, the under-shoot disappears with further grid refinement.

\begin{figure}
\centering
\begin{tabular}{cc}
\includegraphics[width=0.4\columnwidth]{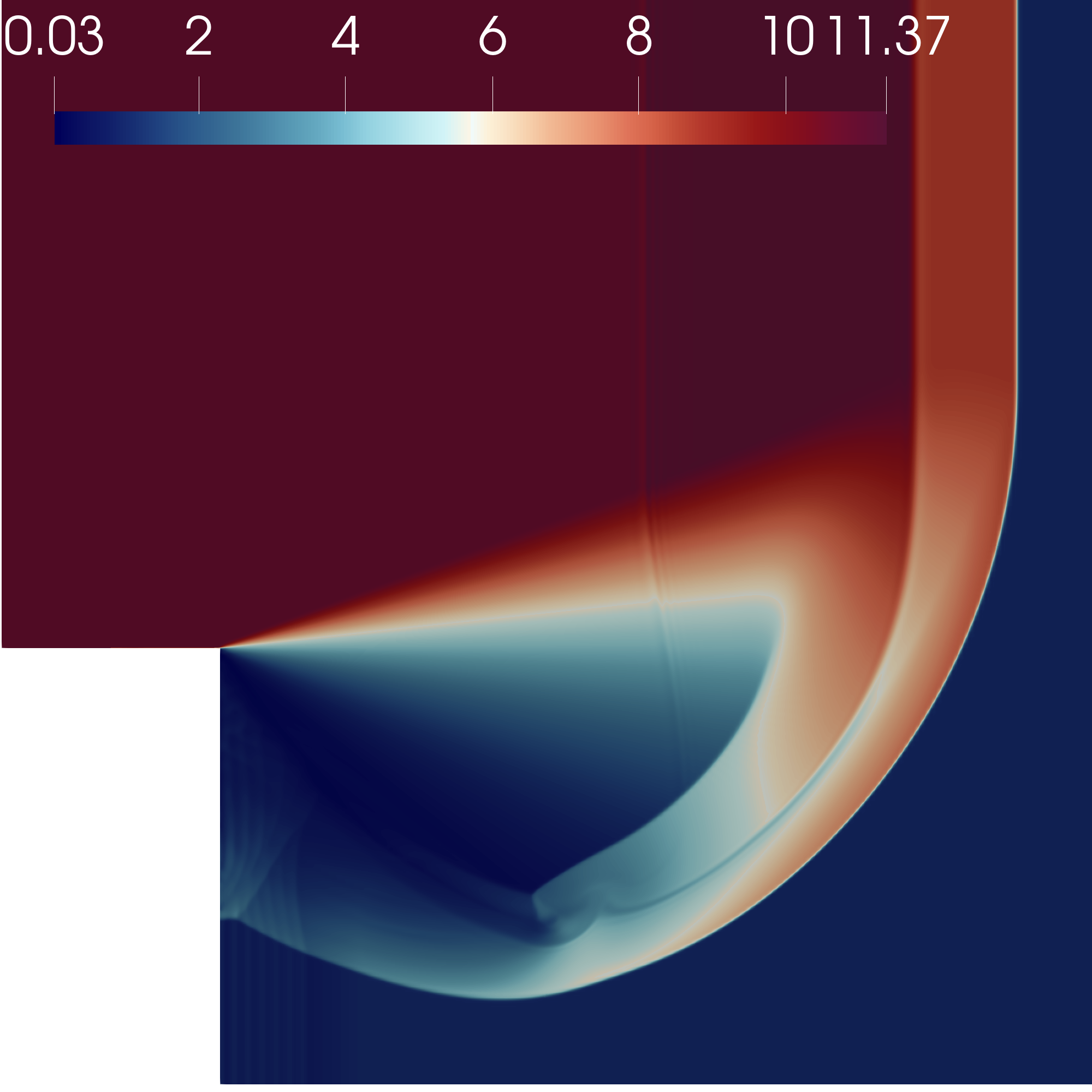} &
\includegraphics[width=0.4\columnwidth]{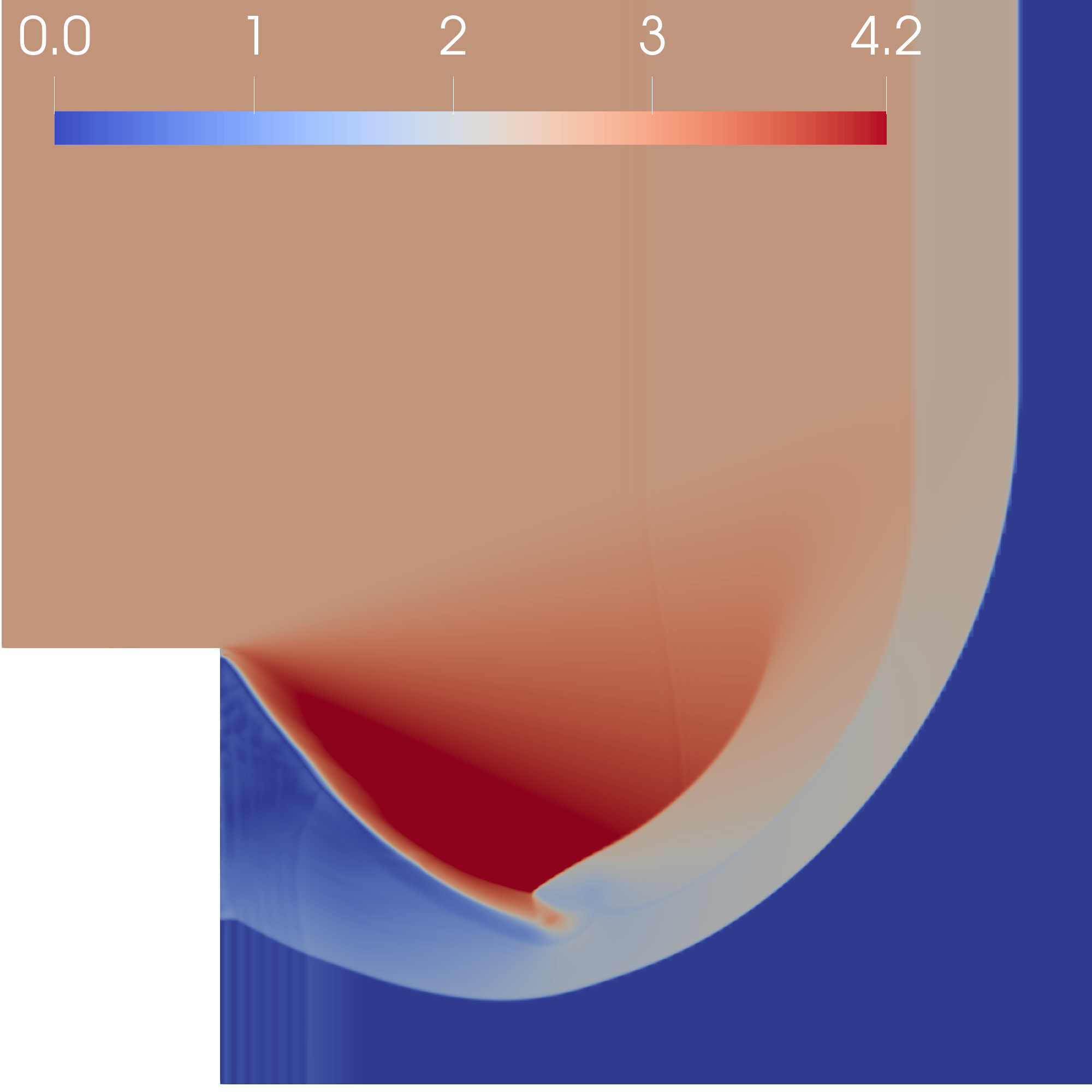} \\
(a) Density $\rho$ & (b) Mach number
\end{tabular}
\caption{Shock diffraction problem~\eqref{eq:shock.diffraction.ic} of reactive Euler's equations~\eqref{eq:reactive.euler} at $t=0.6$ using $\Delta x = \Delta y = 0.025$.}\label{eq:shock.diffraction}
\end{figure}

For the next test case, we consider a physically relevant non-stiff problem of shock diffraction over a step from~\cite{wang2012}. The physical domain is $[0,5]^2$ with a step placed at $[0,1] \times [0,2]$.
The parameters in~\eqref{eq:reactive.euler} are $A = 2566, T_A = 50, q_0 = 50$ and $\gamma = 1.2$ with the initial condition
\begin{equation}
(\rho, v_1, v_2, E, z) = \begin{cases}
(11, 6.18, 0, 970, 1), \quad & x < 0.5,\\
(1, 0, 0, 55, 1), & x \ge 0.5.
\end{cases}
\label{eq:shock.diffraction.ic}
\end{equation}
The reflective boundary conditions are applied at the walls and step, inflow boundary conditions at the left, and the outflow boundary condition is applied at the right boundary.
The simulation is run until time $t = 0.6$ using grid resolution $\Delta x = \Delta y = 0.025$ with polynomial degree $N=3$ and the results are shown in Figure~\ref{eq:shock.diffraction}, showing agreement with the results from~\cite{wang2012}.
\subsection{Ten moment equations} \label{sec:ten.moment}
\begin{figure}
\centering
\includegraphics[width=0.35\columnwidth]{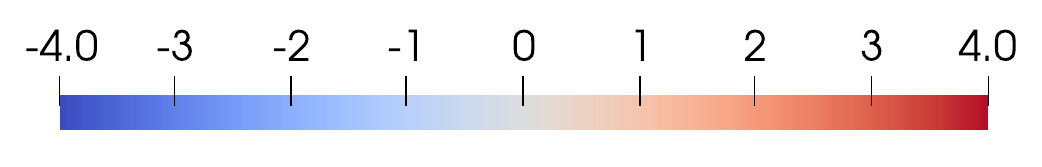}\\
\begin{tabular}{cc}
\includegraphics[width=0.35\columnwidth]{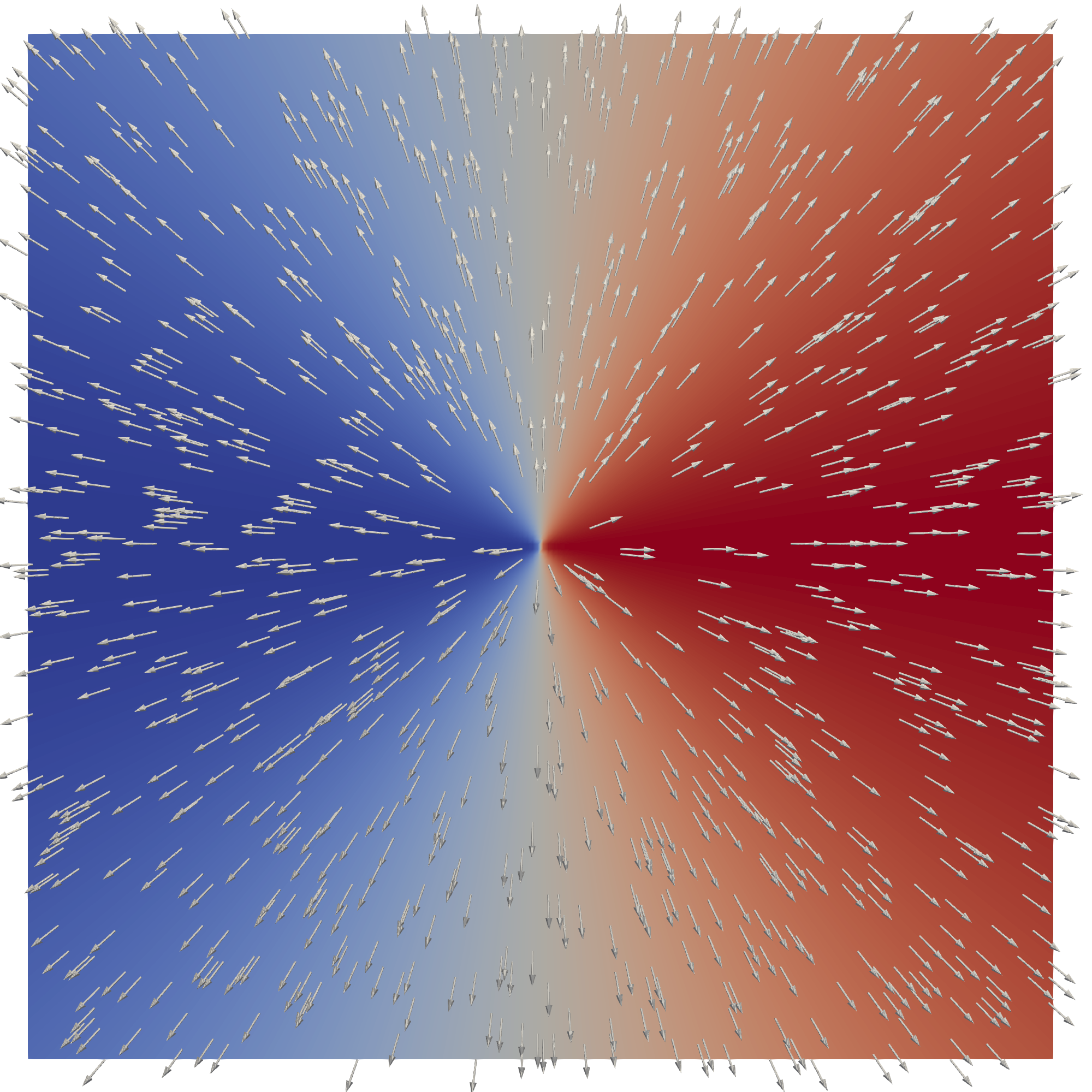}\hspace{1em} &
\includegraphics[width=0.35\columnwidth]{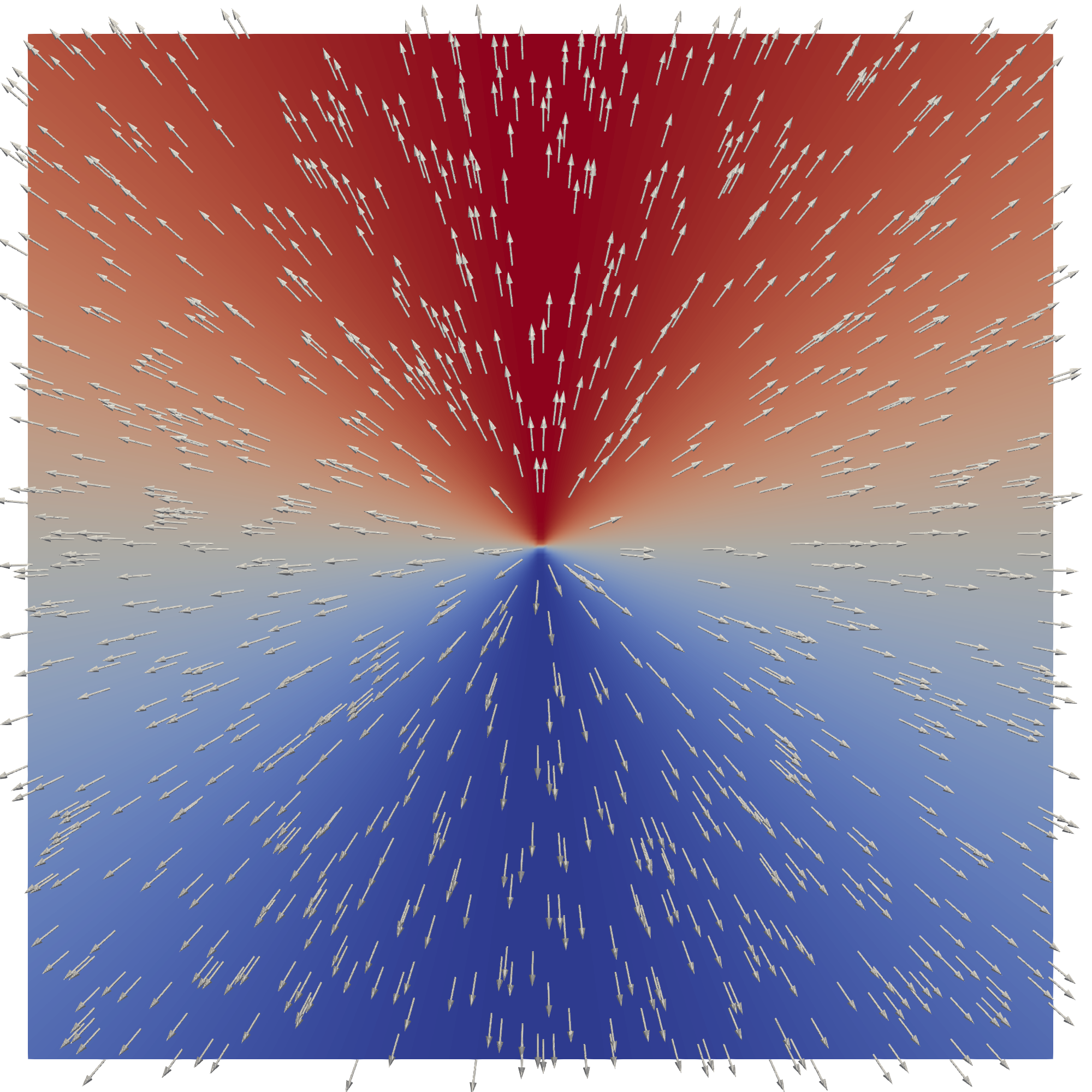} \\
(a) $v_x$ & (b) $v_y$
\end{tabular}
\caption{Initial velocity profile for the near-vacuum test of ten moment equations~\eqref{eq:tenmom.ic.2d}.}
\label{fig:ten.moment.2d.ic}
\end{figure}

\begin{figure}
\centering
\begin{tabular}{cccc}
% Color bars
\includegraphics[width=0.22\columnwidth]{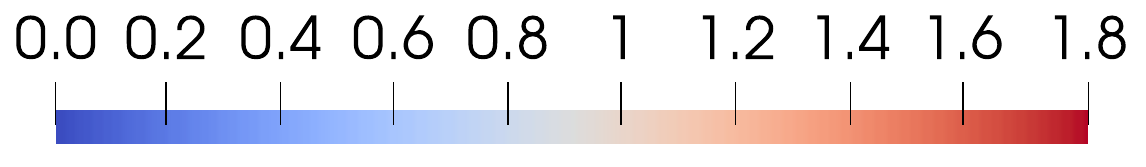} &
\includegraphics[width=0.22\columnwidth]{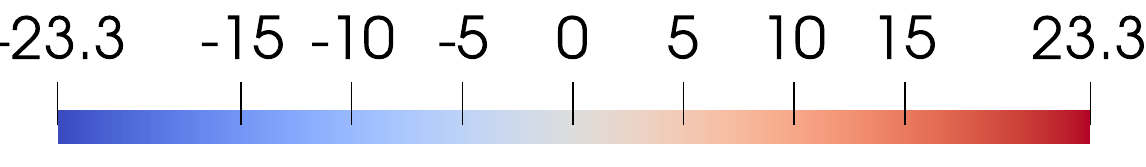} &
\includegraphics[width=0.22\columnwidth]{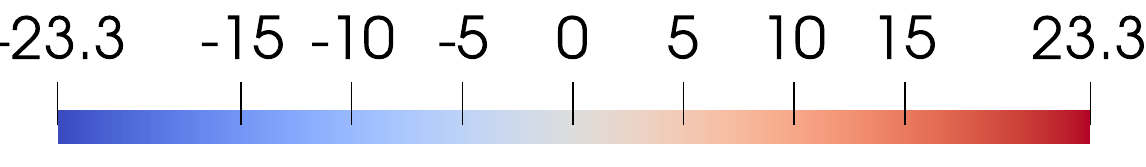} &
\includegraphics[width=0.22\columnwidth]{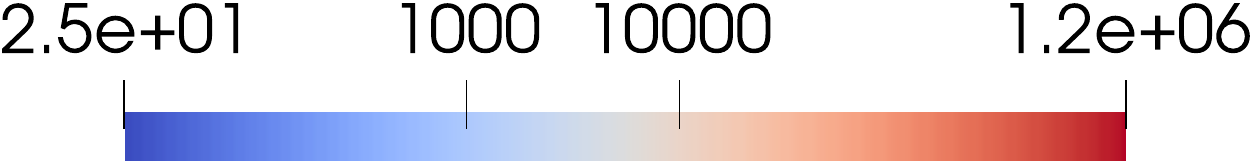} \\
\includegraphics[width=0.22\columnwidth]{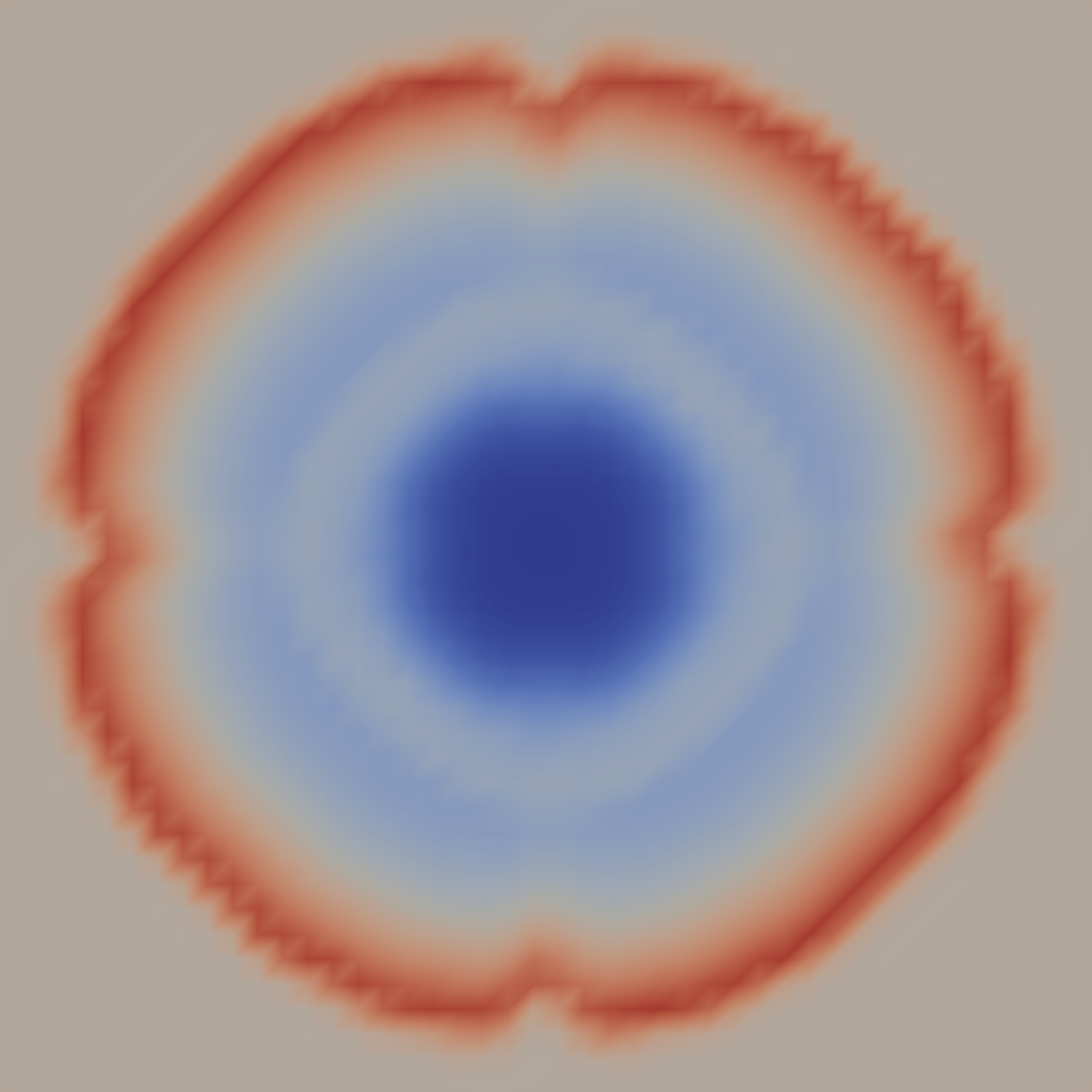} &
\includegraphics[width=0.22\columnwidth]{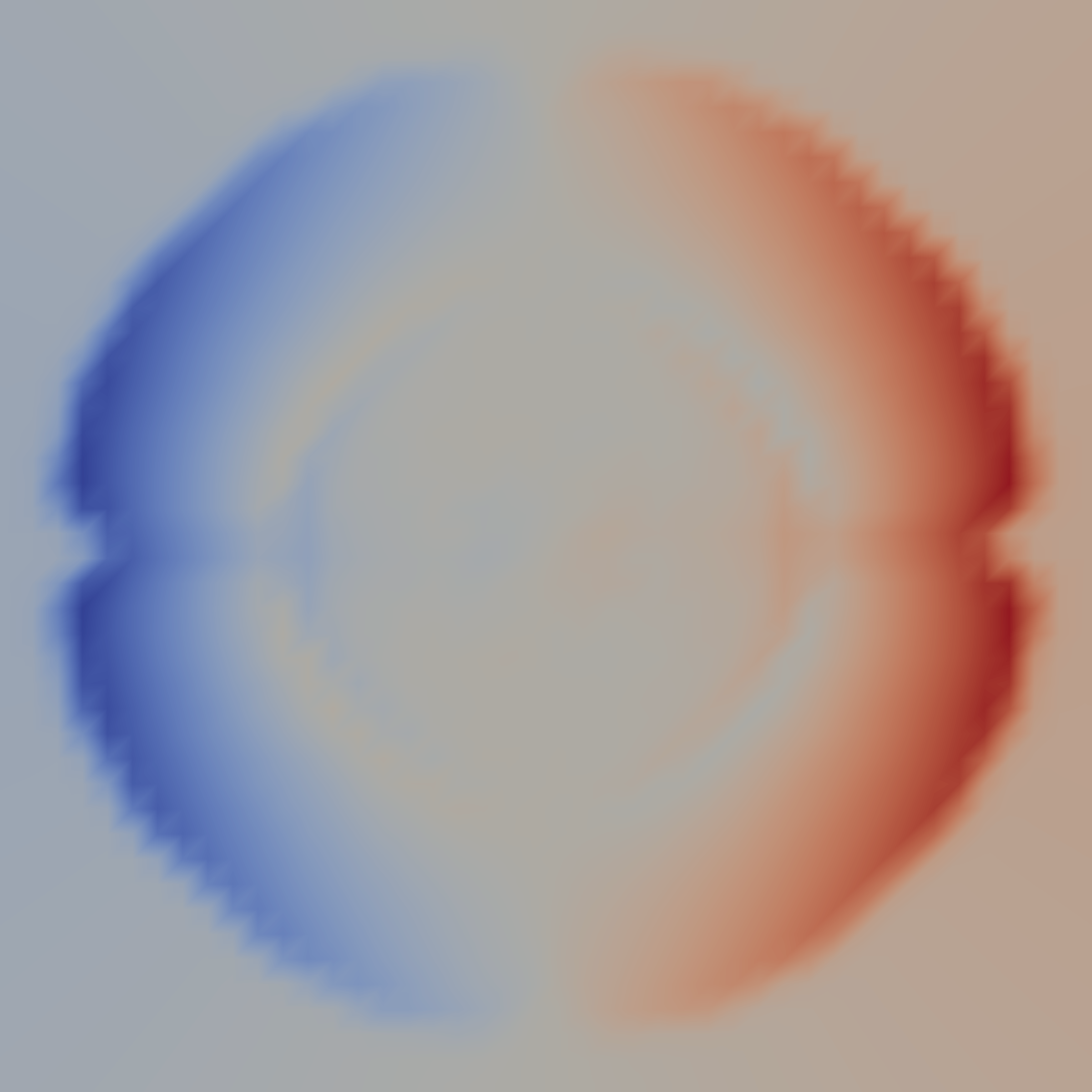} &
\includegraphics[width=0.22\columnwidth]{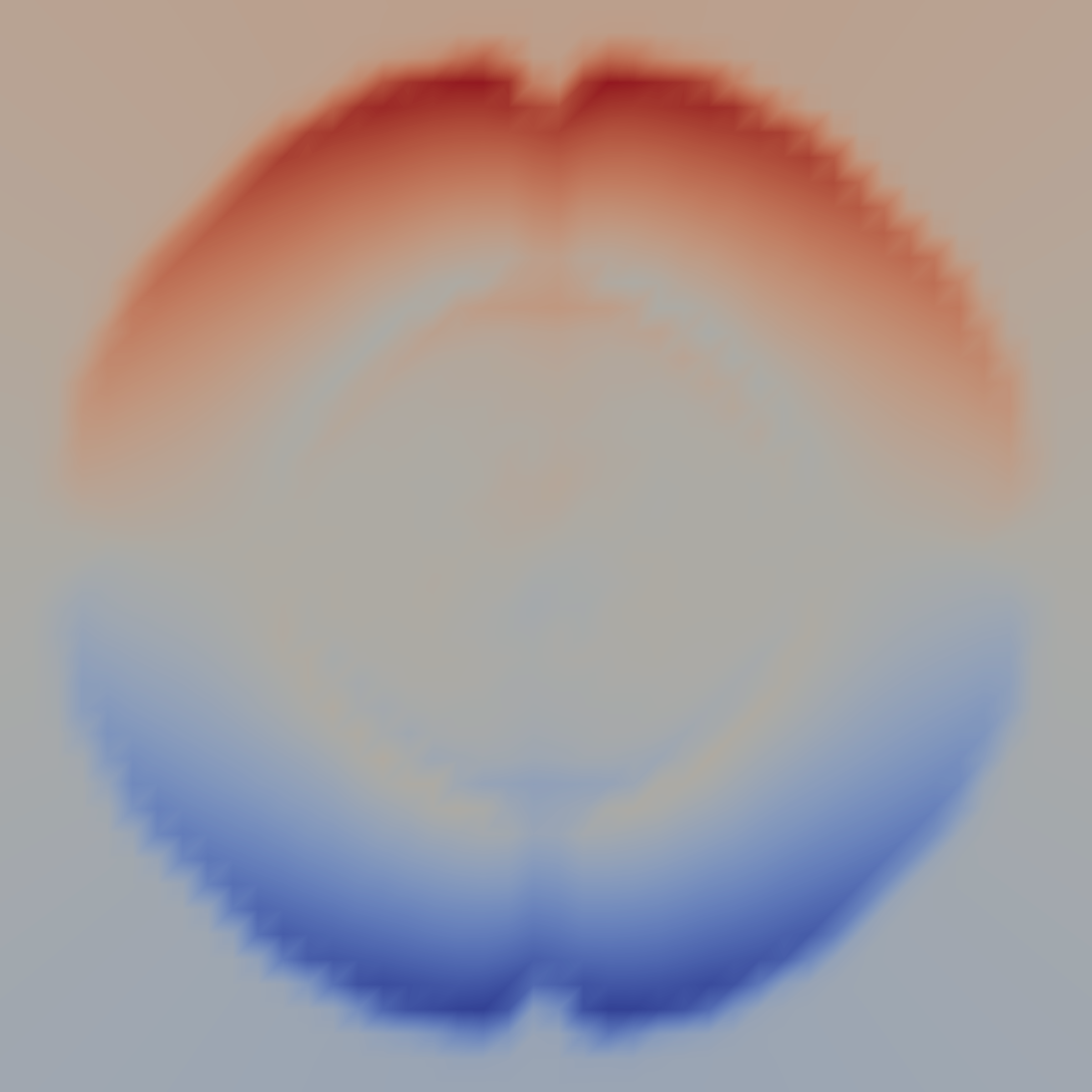} &
\includegraphics[width=0.22\columnwidth]{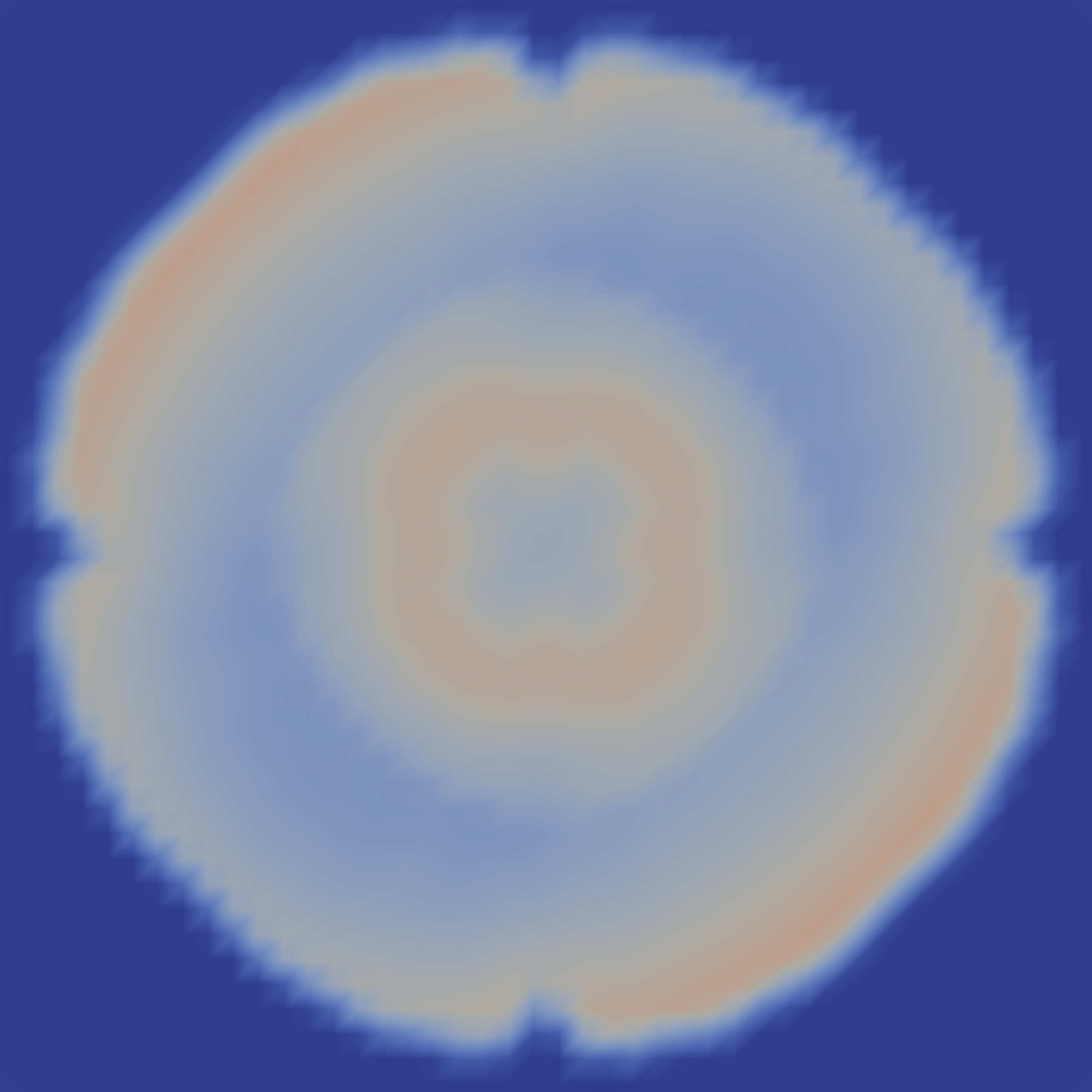} \\
(a) $\rho, \Delta x = 0.08$ & (b) $v_1, \Delta x = 0.08$ & (c) $v_2, \Delta x = 0.08$ & (d) $\det(P), \Delta x = 0.08$ \\
\includegraphics[width=0.22\columnwidth]{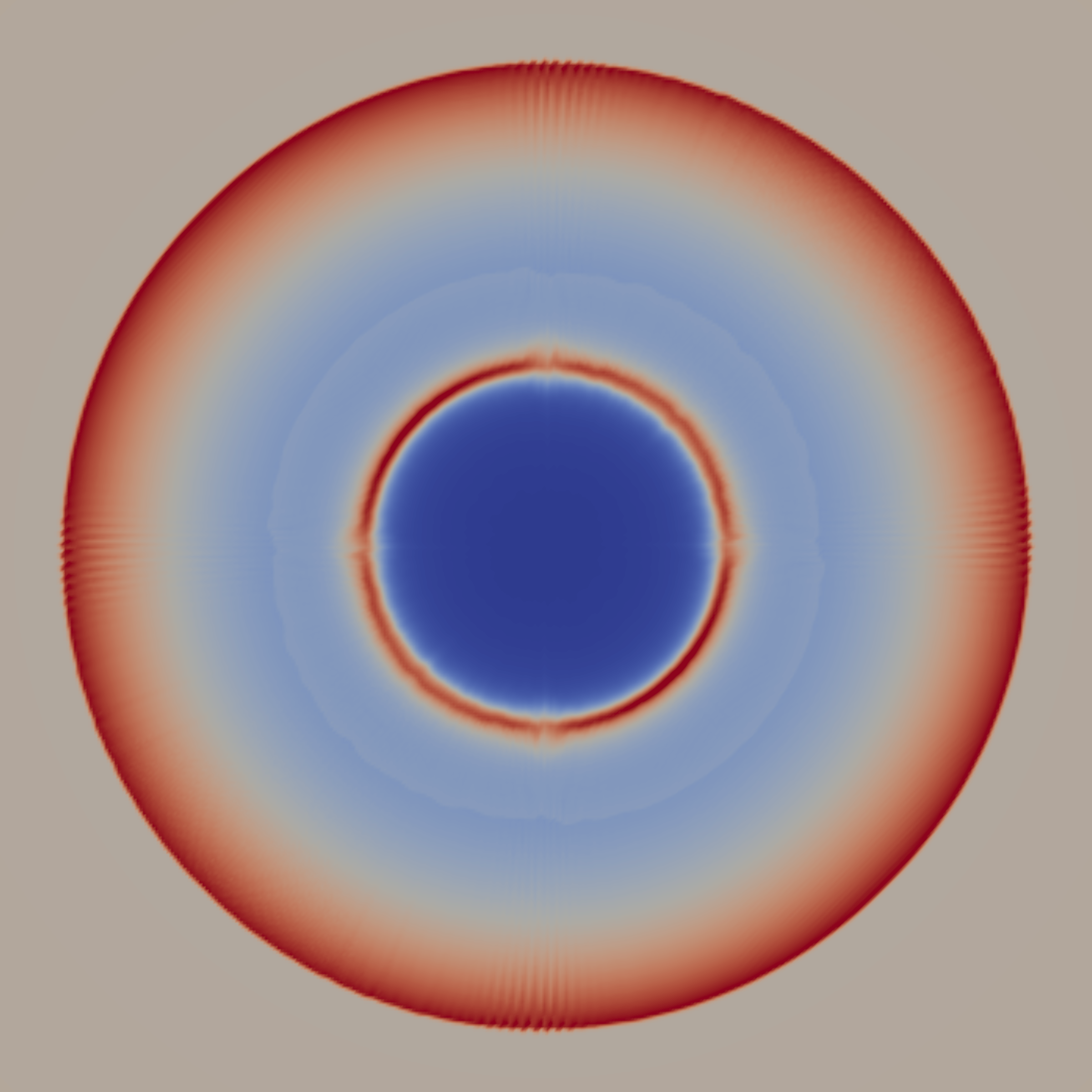} &
\includegraphics[width=0.22\columnwidth]{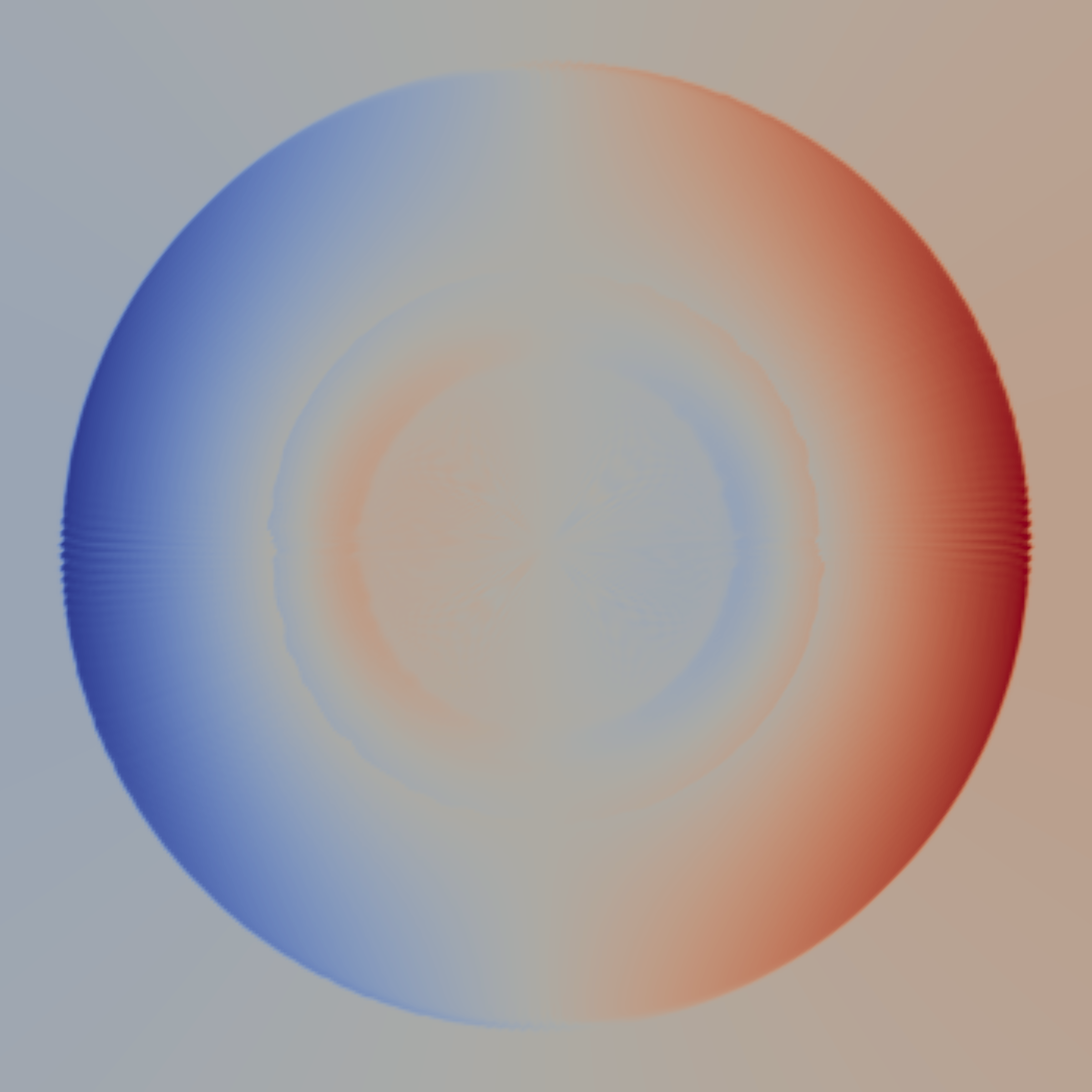} &
\includegraphics[width=0.22\columnwidth]{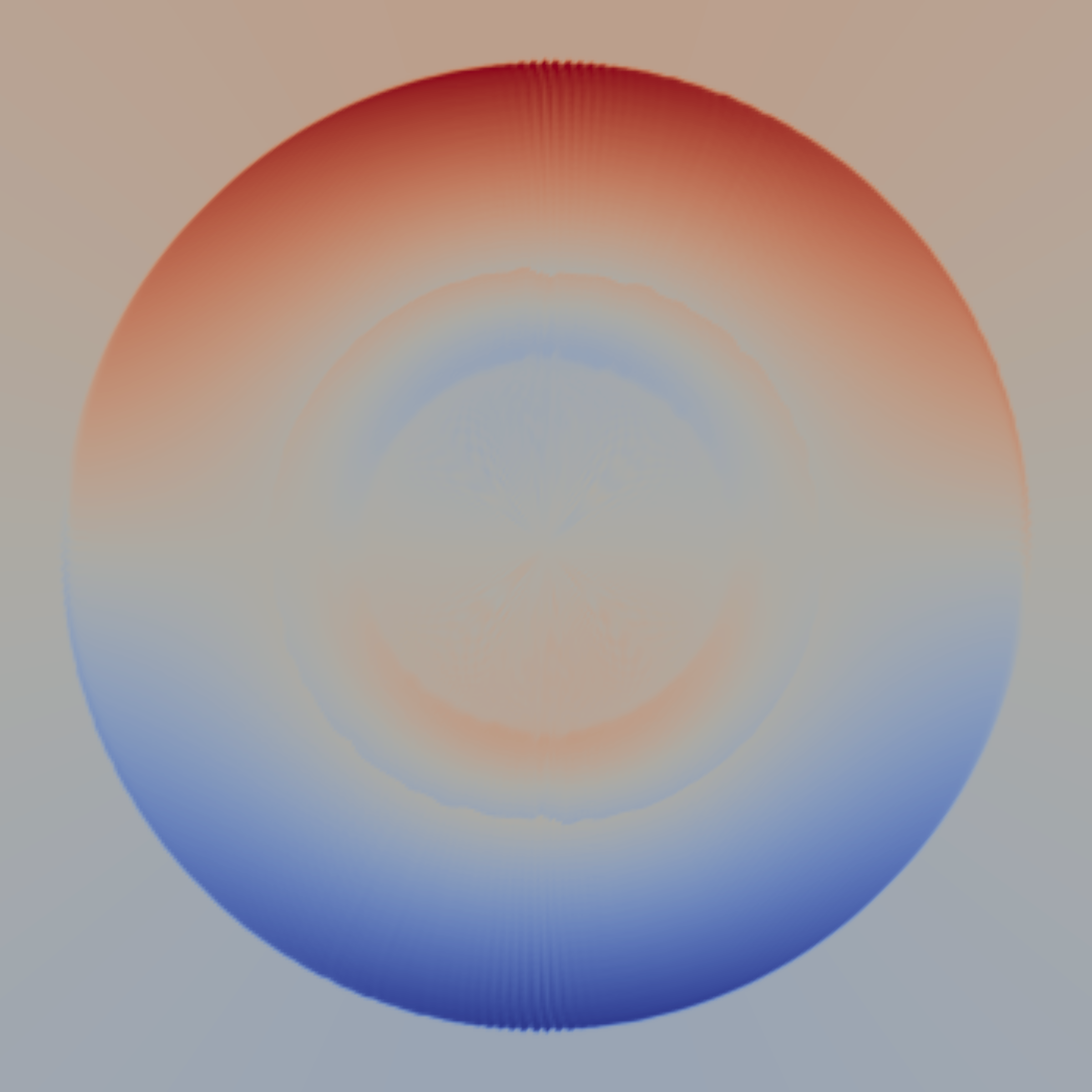} &
\includegraphics[width=0.22\columnwidth]{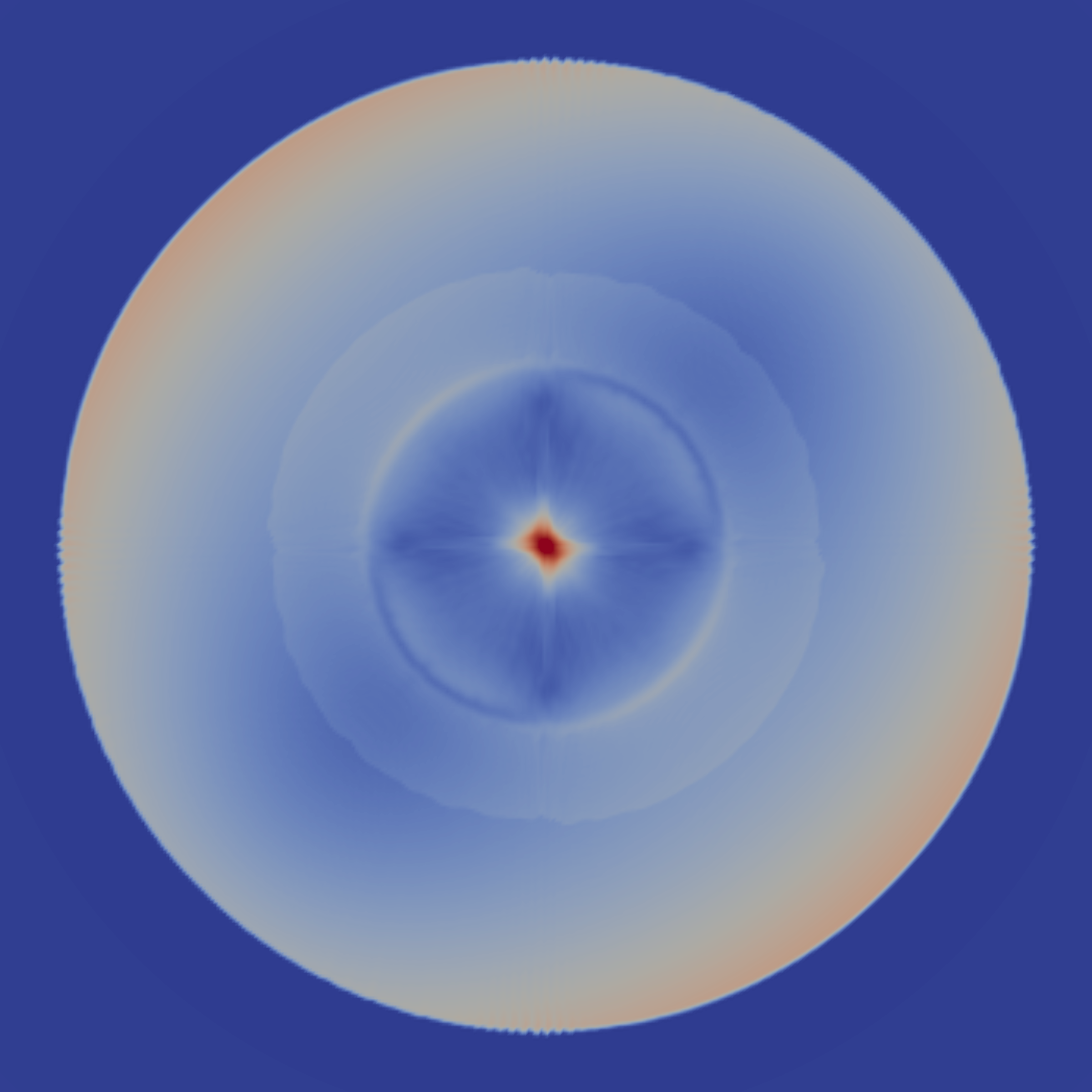} \\
(e) $\rho, \Delta x = 0.01$ & (f) $v_1, \Delta x = 0.01$ & (g) $v_2, \Delta x = 0.01$ & (h) $\det(P), \Delta x = 0.01$
\end{tabular}
\caption{Near-vacuum test of ten moment equations~\eqref{eq:tenmom.ic.2d} using source terms~\eqref{eq:tenmom.source} with stiffness parameter $K = 10^6$ at time $t=0.02$ using AGSA-IMEX(3,4,2)~\cite{boscarino2013} and polynomial degree $N=3$ for two different grid resolutions $\Delta x = 0.08$ (first row) and $\Delta x = 0.01$ (second row).}
\label{fig:ten.moment.2d}
\end{figure}

The ten moment problem~\cite{Berthon2015} can be seen as a generalization of the compressible Euler's equations where the energy and pressure are tensor quantities.
Specifically, the energy tensor is defined by the ideal equation of state $\bE=\frac 12 \bp + \frac 12 \rho \bv \otimes \bv$ where $\rho$ is the density, $\bv$ is the velocity vector, $\bp$ is the symmetric pressure tensor.
The conservative variables are $\uu = (\rho,\rho v_1,\rho v_2, E_{11}, E_{12}, E_{22})$, and the fluxes in~\eqref{eq:general.hyperbolic.equation} are given by
\begin{equation}
\pf_1 = \left[\begin{array}{c}
\rho v_1\\
\cR_{11} + \rho v_1^2\\
\cR_{12} + \rho v_1 v_2\\
\left( \cE_{11} + \cR_{11} \right) v_1\\
\Eonetwo v_1 + \half  ( \cR_{11} v_2 + \cR_{12} v_1 )\\
\Etwotwo v_1 + \cR_{12} v_2
\end{array}\right], \quad \pf_2 = \left[\begin{array}{c}
\rho v_2\\
\cR_{12} + \rho v_1 v_2\\
\cR_{22} + \rho v_2^2\\
\cE_{11} v_2 + \cR_{12} v_1\\
\cE_{12} v_2 + \half  ( \cR_{12} v_2 + \cR_{22} v_1 )\\
\left( \cE_{22} + \cR_{22} \right) v_2
\end{array}\right]. \label{eq:tmp}
\end{equation}
In order to construct a test problem with stiff source terms similar to~\eqref{eq:advection.stiff.source}, the source terms from~\cite{Berthon2015} are modified to include a scaling factor $K$ such that $\bss_K = \bss_K^x + \bss_K^y$, where
\begin{equation}\label{eq:tenmom.source}
\bss_K^{x} = K \left[\begin{array}{c}
0\\
- \half \rho \partial_x W\\
0\\
- \half \rho v_1 \partial_x W\\
- \frac{1}{4} \rho v_2 \partial_x W\\
0
\end{array}\right], \qquad  \bss_K^{y} = K \left[\begin{array}{c}
0\\
0\\
- \half \rho \partial_y W\\
0\\
- \frac{1}{4} \rho v_1 \partial_y W\\
- \half \rho v_2 \partial_y W
\end{array}\right], K \gg 1,
\end{equation}
where $W (x, y, t)$ is a given function, which models electron quiver energy in the laser~\cite{Berthon2015}.
The admissibility set for the ten moment equations is given by
\[
\Uad = \left\{ \uu = (\rho, \rho \bv, \bE) : \ad_1 = \rho > 0, \quad \ad_2 = \text{trace}\left(\bE - \dfrac{\rho \bv \otimes \bv}{2}\right) > 0,\quad \ad_3 = \text{det}\left(\bE - \dfrac{\rho \bv \otimes \bv}{2}\right) > 0 \right\} .
\]
The admissibility constraint $\ad_3$ is not concave as, taking $\rho = 1, \bv = \bzero$, the concavity of $\ad_3$ will imply that
\begin{equation} \label{eq:det.concave}
(E_{11}, E_{22}, E_{12}) \mapsto E_{11} E_{22} - E_{12}^2
\end{equation}
is concave, which is not true.
Being non-concave, the constraint $\ad_3$ cannot be enforced by the formula~\eqref{eq:theta.formula}, and thus the non-linear equation~\eqref{eq:theta.eqn} is solved using Newton's method to enforce $\ad_3$.
We now test the scheme on a challenging near-vacuum problem from~\cite{meena2017} with the modified stiff source term~\eqref{eq:tenmom.source}.
The near-vacuum state arises because the solution consists of two rarefactions~\cite{meena2017}, leading to low density values in the solution.
The initial condition in one spatial dimension is given by
\begin{equation}
(\rho, v_1, v_2, P_{11}, P_{12}, P_{22}) = \begin{cases}
(1, -4, 0, 9, 7, 9), \quad & x_1 < 0,\\
(1, 4, 0, 9, 7, 9), & x_1 \ge 0,
\end{cases} \label{eq:tenmom.ic.1d}
\end{equation}
with the physical domain $[0,4]$.
The electron quiver energy function in the source term~\eqref{eq:tenmom.source} is $W(x) = \frac{1}{40} \exp(-20(x-2)^2)$, and the stiffness parameter~\eqref{eq:tenmom.source} is taken to be $K = 10^5$.
The simulation is run till time $t=0.1$ using polynomial degree $N=3$ and time discretization SSP3-IMEX(4,3,3) of ~\cite{pareschi2005}~\eqref{eq:ssp33.butcher}.
\begin{figure}
\centering
\begin{tabular}{cc}
\includegraphics[width=0.45\columnwidth]{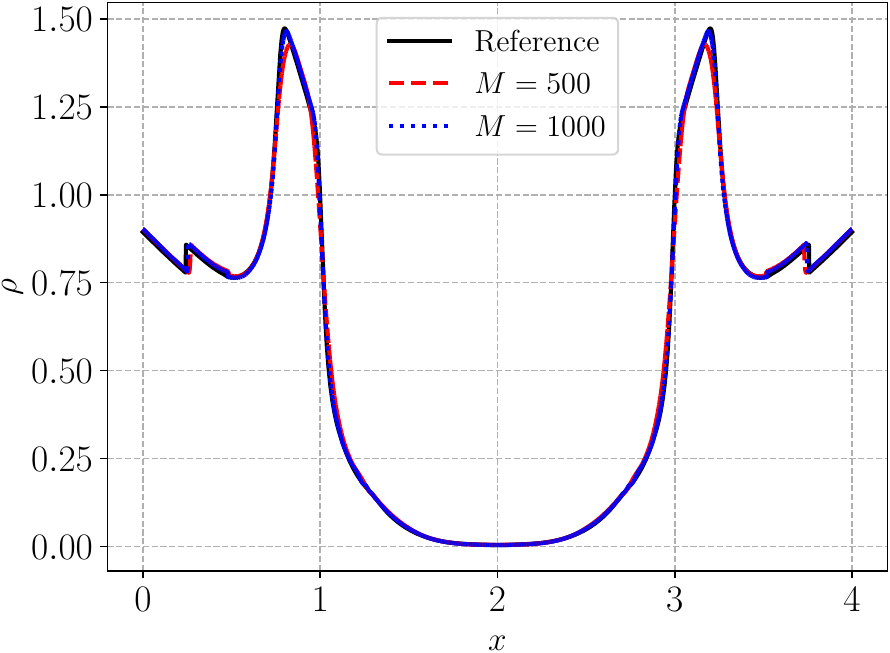} &
\includegraphics[width=0.45\columnwidth]{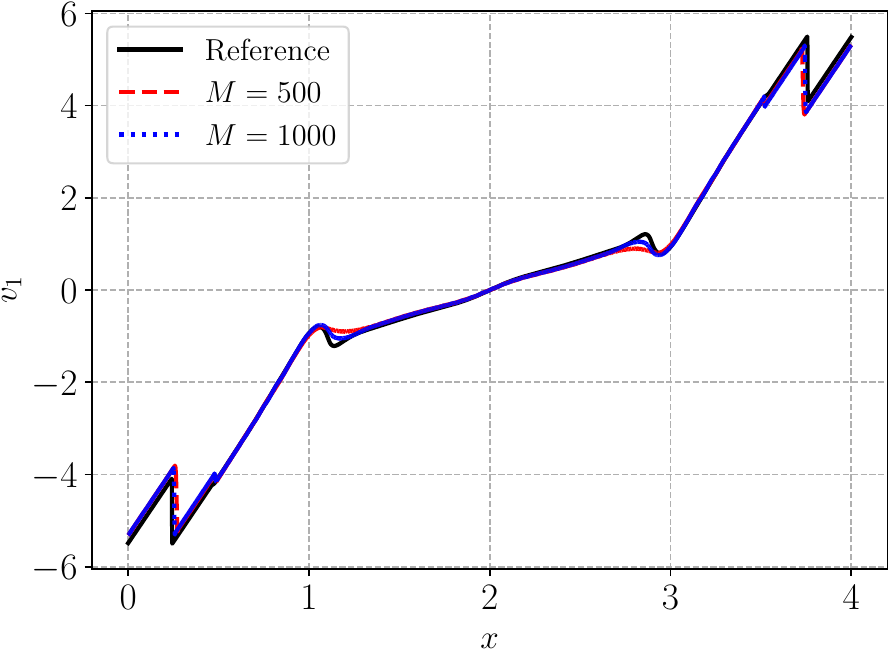} \\
(a) $\rho$ & (b) $v_1$
\end{tabular}
\caption{1-D near-vacuum test of ten moment equations~\eqref{eq:tenmom.ic.1d} using source terms~\eqref{eq:tenmom.source} with stiffness parameter $K = 10^5$ at time $t=0.1$ using SSP3-IMEX(4,3,3)~\cite{pareschi2005} and polynomial degree $N=3$ for $M=500, 1000$ elements showing (a) $\rho$, (b) velocity $v_1$. The reference solution is computed using a first order finite volume scheme on a grid of $10^5$ cells.}
\label{fig:ten.moment.1d}
\end{figure}
The problem is extended to two dimensions in the normal direction by taking the domain to be $[-2,2]^2$ and the initial condition to be
\begin{equation}
(\rho, \bv, P_{11}, P_{12}, P_{22}) = (1, 4 \frac{\bx}{\|\bx\|}, 9, 7, 9).
\label{eq:tenmom.ic.2d}
\end{equation}
The initial velocity profile is shown in Figure~\ref{fig:ten.moment.2d.ic}.
The electron quiver energy function is $W(x,y) = \frac{1}{40} \exp(-20(x^2 + y^2))$, and the stiffness parameter~\eqref{eq:tenmom.source} is $K = 10^6$, making it a tougher problem than the 1-D case.
This test case is quite challenging due to the near-vacuum condition in the initial condition~\eqref{eq:tenmom.ic.2d} along with the stiff source term~\eqref{eq:tenmom.source}.
Thus, an explicit method or an IMEX method without the admissibility enforcing procedure (Algorithms~\ref{alg:flux.limiting}, \ref{alg:alpha.limiting}) gives inadmissible solutions very quickly.
The simulation is run till time $t = 0.02$ using polynomial degree $N=3$ and time discretization AGSA-IMEX(3,4,2) of ~\cite{boscarino2013}~\eqref{eq:ssp33.butcher}.
The AGSA-IMEX(3,4,2) is shown to demonstrate the flexibility of the cRKFR IMEX scheme to work with different time discretizations, but other time discretization including HT (1,1,2)~\eqref{eq:ht112.butcher} also work for this problem.
The results are shown in Figure~\ref{fig:ten.moment.2d} for two different grid resolutions $\Delta x = 0.08$ and $\Delta x = 0.01$.
Since the time step size is determined by the advective time step restriction~\eqref{eq:time.step.2d}, the coarse grid simulation with $\Delta x = 0.08$ uses a larger time step size.
This makes the coarse grid simulation a good test for the cRKFR IMEX scheme's ability to handle stiff source terms while using large time steps.
In addition, coarse grid simulations are practical because of the high computational cost of fine grids, especially in higher spatial dimensions.

\subsection{Shear shallow water equations} \label{sec:ssw}
\begin{figure}
\centering
\begin{tabular}{cc}
\includegraphics[width=0.45\columnwidth]{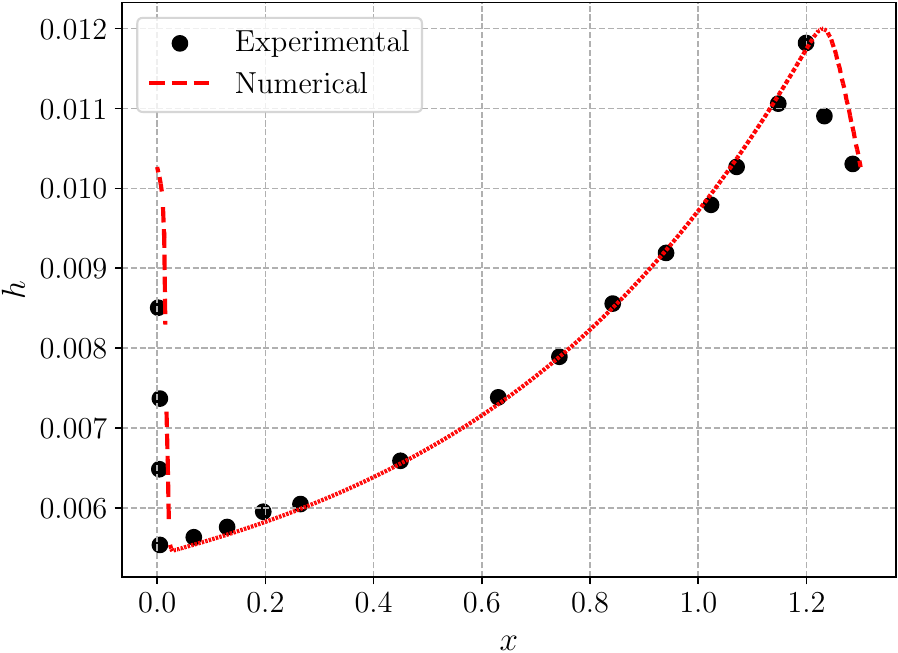} &
\includegraphics[width=0.45\columnwidth]{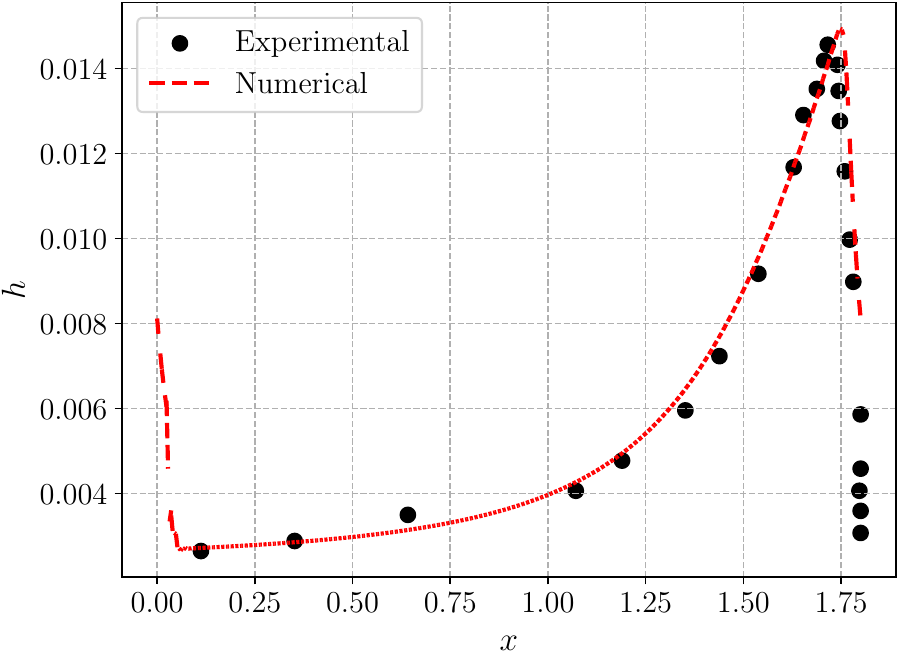} \\
(a) & (b)
\end{tabular}
\caption{Simulation of 1-D roll wave problem, height profile for (a) Case 1~\eqref{eq:roll.wave.ic.1d} and (b) Case 2~\eqref{eq:roll.wave.ic.2d} showing comparison against Brock's experiments~\cite{Brock1969,Brock1970} using polynomial degree $N=3$ and 200 elements for case 1 and 150 elements for case 2.}
\label{fig:roll.wave}
\end{figure}

\begin{figure}
\centering
\includegraphics[width=0.9\columnwidth]{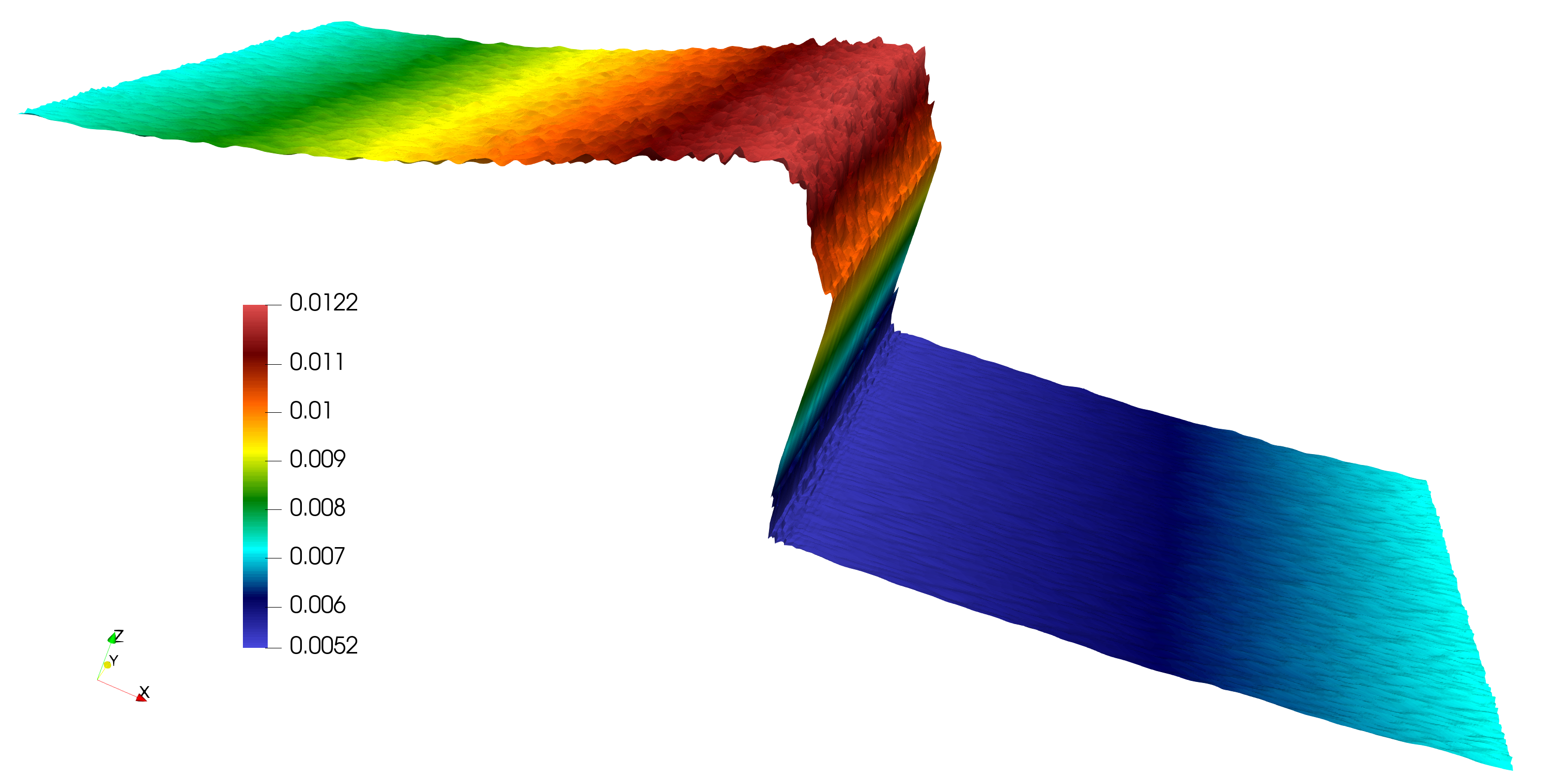}
\caption{Simulation of 2-D roll wave problem~\eqref{eq:roll.wave.ic.2d} showing elevated water height using polynomial degree $N=3$ on a $130 \times 50$ mesh.}
\label{fig:roll.wave.2d}
\end{figure}

The next test case consists of the non-conservative shear shallow water equations from~\cite{chandrashekar2020ssw}.
The unknown variables are given by the fluid height $h$, momentum components $hv_1, hv_2$ and the symmetric Reynolds stress tensor $\bSR$ whose components are $\SR_{11}, \SR_{12}, \SR_{22}$.
Similar to the ten moment equations~\eqref{eq:tmp}, the energy is a tensor that relates to the Reynolds stress tensor $\bSR$ and pressure tensor $\bSP$ as
\[
\bSR := h \bSP, \qquad \bSE := \half \bSR + h \half \bv \otimes \bv.
\]
The fluxes~\eqref{eq:general.hyperbolic.equation} are given by
\[
\fx = \left[ \begin{array}{c}
hv_1\\
\SR_{11} + hv_1^2 + \half gh^2\\
\SR_{12} + hv_1 v_2\\
(\SE_{11} + \SR_{11}) v_1\\
\SE_{12} v_1 + \half (\SR_{11} v_2 + \SR_{12} v_1)\\
\SE_{22} v_1 + \SR_{12} v_2
\end{array} \right], \quad \fy = \left[ \begin{array}{c}
hv_2\\
\SR_{12} + hv_1 v_2\\
\SR_{22} + hv_2^2 + \half gh^2\\
\SE_{11} v_2 + \SR_{12} v_1\\
\SE_{12} v_2 + \half (\SR_{12} v_2 + \SR_{22} v_1)\\
(\SE_{22} + \SR_{22}) v_2
\end{array} \right].
\]
The non-conservative matrices and source terms in~\eqref{eq:general.hyperbolic.equation} are given by
\[
\B_1 = \left[ \begin{array}{c}
0\\
0\\
0\\
ghv_1\\
\half ghv_2\\
0
\end{array} \right], \qquad \B_2 = \left[ \begin{array}{c}
0\\
0\\
0\\
0\\
\half ghv_1\\
ghv_2
\end{array} \right], \qquad \bss = \left[ \begin{array}{c}
0\\
- gh \df{\btop}{x_1} - C_f | \vel |v_1\\
- gh \df{\btop}{x_2} - C_f | \vel |v_2\\
- ghv_1 \df{\btop}{x_1} + \half h \Diss_{11} - C_f | \vel |v_1^2\\
- \half ghv_2 \df{\btop}{x_1} - \half ghv_1 \df{\btop}{x_2} + \half h
\Diss_{12} - C_f | \vel |v_1 v_2\\
- ghv_2 \df{\btop}{x_2} + \half h \Diss_{22} - C_f | \vel |v_2^2
\end{array} \right],
\]
where the $g$ is the gravitational acceleration 9.81, $C_f$ is the Chezy coefficient, $\btop(x_1, x_2)$~\cite{ivanova2017} is the bottom topography which we take in our tests to be $b(x_1,x_2) = -x\tan(\theta_r)$, where $\theta_r$ is the inclination angle, and $\Diss_{11}, \Diss_{12}, \Diss_{22}$ are the components of the dissipation tensor
\[
\Diss = - \frac{2 \alpha |\vel|^3}{h} \bSP, \qquad
\alpha = \max \left(0, C_r \frac{T - \phi h^2}{T^2}\right),\qquad T = \text{trace}(\bSP) = \SP_{11} + \SP_{22} ,
\]
and $C_r, \phi$ are model constants that are calibrated using experiments.
The initial condition for the 1-D roll wave problem from~\cite{chandrashekar2020ssw} is given by
\begin{equation}
\begin{gathered}
(h ,v_1, v_2, \SR_{11}, \SR_{12}, \SR_{22}) =
(
h_0 (1 + a \sin(2 \pi x/L_x)),
\sqrt{g h_0 \tan(\theta_r)/C_f},
0,
\phi h^2/2,
0,
\phi h^2/2
),
\end{gathered}
\label{eq:roll.wave.ic.1d}
\end{equation}
and the physical domain to be $[0, L_x]$.
For test case called case 1 in~\cite{chandrashekar2020ssw}, the physical parameters are taken to be
\begin{equation}
C_f = 0.0036,\ C_r = 0.00035,\  h_0 = 0.00798,\ \theta_r = 0.05011,\ a = 0.05,\ L_x = 1.3,\ \phi = 22.76.
\label{eq:case1.params}
\end{equation}
In case 2 of~\cite{chandrashekar2020ssw}, the parameters are taken to be
\begin{equation}
C_f = 0.0038,\  C_r = 0.002,\ h_0 = 0.00533,\ \theta_r = 0.119528,\ a = 0.05,\ L_x = 1.8,\ \phi = 153.501.
\label{eq:case2.params}
\end{equation}
For the 2-D problem, the height in~\eqref{eq:roll.wave.ic.1d} is replaced by
\begin{equation}
h(x_1,x_2) = h_0 (1 + a \sin(2 \pi x_1/L_{x_1}) + a \sin(2 \pi x_2/L_{x_2})),
\label{eq:roll.wave.ic.2d}
\end{equation}
in the 2-D domain $[0,L_x] \times [0,L_y]$, while other variables remain the same as in~\eqref{eq:roll.wave.ic.1d}.
The parameters in the 2-D test are taken as in case 1~\eqref{eq:case1.params} with $L_y = 0.5$.
The height profiles for case 1 and case 2~(\ref{eq:case1.params},~\ref{eq:case2.params}) in 1-D using the cRKFR scheme are shown in Figures~\ref{fig:roll.wave}a, b respectively.
Case 1 is run till $t=26.99$ using polynomial degree $N=3$ and 200 elements, and case 2 till $t=26.35185$ using polynomial degree $N=3$ and 150 elements.
In both cases, as in~\cite{chandrashekar2020ssw}, Figure~\ref{fig:roll.wave} shows a comparison of the height profile obtained by the cRKFR scheme against data from Brock's experiments~\cite{Brock1969,Brock1970} and an agreement similar to~\cite{chandrashekar2020ssw} is observed.
This shows that our non-conservative treatment is working well.
In Figure~\ref{fig:roll.wave.2d}, the 2-D results are shown at time $t=36$ using polynomial degree $N=3$ on a $130 \times 50$ mesh.
The solution resembles the results of~\cite{chandrashekar2020ssw} showing similar turbulent structures.
The advantage of using the shear shallow water model over the standard shallow water equations is also noted in these experiments as, along with a hydraulic jump, there is always a smooth roll wave profile immediately behind it.
A standard shallow water model would only be able to capture the hydraulic jump~\cite{chandrashekar2020ssw}.
In these tests, because of the low water height values, the positivity enforcing procedure (Algorithms~\ref{alg:flux.limiting}, \ref{alg:alpha.limiting}) was essential to prevent physically inadmissible solutions.

\subsection{Magnetohydrodynamics equations} \label{sec:mhd}
We solve the ideal GLM-MHD equations from~\cite{derigs2018,derigs2016} (specifically, equation (3.16) of~\cite{derigs2018})
% \begin{align}
\begin{equation}
\begin{split}
     \pd{}{t} \left(\begin{array}{c}
       \rho\\
       \rho \uu\\
       E\\
       \bB\\
       \psi
     \end{array}\right) &+ \nabla \cdot \left(\begin{array}{c}
       \rho \uu\\
       \rho (\uu \otimes \uu) + (p + \| \bB \|^2 /
       2)  I - \bB \otimes \bB\\
       \uu  (E + p + \| \bB \|^2 / 2) - \bB
       (\uu \cdot \bB)\\
       \bB \otimes \uu - \uu \otimes \bB\\
       c_h  \bB
     \end{array}\right) = - \Upsilon_{\text{GLM}},
     \\
     & \Upsilon_{\text{GLM}} := (\nabla \cdot \bB)
     \left[\begin{array}{c}
       0\\
       \bB\\
       \uu \cdot \bB\\
       \uu\\
       0
     \end{array}\right] + (\nabla \psi) \cdot \left[\begin{array}{c}
       \bzero\\
       0\\
       \uu \psi\\
       0\\
       \uu
     \end{array}\right],
\end{split} \label{eq:glm.mhd}
\end{equation}
% \end{align}
where $\rho, \rho \uu, E$ are the mass, momenta, and total specific energy of the plasma system, $p$ is the thermal pressure, $I$ is the identity matrix, $\bB$ is the magnetic field, and $\psi$ is the divergence-correcting field and $c_h$ is the hyperbolic divergence cleaning speed~\cite{derigs2018,derigs2016}.
The thermal pressure relates to the other variables through the ideal equation of state
\[
p = (\gamma - 1) \epsilon, \qquad \epsilon = E - \frac 12 \rho \norm{u}^2 - \frac 12 \norm{\bB}^2 - \frac 12 \psi^2.
\]
The first term in $\Upsilon_{\text{GLM}}$ is the Godunov-Powell source term~\cite{godunov2025,Powell1997} that makes the system hyperbolic and symmetrizable.
The $\psi$ is obtained as a generalized Lagrange multiplier (GLM) to control the divergence of the magnetic field $\bB$.
If the initial solution does not satisfy the divergence-free constraint, then the divergence cleaning will ensure that the deviations will decay instead of increasing or being convected by the flow.
Thus, the system~\eqref{eq:glm.mhd} ensures that any magnetic divergence caused by inaccuracies of a numerical method for the ideal GLM-MHD system remains small.
As $c_h \to \infty$, the GLM-MHD equations converge to the divergence constraint.
We choose $c_h$ to be just so large that the time step restriction from the hyperbolic GLM-MHD equations is not more severe than that from the ideal MHD equations~\cite{derigs2018,derigs2016}.

The fluxes and non-conservative products for the GLM-MHD equations~\eqref{eq:glm.mhd} are incorporated into \texttt{Tenkai.jl}~\cite{tenkai} using \texttt{Trixi.jl}~\cite{schlottkelakemper2021purely,schlottkelakemper2020trixi,ranocha2022adaptive} through Julia's package manager.
\subsubsection{Alfvén wave convergence test} \label{sec:mhd.convergence}

\begin{figure}
\centering
\begin{tabular}{cc}
\includegraphics[width=0.38\columnwidth]{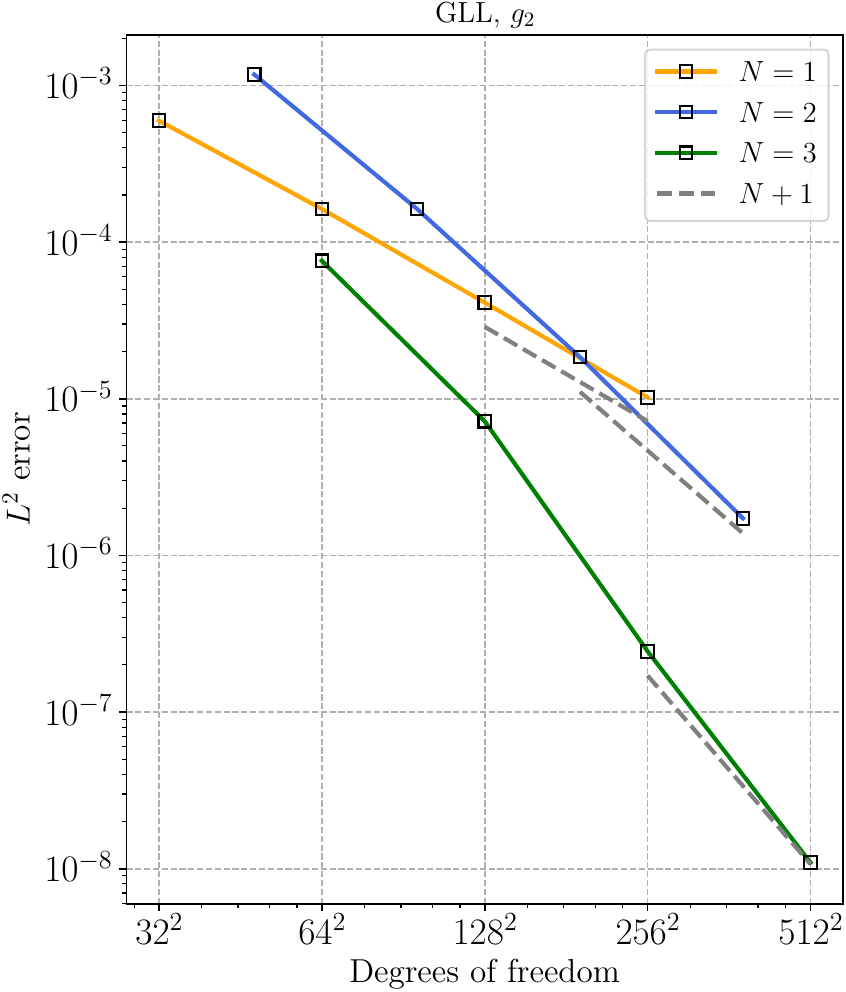} &
\includegraphics[width=0.38\columnwidth]{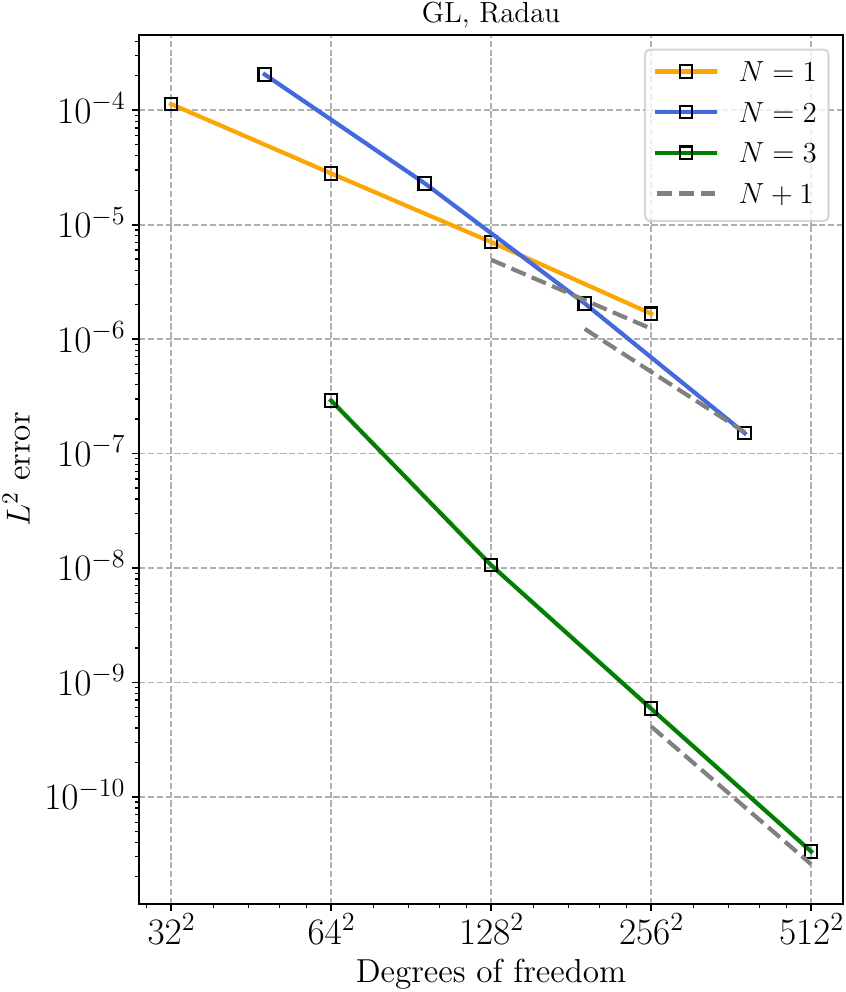} \\
(a) & (b)
\end{tabular}
\caption{Alfvén wave MHD convergence test (a) $g_2$ correction and GLL points, (b) Radau correction and GL points.} \label{fig:alfven.wave}
\end{figure}

The smooth Alfvén wave test from~\cite{toth2000} is used to verify the order of accuracy of the scheme~\cite{derigs2016,derigs2018}.
It consists of propagation of a circularly polarized wave.
The exact solution is given by
\begin{equation*}
\begin{gathered}
\rho = 1, \quad \bv = v_\perp (- \sin \alpha, \cos \alpha, 0.1 \cos(2 \pi (x_\perp + t))), \quad p = 0.1, \\
\bB = (\cos \alpha + v_1, \sin \alpha + v_2, v_3),
\end{gathered}
\end{equation*}
where
\[
x_\perp = x \cos \alpha + y \sin \alpha, \quad B_{\|} = 1, \quad v_\perp = 0.1 \sin(2 \pi x_\perp),
\]
and $\alpha = \pi/4$ is the angle of propagation with respect to the $x$-axis.
The domain is $[0, 1/\cos \alpha] \times [0, 1 / \sin \alpha] = [0, \sqrt 2]^2$ with periodic boundary conditions.
The simulation is run till $t=2$ and the convergence analysis is performed for polynomial degrees $N=1,2,3$ showing the optimal order of accuracy in Figure~\ref{fig:alfven.wave}.
It is also observed that the GL points with Radau correction give a lot better accuracy than the GLL points with $g_2$ correction.

\subsubsection{Orszag-Tang vortex}

\begin{figure}
\centering
\begin{tabular}{cc}
\includegraphics[width=0.45\columnwidth]{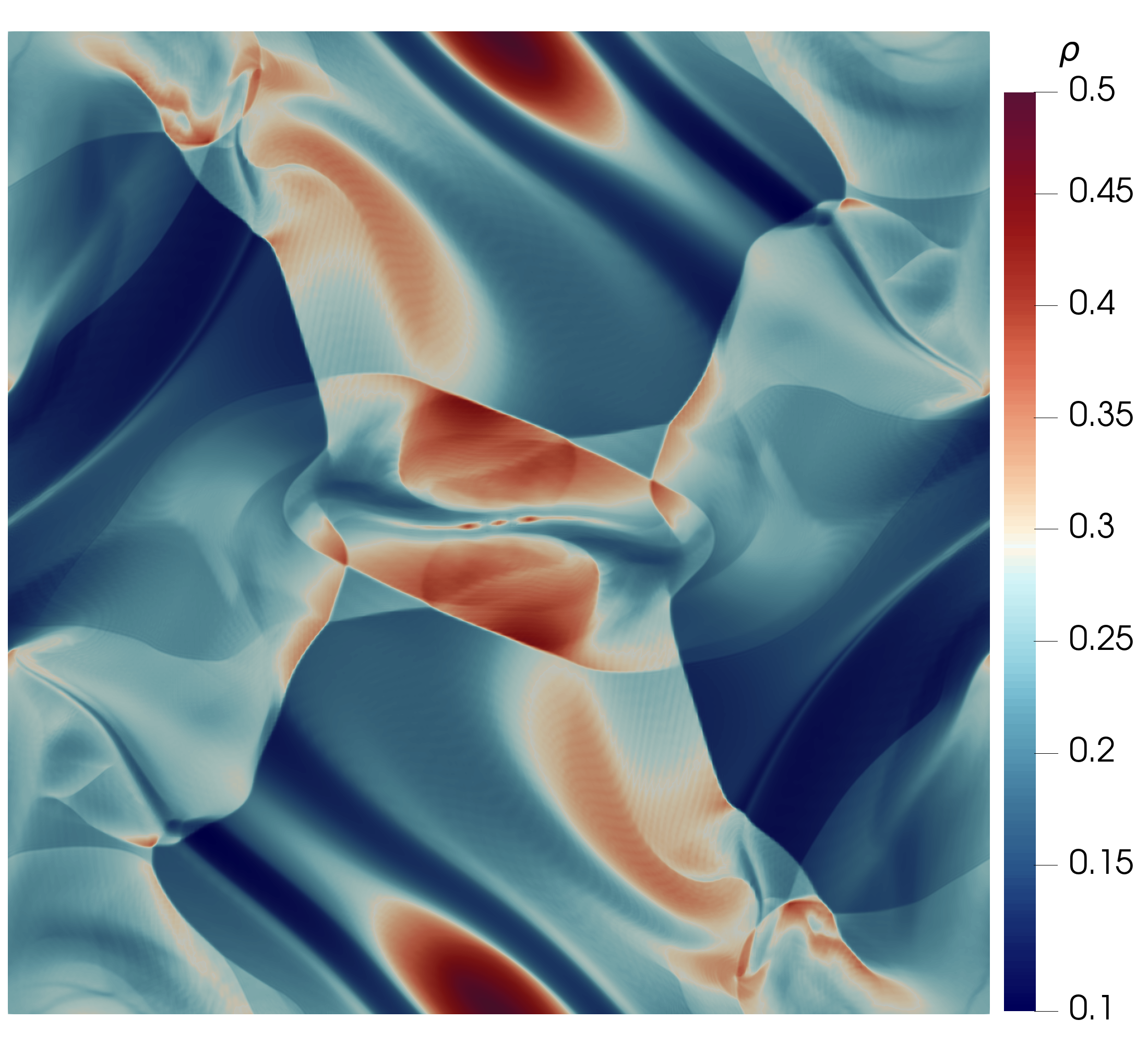} &
\includegraphics[width=0.45\columnwidth]{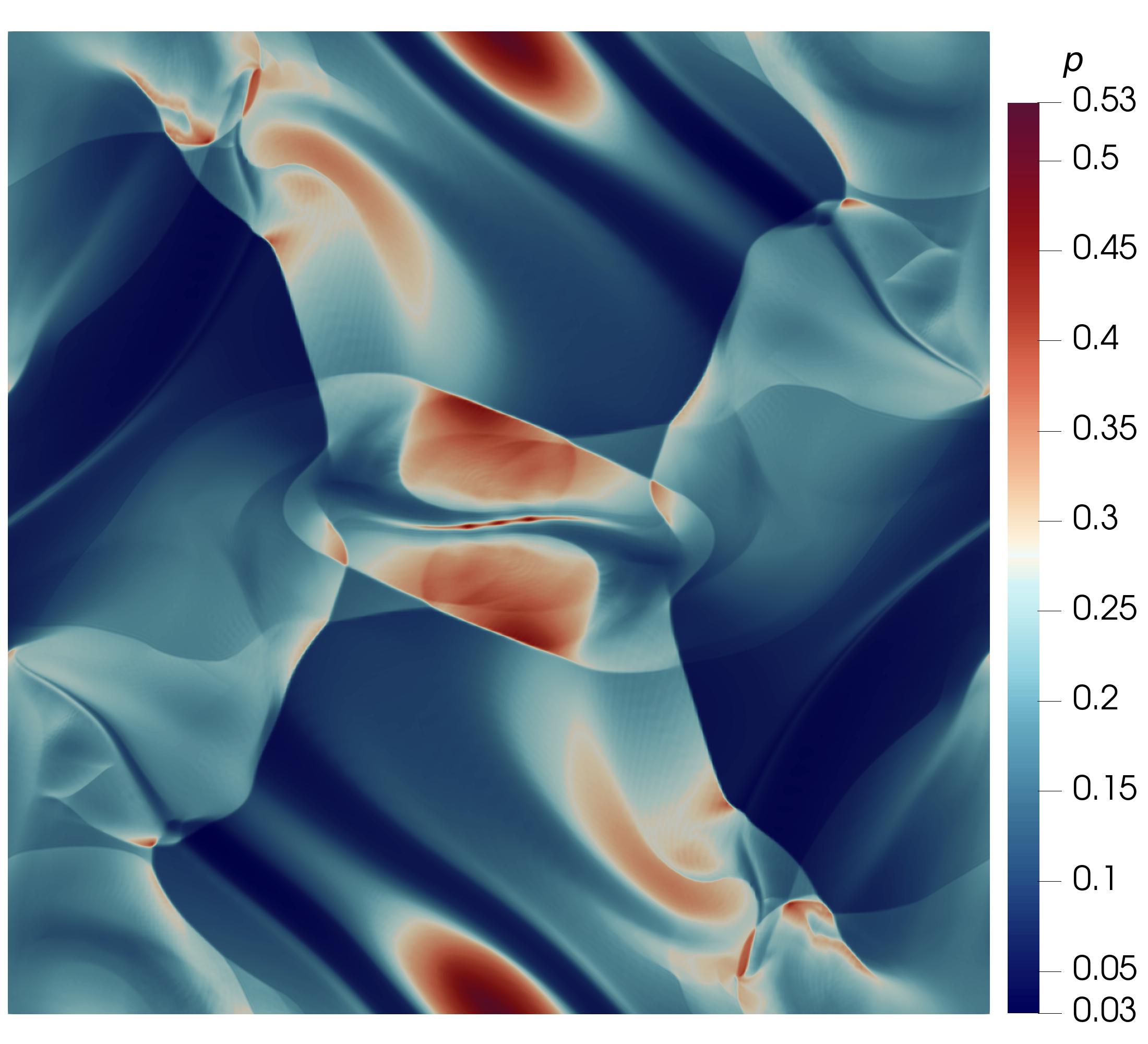} \\
(a) & (b)
\end{tabular}
\caption{Orszag-Tang vortex~\eqref{eq:orszag.tang.ic} test (a) density and (b) pressure at time $t=0.5$ using $N=4$ and $1025 \times 1025$ degrees of freedom.}
\label{fig:tang.vortex}
\end{figure}

\begin{figure}
\centering
\begin{tabular}{cc}
\includegraphics[width=0.45\columnwidth]{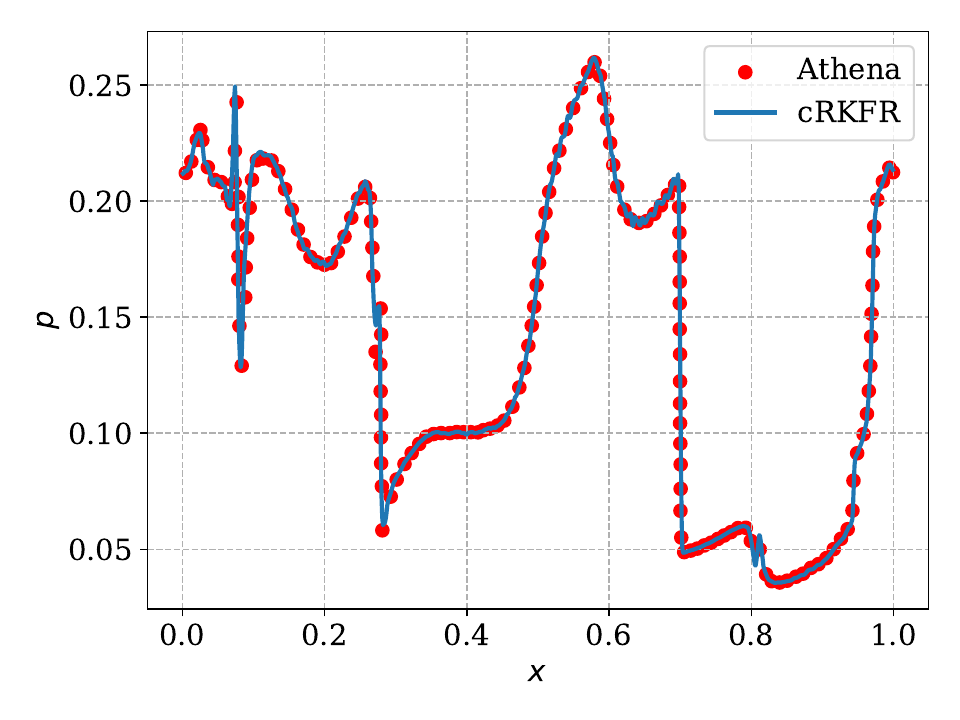} &
\includegraphics[width=0.45\columnwidth]{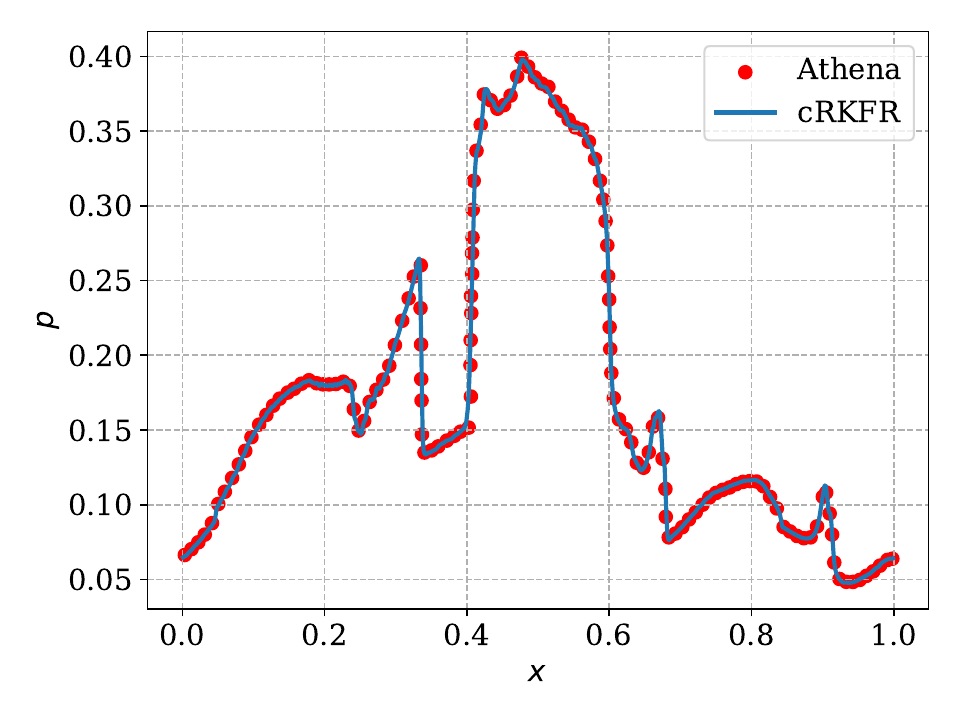} \\
(a) $y=0.3125$ & (b) $y=0.4277$
\end{tabular}
\caption{Pressure line cuts at $y=0.3125$ and $y=0.4277$ for the Orszag-Tang vortex~\eqref{eq:orszag.tang.ic} test at time $t=0.5$ using $N=4$ and $1025 \times 1025$ degrees of freedom.}
\label{fig:tang.line.cuts}
\end{figure}

The Orszag-Tang vortex problem~\cite{orszag1979,derigs2016} is a widely used benchmark test for MHD codes, and is also mentioned to be suitable for MHD turbulence in~\cite{derigs2018}.
It includes dissipation of kinetic and magnetic energy, magnetic reconnection, the formation of high-density jets, dynamic alignment and the emergence and manifestation of small-scale structures~\cite{derigs2018}.
The initial condition consists of smooth non-random data given by
\begin{equation}
\begin{gathered}
\rho = 1, \quad \bv = (- \sin(2 \pi y), \sin(2 \pi x), 0), \quad
p = \gamma^{-1}, \quad \bB = (-\gamma^{-1} \sin 2 \pi y, \gamma^{-1} \sin 4 \pi x, 0),
\end{gathered}
\label{eq:orszag.tang.ic}
\end{equation}
with $\gamma = 5/3$ and periodic boundary conditions on the domain $[0,1]^2$.
The flow evolves gradually becoming increasingly complex with the formation of shocks, demonstrating compressible MHD turbulence~\cite{derigs2018}.
The simulation is run till time $t=0.5$ using polynomial degree $N=4$ with $1025 \times 1025$ degrees of freedom ($205 \times 205$ elements), chosen to compare with reference data from~\cite{ramirez2021} generated using Athena~\cite{stone2008} with $1024 \times 1024$ degrees of freedom.
In Figure~\ref{fig:tang.vortex}, the density and pressure contours are shown, while in Figure~\ref{fig:tang.line.cuts}, line cuts of pressure at $y=0.3125$ and $y=0.4277$ are shown to compare with the reference solution generated using Athena with a resolution of $1024 \times 1024$ degrees of freedom.
A reasonable agreement similar to~\cite{ramirez2021} is observed, validating our non-conservative product treatment.

%------------------------------------------------------------------------------
\subsubsection{Rotor test}
\begin{figure}
\centering
\begin{tabular}{ccc}
\includegraphics[width=0.3\columnwidth]{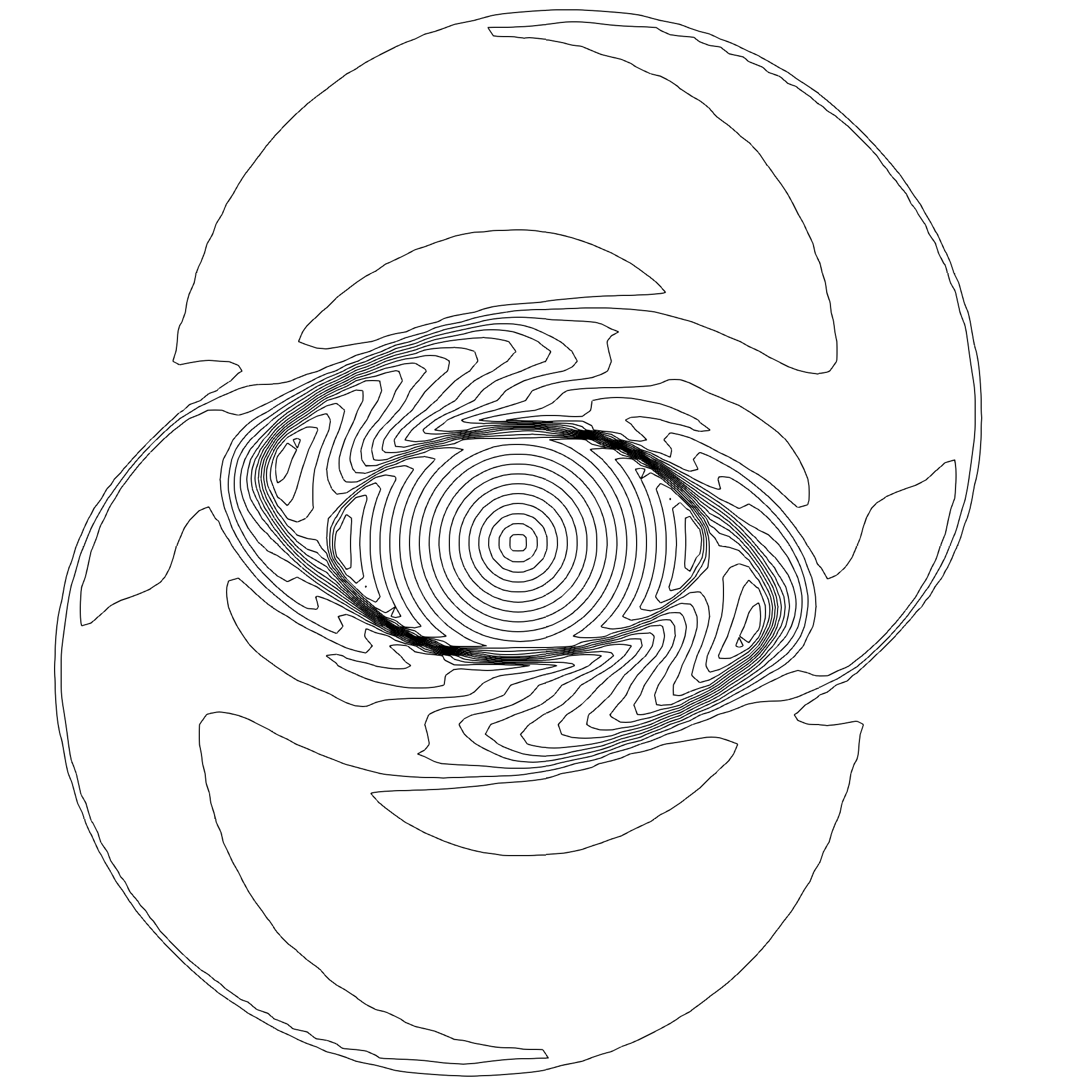} &
\includegraphics[width=0.3\columnwidth]{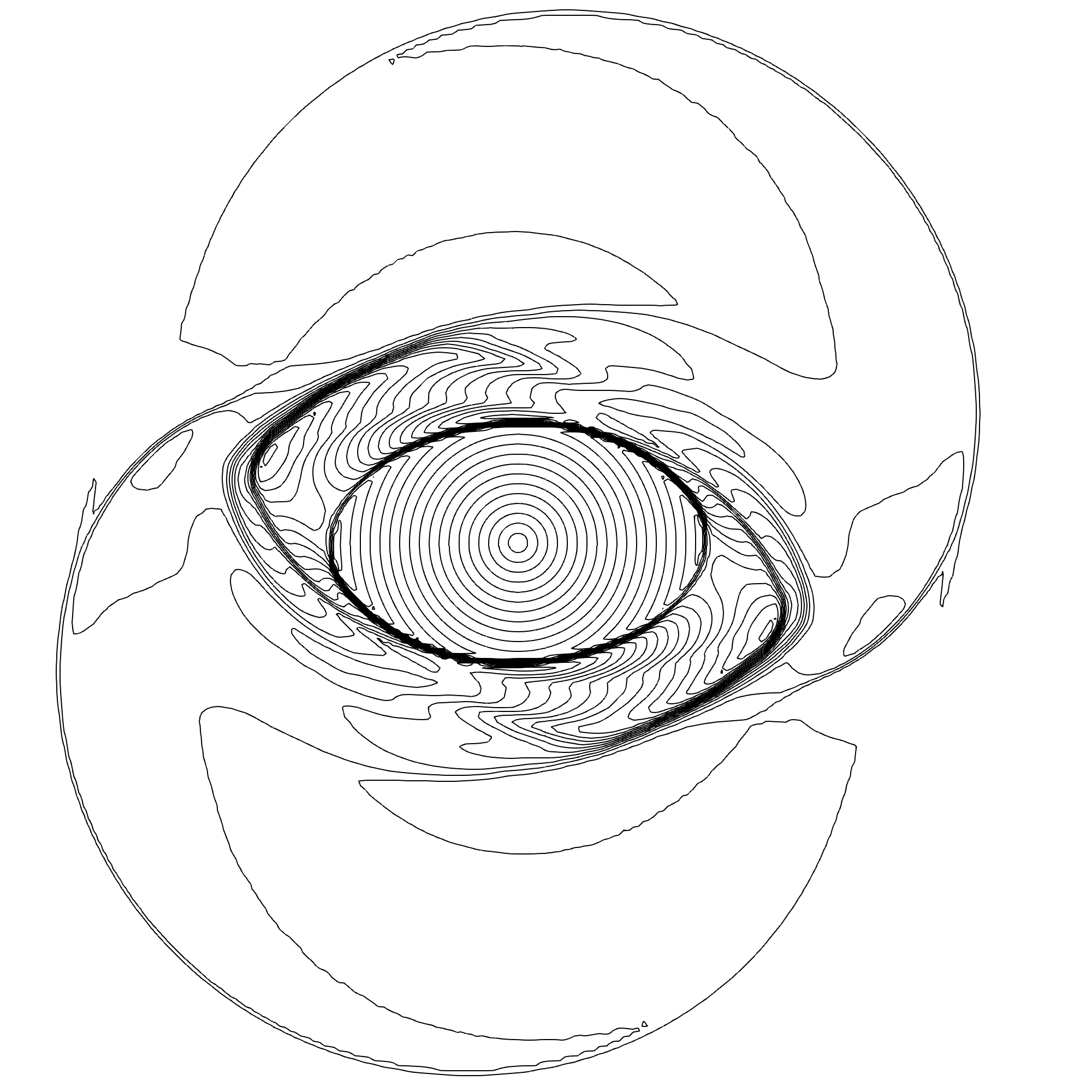} &
\includegraphics[width=0.3\columnwidth]{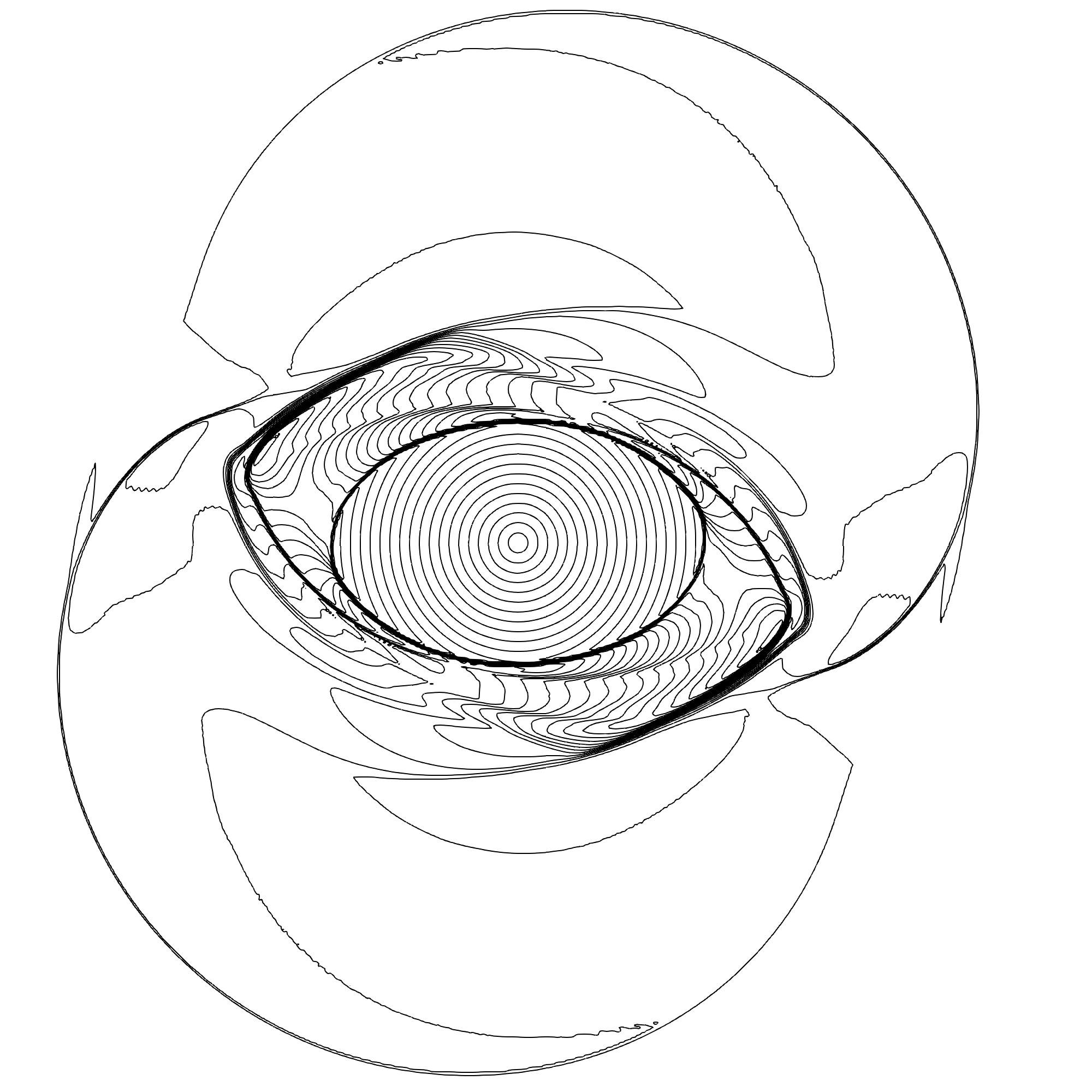} \\
(a) $128^2$ & (b) $256^2$ & (c) $512^2$
\end{tabular}
\caption{Rotor test~\eqref{eq:rotor.ic} Mach number contours at time $t=0.15$ using polynomial degree $N=4$ for grid resolutions (a) $128^2$, (b) $256^2$ and (c) $512^2$.}
\label{fig:rotor.mach}
\end{figure}

The MHD rotor test from~\cite{balsara1999} consists of a rapidly spinning dense cylinder embedded in a magnetized, homogeneous medium at rest~\cite{derigs2018}.
In this test, centrifugal forces cause the dense cylinder to not be in equilibrium.
The rotor spins with the initially prescribed rotating velocity, causing the initially uniform magnetic field to wind up.
This wrapping by the magnetic field leads to strong toroidal Alfvén waves launched into the ambient fluid.
The initial conditions used are given by
\begin{equation}
\begin{gathered}
\rho, v_1, v_2 =
\begin{cases}
10, -20 (y - 0.5), 20 (x - 0.5) \quad & r \le r_0,\\
1, 0, 0, & r > r_1,\\
1 + 9f(r), -20f(r)(y-0.5), 20f(r)(x-0.5), & r_0 < r \le r_1,
\end{cases}
v_3 = 0, \quad p = 1, \\
B_1 = \frac{5}{\sqrt{4 \pi}}, \quad B_2 = 0, \quad B_3 = 0,
\end{gathered}
\label{eq:rotor.ic}
\end{equation}
where
\[
r = \sqrt{(x-0.5)^2 + (y-0.5)^2}, \quad f(r) = \frac{r_1 - r}{r_1 - r_0}, \quad
r_0 = 0.1, \quad r_1 = 0.115.
\]
Following~\cite{chandrashekar2020}, we validate our scheme by plotting the Mach number contours at time $t=0.15$ for three different grid resolutions $128^2, 256^2$ and $512^2$ in Figure~\ref{fig:rotor.mach} and polynomial degree $N=4$.
As mentioned in~\cite{fengyan2005}, a scheme that does not sufficiently control the divergence error of the magnetic field will show distortions in the Mach number.
In Figure~\ref{fig:rotor.mach}, we see no such distortions at any resolution, indicating that our non-conservative scheme is working fine for this test case.

%------------------------------------------------------------------------------
\subsection{Multi-ion MHD} \label{sec:multi.ion.mhd}
%------------------------------------------------------------------------------

As an equation consisting of several non-conservative terms, we solve the multi-ion magnetohydrodynamics (MHD) equations, first introduced in~\cite{toth2010multi} with recent development using flux differencing DG methods in~\cite{ramirez2025}.
These equations describe the motion of multi-ion plasmas with independent momentum and energy equations from each ion species.
The multi-ion MHD equations are used in their quasi-conservative form, as given in equation (18) of~\cite{ramirez2025}
\begin{equation}
\begin{gathered}
  \label{eq:multi-ion_mod_div}
    % &
    \frac{\partial}{\partial t}
    \begin{pmatrix}
    \rho_k \\ \rho_k \vecr{v}_k  \\ E_{k} \\ \vecr{B} \\
    % \red{
      \psi
      % }
    \end{pmatrix}
    + \Nabla \cdot
    % \underbrace{
    \begin{pmatrix}
    \rho_k \vecr{v}_k  \\
    \rho_k \vecr{v}_k\, \vecr{v}_k^{\,T} +
     p_k \threeMatrix{I}   \\
    \vecr{v}_k\left(\hydroEner_k + p_k \right)  \\
    \bzero \\
    \vecr{0}
    \end{pmatrix}
    % }_{\blocktensor{f}^{\supEuler}}
    +
    % MHD conservative term!!
    %%%%%%%%%%%%%%%%%%%%%%%%%
    \Nabla \cdot
    % \underbrace{
    \begin{pmatrix}
    \vecr{0}  \\
    % Momentum equation
    \bzero \\
    % Energy equation
    \vecr{v}^+_k\,\|\vecr{B}\|^2 - \vecr{B}\left(\vecr{v}^+_k\cdot\vecr{B}\right)\\
    \vecr{v}^+ \vecr{B}^T - \vecr{B} (\vecr{v}^+)^T \\
    \vecr{0}
    \end{pmatrix}
    % }_{\blocktensor{f}^{\supMHD}}
    +
    % \red{
    \Nabla \cdot
    % \underbrace{
    \begin{pmatrix}
    \vecr{0} \\
    \bzero \\
    c_h \psi \vecr{B} \\
    c_h \psi \threeMatrix{I} \\ c_h \vecr{B}\\
    \end{pmatrix}
    % }_{\blocktensor{f}^{\supGLM}}
    % }
    \\
    % &
    +
    % Source term
    %%%%%%%%%%%%%
    % \underbrace{
    \begin{pmatrix}
    0 \\
    r_k \rho_k (\vecr{v}^+ - \vecr{v}_k) \times \vecr{B} \\
    r_k \rho_k \vecr{v}_k \cdot (\vecr{v}^+ - \vecr{v}_k) \times \vecr{B}  \\
    \vecr{0} \\
    0
    \end{pmatrix}
    % }_{\state{g}}
    +
    % Powell non-conservative term
    %%%%%%%%%%%%%%%%%%%%%%%%%%%%%%
    % \underbrace{
    \begin{pmatrix}
    0 \\
    \frac{r_k \rho_k}{n_e e} \vecr{B} \\
    \vecr{v}^+ \cdot \vecr{B} \\
    \vecr{v}^+ \\
    0
    \end{pmatrix}
    (\Nabla \cdot \vecr{B})
    % }_{\noncon^{\Powell} \coloneqq \stateG{\phi}^{\Powell} (\Nabla \cdot \vecr{B}) }
    +
    % Multi-ion non-conservative term
    %%%%%%%%%%%%%%%%%%%%%%%%%%%%%%%%%
    \begin{pmatrix}
    0 \\
    \frac{r_k \rho_k}{n_e e} \Nabla \cdot \left(\half \norm{\vecr{B}}^2 \mat{I} - \vecr{B}\vecr{B}^T
    + p_e \mat{I}\right) \\
    \vecr{v}_k^+ \cdot \Nabla p_e
    + \vecr{B} \cdot
    \Nabla \cdot \left(\vecr{v}_k^- \vecr{B}^T - \vecr{B} (\vecr{v}_k^-)^T  \right) \\
    \vecr{0} \\
    0
    \end{pmatrix} \\
    % &
    +
    % \red{
    % \underbrace{
    \begin{pmatrix}
        \vecr{0} \\
        \bzero\\
        \vecr{v}^+ \psi \\
        \bzero \\
        \vecr{v}^+
    \end{pmatrix}
    \cdot
    \Nabla \psi
    % }_{\noncon^{\supGLM}}
    % }
    =
    \state{0},
\end{gathered}
\end{equation}
where $\rho_k, \vecr{v}_k = (v_{k,1}, v_{k,2}, v_{k,3}), \hydroEner_k$ are the mass density, velocity and total energy of ion species $k$, respectively.
The $\vecr{B} = (B_1, B_2, B_3)$ is the magnetic field, $\psi$ is the divergence-cleaning scalar field~\cite{derigs2018,derigs2016}, $r_k$ is the charge-to-mass ratio given by $q_k / m_k$ where $n_k$ denotes the number density and $q_k$ the charge of ion species $k$.
The source terms in the above equations represent the Lorentz force, an integral part of the multi-ion MHD equations.
Similar to the single ion MHD case~\eqref{eq:glm.mhd}, $c_h$ is the hyperbolic divergence cleaning speed.
The $e$ denotes the elementary charge and $n_e$ denotes the number density of electrons given by assuming quasi-neutrality $n_e = \frac{1}{e} \sum_k n_k q_k$.
The electron pressure $p_e$ is estimated with a phenomenological model.
It can be computed as a fraction of the total ion pressure as $p_e = \alpha \sum_k p_k$, where $\alpha$ is heuristically determined depending on the problem of interest.
The $\vecr{v^+}$ denotes the charge-averaged ion velocity given by $\vecr{v^+} = \sum_k \frac{n_k q_k}{n_e e} \vecr{v_k}$.
Calorically perfect plasmas are assumed for each ion species, and the total energy $E_k$ includes the GLM variable $\psi$ in order to avoid erroneous variation of the hydrodynamic energy $\hydroEner_k$~\cite{ramirez2025}
\[
E_k = \hydroEner_k + \frac 12 \norm{\vecr{B}}^2 + \frac 12 \psi^2, \qquad
\hydroEner_k = \frac{p_k}{\gamma - 1} + \frac 12 \rho_k \norm{\vecr{v_k}}^2.
\]
The fluxes and non-conservative products for the multi-ion MHD equations~\eqref{eq:multi-ion_mod_div} are incorporated into \texttt{Tenkai.jl}~\cite{tenkai} using \texttt{Trixi.jl}~\cite{schlottkelakemper2021purely,schlottkelakemper2020trixi,ranocha2022adaptive} through Julia's package manager.
\subsubsection{Convergence test}

\begin{figure}
  \centering
  \begin{tabular}{cc}
    \includegraphics[width=0.38\columnwidth]{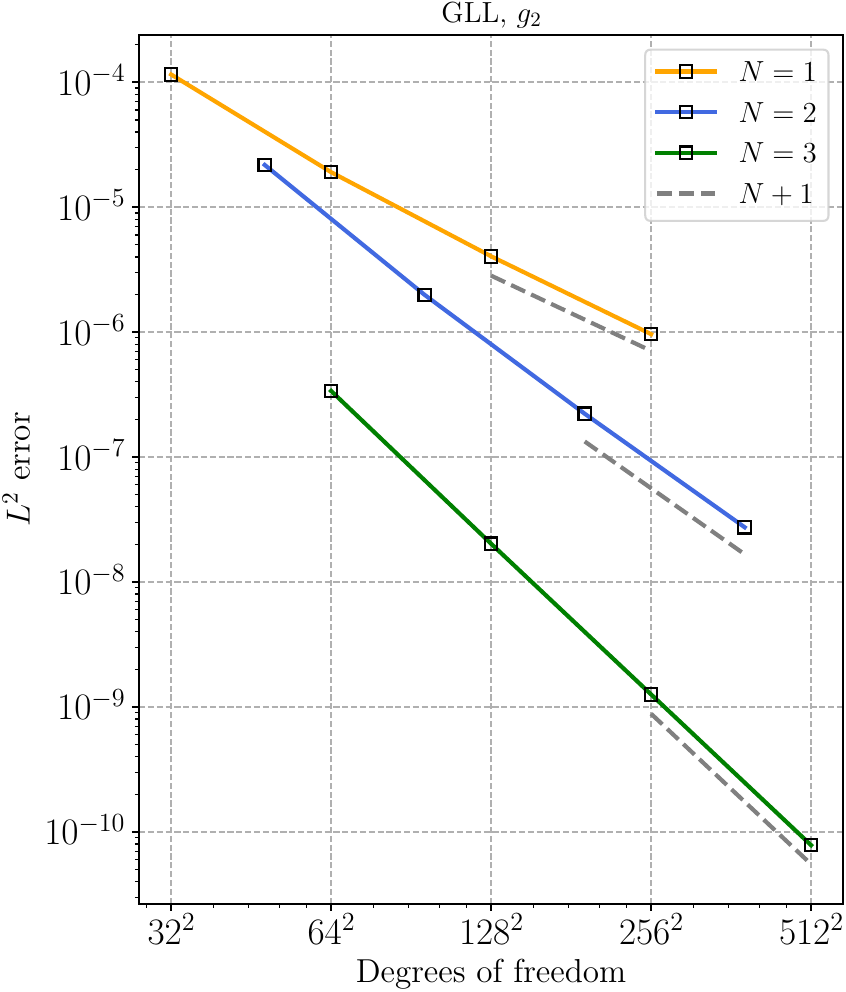} &
    \includegraphics[width=0.38\columnwidth]{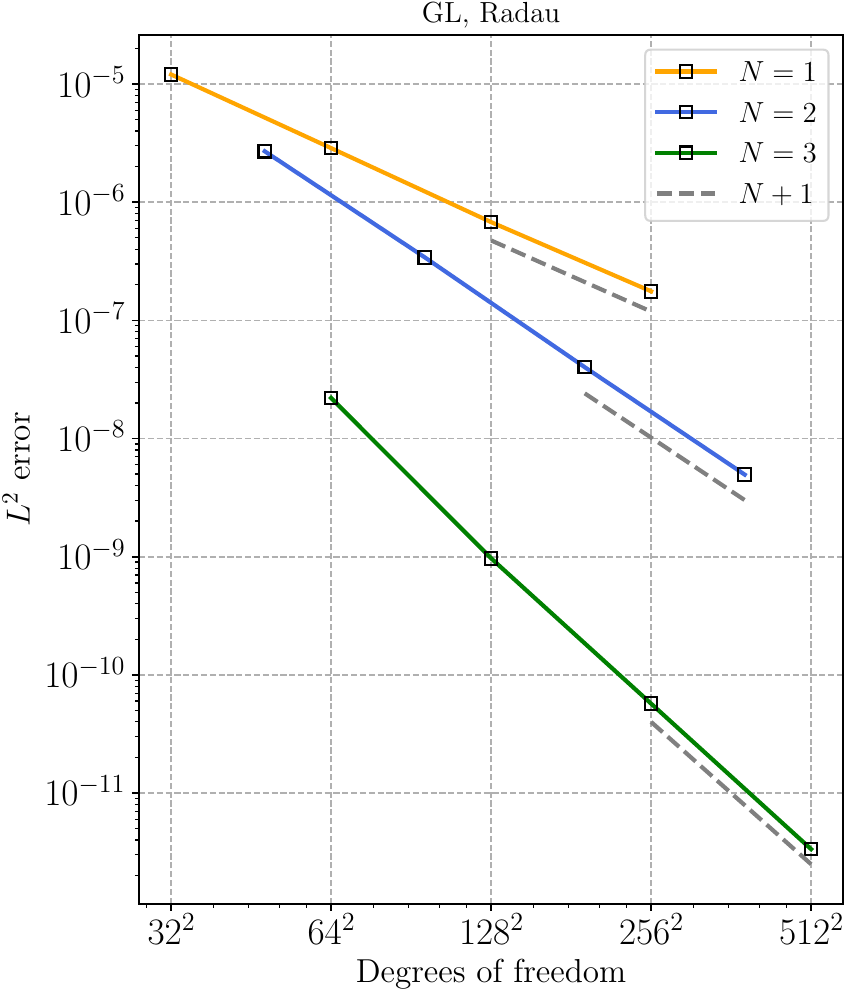} \\
    (a) & (b)
  \end{tabular}
  \caption{Multi-ion convergence test using (a) $g_2$ correction and GLL points, (b) Radau correction and GL points.}
  \label{fig:multi.ion.convergence}
\end{figure}

We test the order of accuracy of the scheme using the manufactured solution test from~\cite{toth2010multi,ramirez2025} which relies on choosing a smooth exact solution and modifying the source term appropriately.
The expressions for the manufactured solution are not provided to save space, but can be found in equation (91) of~\cite{ramirez2025} or in our reproducibility repository~\cite{babbar2025crknonconsrepro}.
We would like to highlight that the manufactured solution consists of terms for all the variables ensuring that the non-conservative terms are well tested.
% \begin{align}\label{eq:mansol}
%     \rho_1 (\bx, t) &= \chi_1,
%     &
%     \rho_2  (\bx, t) &= \chi_2, \nonumber\\
%     \rho_1 v_{1,1} (\bx, t) &= \chi_1, &
%     \rho_2 v_{2,1} (\bx, t) &= \chi_2,\nonumber\\
%     \rho_1 v_{1,2} (\bx, t) &= \chi_1,&
%     \rho_2 v_{2,2} (\bx, t) &= \chi_2,\nonumber\\
%     \rho_1 v_{1,3} (\bx, t) &= 0.1 \chi_1,&
%     \rho_2 v_{2,3} (\bx, t) &= 0.1 \chi_2,\nonumber\\
%     E_1 (\bx, t) &= 2 \chi_1^2+\chi_1, &
%     E_2 (\bx, t) &= 2 \chi_2^2+\chi_2,\nonumber\\
%     B_1 (\bx, t) &= 0.25 \chi, &
%     B_2 (\bx, t) &= -0.25 \chi,\nonumber\\
%     B_3 (\bx, t) &= 0.1 \chi, &
%     \psi(\bx, t) &=0,
% \end{align}
% where
% \begin{align*}
%     \chi &= \underbrace{0.1 \sin(\pi (x + y - t))}_{:= \chi_0} + 2,\\
%     \chi_1 &= 0.04 \sin(\pi (x + y - t)) + 1,\\
%     \chi_2 &= \chi-\chi_1.
% \end{align*}
Further, the heat capacity ratios for the two ion species are $\gamma_1 = 2$ and $\gamma_2 = 4$, and the charge-to-mass ratios are $r_1 = 2$ and $r_2 = 1$.
The electric pressure ratio is computed as $p_e = 0.2 (p_1 + p_2)$.
The physical domain is $[-1,1]^2$ with periodic boundary conditions and the final time is $t=1$.
The convergence results are shown in Figure~\ref{fig:multi.ion.convergence} for degrees $N=1,2,3$, showing the optimal order of accuracy for all degrees.

\subsubsection{Kelvin-Helmholtz instability}

\begin{figure}
\centering
\begin{tabular}{ccccc}
\includegraphics[width=0.053\columnwidth]{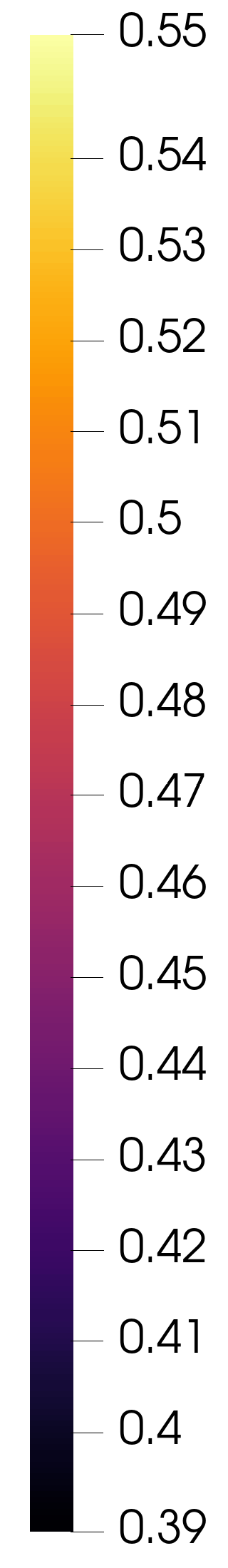} &
\includegraphics[width=0.20\columnwidth]{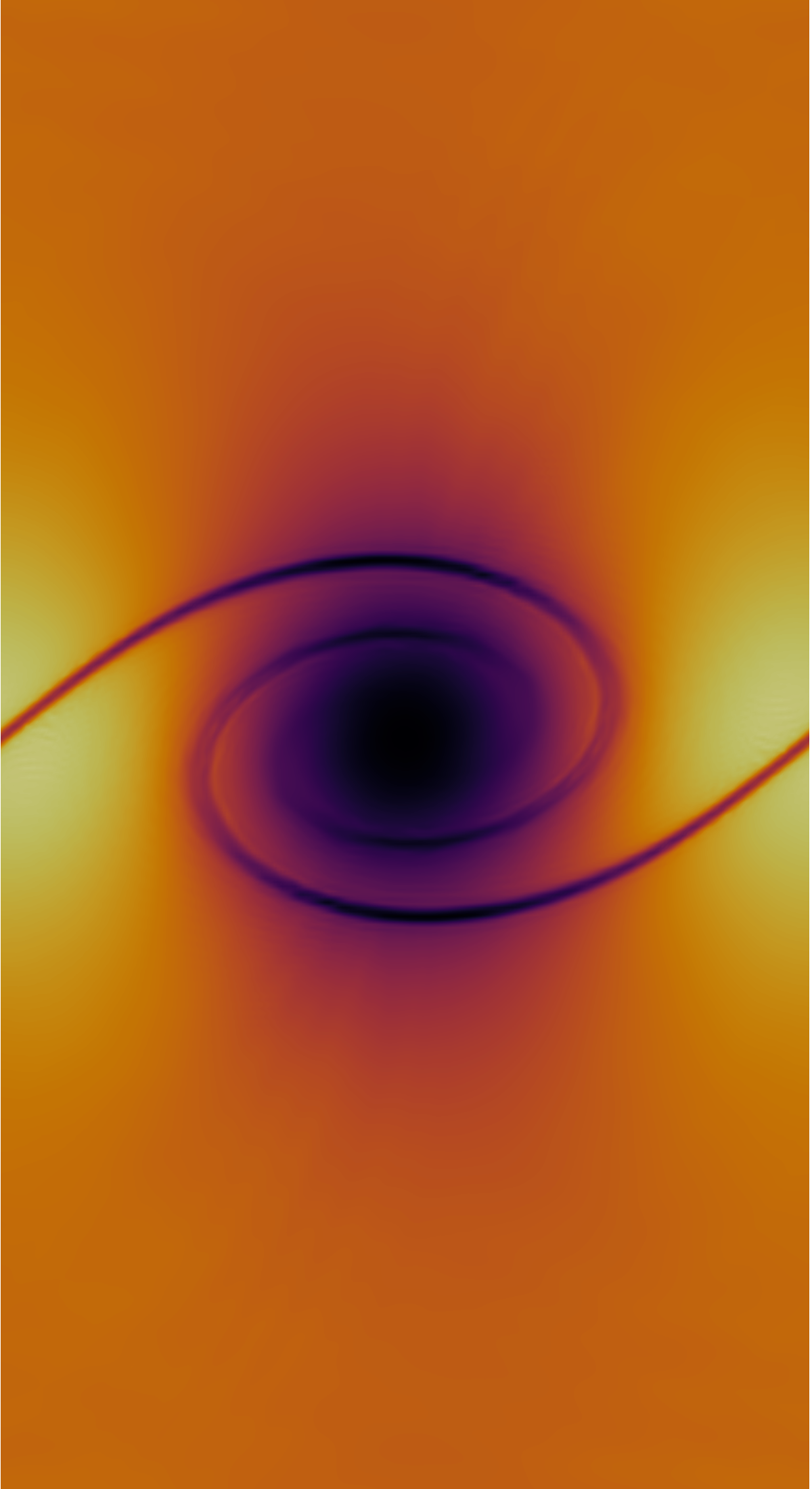} &
\includegraphics[width=0.20\columnwidth]{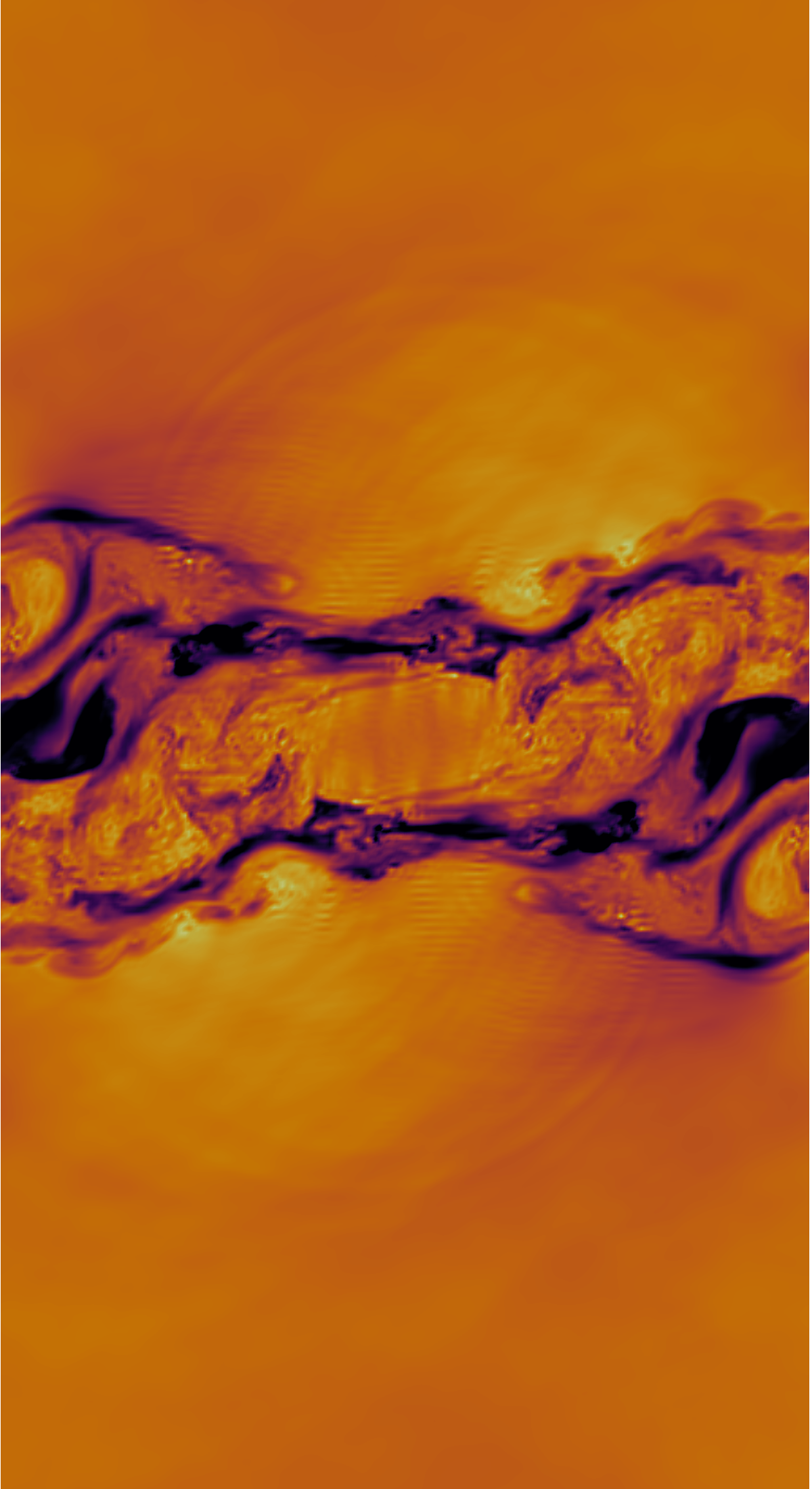} &
\includegraphics[width=0.20\columnwidth]{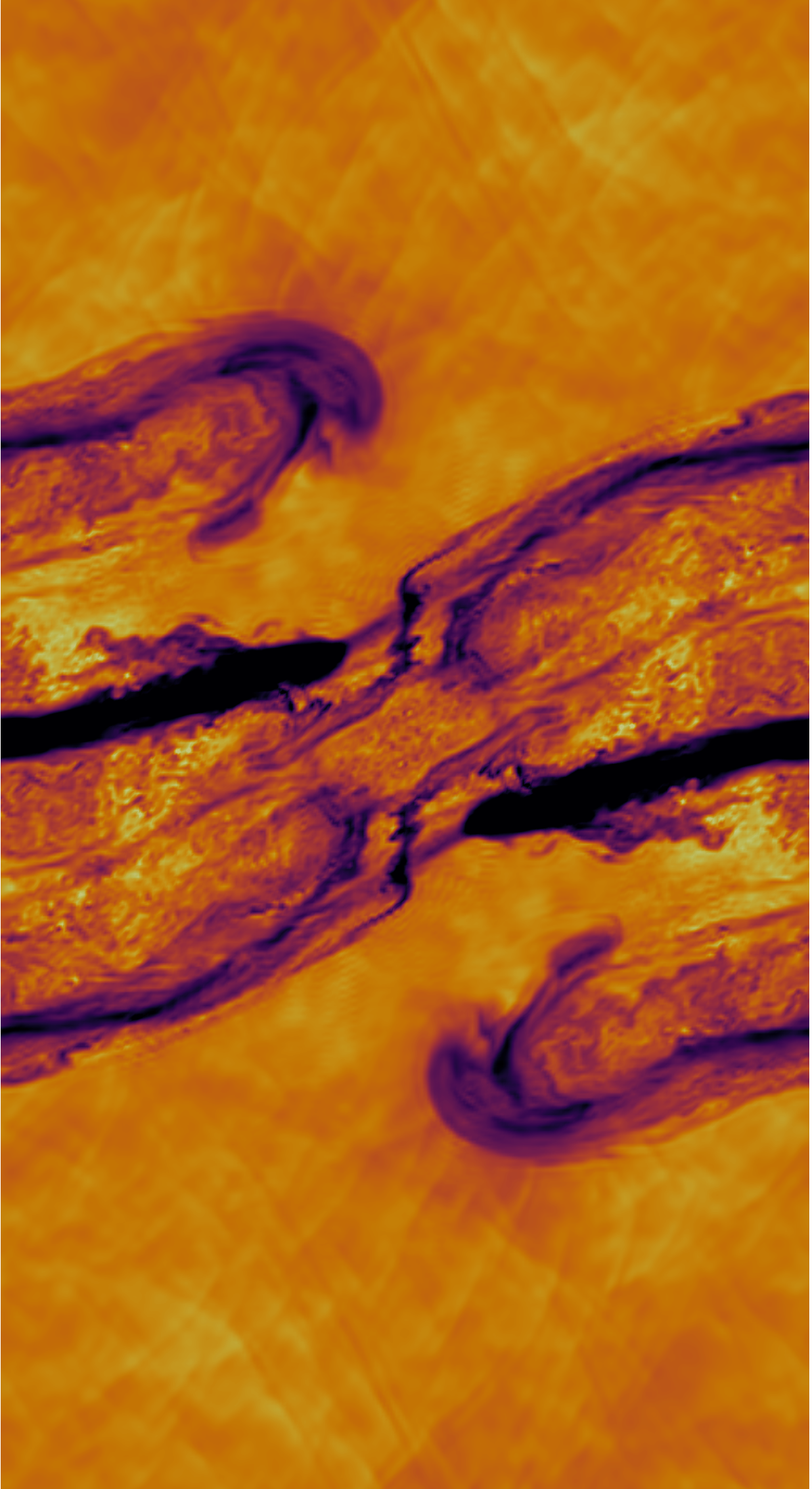} &
\includegraphics[width=0.20\columnwidth]{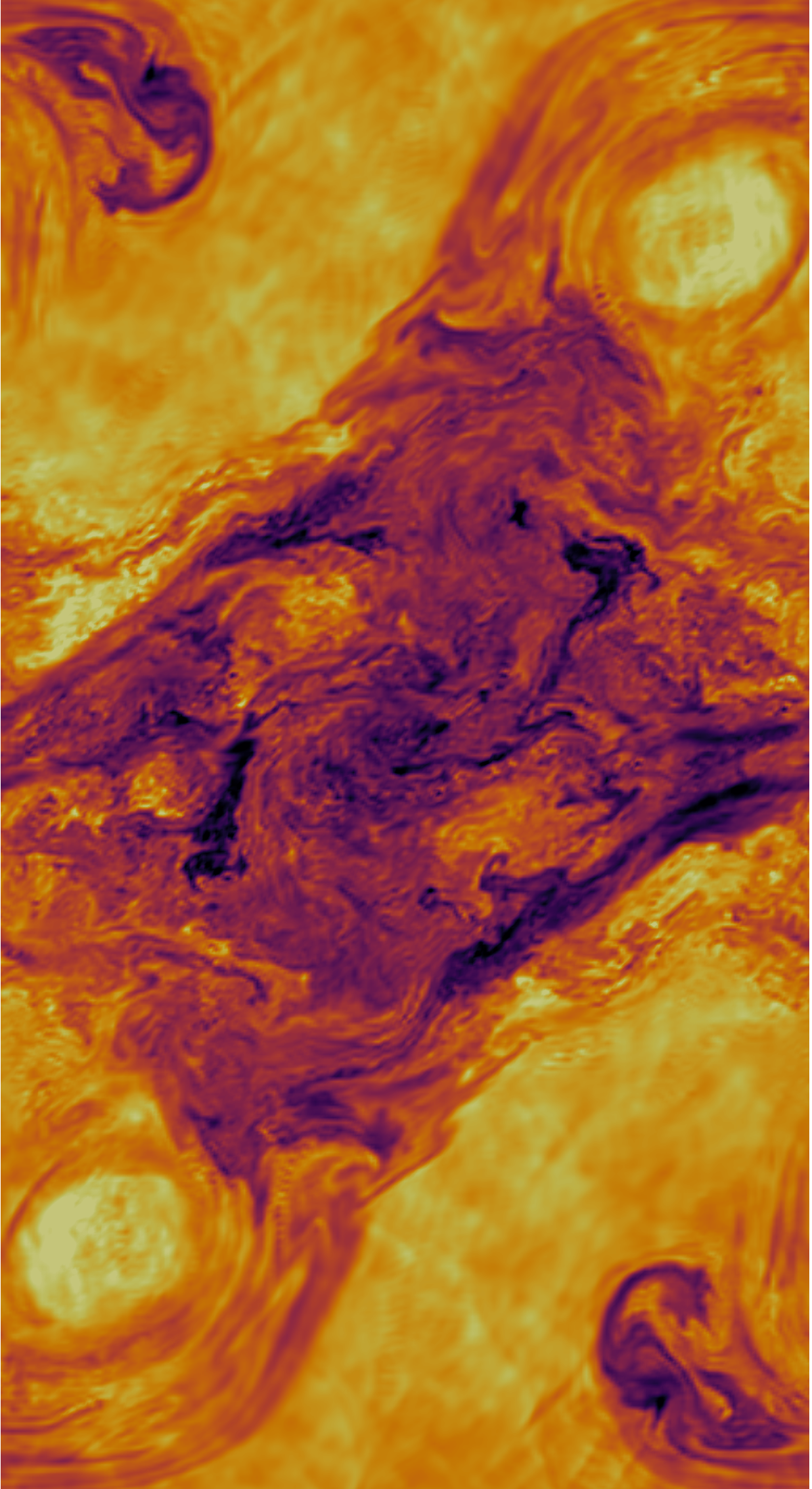} \\
& $t=5$ & $t=8$ & $t=12$ & $t=20$
\end{tabular}
\caption{Multi-ion MHD Kelvin-Helmholtz instability, density of ion species 1 at different times using polynomial degree $N=3$ and $64 \times 128$ elements.}
\label{fig:khi.multi.ion}
\end{figure}

We now consider the two ion Kelvin-Helmholtz instability tests from~\cite{ramirez2025}.
Two ion species are considered with heat capacity ratios $\gamma_1 = 5/3$ (monoatomic) and $\gamma_2 = 1.4$ (diatomic), charge to mass ratios $r_1 = 1$ and $r_2 = 0.5$.
The collision source terms between the ion species are neglected, and a trivial electron pressure gradient $\Nabla p_e = 0$ is taken.
As in the previous test, the reader is referred to equation (98) of~\cite{ramirez2025} or to our reproducibility repository~\cite{babbar2025crknonconsrepro} for the complete expressions of the initial condition.
% The initial condition is given by
% \begin{align}\label{eq:khi_condini}
% \rho_k(x,y,t=0) &= \half, &p_k(x,y,t=0) &= \frac{1}{\gamma_k}, &\psi(x,y,t=0) &= 0, \nonumber\\
% v_{k,1} (x,y,t=0) &= \half \tanh \left( \frac{y}{y_0} \right),
% &v_{k,2}(x,y,t=0) &= v_{2,0} \sin (2 \pi x) \exp \left(- \frac{y^2}{\sigma^2} \right), &
% v_{k,3}(x,y,t=0) &= 0, \nonumber \\
% B_1 (x,y,t=0) &= c_a \cos \theta,
% &B_2(x,y,t=0) &= 0, &B_3(x,y,t=0) &= c_a \sin \theta,
% \end{align}
% with $k=1,2$, $y_0=1/20$ is the steepness of the shear, $c_a=0.1$ is the Alfvén speed, and $\theta = \pi/3$ is the angle of the initial magnetic field.
% Moreover, the parameters of the perturbation are $v_{2,0}=0.01$ and $\sigma=0.1$.
The physical domain is $[0,1] \times [-1,1]$ with periodic boundary conditions in the $x$-direction and solid wall boundary conditions in the $y$-direction.
The simulation is run till $t=20$ using polynomial degree $N=3$ and $64 \times 128$ elements. The density of ion species 1 at different times is shown in Figure~\ref{fig:khi.multi.ion}.
The solution at time $t=5$ is in good agreement with the results shown in~\cite{ramirez2025}, thus validating our non-conservative scheme for the multi-ion MHD equations.
We have run the simulation with blending coefficient restricted to $\alpha_{\max} = 0.3$, although similar results are obtained with other values of $\alpha_{\max}$.
However, the simulation crashes when no blending scheme or when $\alpha_{\max} \le 0.1$ is used.
In~\cite{ramirez2025}, the authors use an entropy stable flux difference scheme and are able to run the simulation without any blending.
Thus, development of entropy stable schemes for cRKFR schemes on non-conservative systems is an important area of research.
The indicating variable~\cite{hennemann2021,babbar2024admissibility} used to compute the blending coefficient is the product of sum of densities and pressures $(\sum_k \rho_k)(\sum_k p_k)$ of all species $k$, which is the natural extension of the indicator used for compressible Euler equations in~\cite{hennemann2021,babbar2024admissibility}.

\section{Conclusions} \label{sec:conclusions}

We extended the compact Runge-Kutta flux reconstruction method of~\cite{babbar2025crk} in a time averaged framework to handle hyperbolic systems with non-conservative products and stiff source terms.
The stiff source terms are handled by using an IMplicit EXplicit (IMEX) time integration discretization, where the non-stiff advection terms are treated explicitly, while the stiff source terms are treated by a locally implicit approach.
Only one inter-element numerical flux computation per time step is required, which is a feature retained from the conservative case in~\cite{babbar2025crk}.
Although the hyperbolic system consists of advection terms that cannot be written as the divergence of a flux, inter-element communication is still performed by using \textit{discontinuous numerical fluxes} at the element interfaces.
Similar numerical fluxes are then used to develop a first order finite volume scheme on subcells within each element, which is then used to develop a subcell based blending limiter to handle non-smooth solutions of hyperbolic equations.
Using the flux limiter of~\cite{babbar2024admissibility} extended to non-conservative systems, and additional limiting with the subcell based lower order evolution, an admissibility-preserving scheme is obtained.
Numerical results for a variety of non-conservative hyperbolic systems with and without stiff source terms have been shown to demonstrate the performance of the scheme.
The method and code are developed in such a way that the user only needs to provide the system of equations and an estimate of the wave speed for the system to use the scheme.
Some open problems that remain to be addressed are the development of path conservative methods within this framework, and the development of well-balanced cRK methods.

\section*{Acknowledgments}

AB and HR were supported by the Deutsche Forschungsgemeinschaft
(DFG, German Research Foundation, project number 528753982
as well as within the DFG priority program SPP~2410 with project number 526031774).
AB was also supported by the Alexander von Humboldt Foundation.

We thank Praveen Chandrashekar for discussions on, and assistance with the implementation of, the shear shallow water equations, Marco Artiano for adding Jin-Xin relaxation to Trixi.jl, Andrés M.\ Rueda-Ramírez for making his multi-ion MHD code publicly available, and Qifan Chen for the constructive feedback on an earlier version of the paper.

\appendix
\section{Stiff source term analysis} \label{app:stiff.source}

We follow~\cite{svard2011} to motivate the application of IMEX model using the simple model
problem consisting of an advection equation with a linear stiff source term
\begin{equation}
\begin{aligned}
u_t + au_x & = - Ku,\\
u (x, 0) & = u_0 (x),
\end{aligned} \label{eq:advection.stiff.source}
\end{equation}
where $K$ is a large constant.
The exact solution is
\[ u (x, t) = u_0 (x - t) \exp (- Kt), \]
which is of a decaying nature.
Assuming $a > 0$, the first order explicit
method with an upwind discretization is given by
\[ u_j^{n + 1} = u_j^n - a \Delta t \frac{u_j^n - u_{j - 1}^n}{\Delta x} -
   K \Delta t u_j^n, \]
while one with the implicit discretization for the source terms is given by
\[ u_j^{n + 1} = u_j^n - a \Delta t \frac{u_j^n - u_{j -
   1}^n}{\Delta x} - K \Delta t u_j^{n + 1} .
\]
To perform a von Neumann analysis, we substitute a solution consisting of a
single Fourier node
\[ u_j^n = \hat{u}^n \exp (ij \xi), \quad \xi \in [- \pi, \pi], \qquad i =
   \sqrt{- 1} .
\]
Thus, the evolution of $u_j^{n + 1} = \hat{u}^{n + 1} \exp (ij \xi)$ for the
explicit scheme is given by

\begin{equation}
  \begin{aligned}
    \hat{u}^{n + 1} & = & \left( 1 - a \frac{\Delta t}{\Delta x}
    \bigl( 1 - \exp (- i \xi) \bigr) \right)  \hat{u}^n - K \Delta t \hat{u}^n & = &
    \left( 1 - a \frac{\Delta t}{\Delta x} \bigl( 1 - \exp (- i \xi) \bigr) - K
    \Delta t \right)  \hat{u}^n,
  \end{aligned} \label{eq:expl.evolution}
\end{equation}
while the evolution for the IMEX scheme is given by
\begin{equation}
  \begin{aligned}
    \hat{u}^{n + 1} & = \left( 1 - a \frac{\Delta t}{\Delta x} \bigl( 1 - \exp (- i \xi) \bigr) \right)  \hat{u}^n - K \Delta t \hat{u}^{n + 1}\\
    \Rightarrow \hat{u}^{n + 1} & = \frac{1 - a \frac{\Delta
    t}{\Delta x} \bigl( 1 - \exp (- i \xi) \bigr)}{1 + K \Delta t}  \hat{u}^n
  \end{aligned} \label{eq:implicit.evolution}
\end{equation}
The IMEX scheme~\eqref{eq:implicit.evolution} will show a decaying behavior with just the CFL condition $a \frac{\Delta t}{\Delta x} < 1$, while the explicit evolution~\eqref{eq:expl.evolution} will additionally require $\Delta t$ = $\mathcal{O}\left(\frac 1K\right)$.

\section{Discontinuous Galerkin and flux reconstruction methods for non-conservative hyperbolic equations} \label{app:dg.fr}
It is well known that the discontinuous Galerkin (DG) method can be cast into the flux reconstruction (FR) framework when the DG scheme uses the same quadrature nodes and solution points, and when the solution points and correction functions are chosen appropriately~\cite{Huynh2007,Grazia2014}.
In this section, we derive a DG method for hyperbolic systems with non-conservative products, and then show how it is equivalent to the FR scheme~\eqref{eq:general.non.conservative.fr.one}.
The hyperbolic system with non-conservative products is given by
\begin{equation}
\label{eq:non.cons.appendix}
\partial_t \uu + \partial_x \pf(\uu) + \bB (\uu)  \partial_x \uu =
   \bzero,
\end{equation}
where $\pf' + \bB$ is a diagonalizable matrix with real eigenvalues.
The DG method in the \textit{weak form} is derived by performing a \textit{formal integration by parts} with the Lagrange polynomial $\ell_p \in \poly_N$~\eqref{eq:defn.lagrange} to obtain
\begin{equation}\label{eq:weak.form.non.cons}
\begin{aligned}
& \int_{\Omega_e}\partial_t\uud_h\ell_p\ud x
+ \int_{\Omega_e}\partial_x\pf(\uud_h)\ell_p \ud x
+ \int_{\Omega_e}\bB(\uud_h)\partial_x\uud_h\ell_p\ud x
% = \partial_t\int_{\Omega_e}\uu\ell_p\ud x
% + \int_{\Omega_e}\pf(\uu)_x\ell_p \ud x
% + \int_{\Omega_e}\bB(\uu)\uu_x\ell_p\ud x
\\
&= \int_{\Omega_e}\partial_t\uud_h\ell_p\ud x
- \int_{\Omega_e} \discf \partial_x \ell_p \ud x
- \int_{\Omega_e}\partial_x(\bB(\uud_h)\ell_p)\uud_h\ud x
\\
& \qquad
+ (\fnum_\eph + \bB(\uu_\eph^-)\unum_\eph) \ell_p(x_\eph)
- (\fnum_\emh + \bB(\uu_{\emh}^+)\unum_\emh) \ell_p(x_\emh),
\end{aligned} \quad \forall p=0,\ldots,N,
\end{equation}
where $\uud_h$ denotes the discontinuous numerical solution~\eqref{eq:soln.poly}, $\discf$ denotes the discontinuous flux approximation~\eqref{eq:discts.flux}, $\fnum_\eph$ is the numerical flux defined in~\eqref{eq:combined.num.flux.fr}, whose central part approximates the conservative flux $\pf$ and the dissipative part consists of the wave speed estimate of the whole system, $\unum_\eph$ is the average of solution defined in~\eqref{eq:semi.fv}.
The volume flux term in~\eqref{eq:weak.form.non.cons} is given by the discontinuous flux approximation because the quadrature nodes are the same as the solution points (also see Section 2 of~\cite{babbar2025crk}).
We perform another formal\footnote{
The formal integration by parts step in~\eqref{eq:weak.form.non.cons} uses $\bB (\uu_{\eph^-}), \bB (\uu_{e - \half^+})$ instead of $\bB (\uu)_{\eph^-}, \bB (\uu)_{e - \half^+}$ as the latter will require additional inter-element communication of $\bB (\uu)_{\eph^-}, \bB (\uu)_{e - \half^+}$ instead of simply communicating $\uu_{\eph^{\pm}}$.
In order to maintain this consistency, the same $\bB (\uu_{\eph^-}), \bB (\uu_{e - \half^+})$ are used in~\eqref{eq:strong.form.non.cons}, which makes the integration by parts formal.
For the FR scheme with Gauss-Lobatto-Legendre points, this is irrelevant since $\bB (\uu_{\eph^-}) = \bB (\uu)_{\eph^-}$.} integration by parts on~\eqref{eq:weak.form.non.cons} to obtain the \textit{strong form} of the DG method~\eqref{eq:weak.form.non.cons}
\begin{equation}
\begin{aligned}
& \int_{\Omega_e} \partial_t \uud_h \ell_p \text{\textrm{d}} x
+ \int_{\Omega_e} \partial_x \discf \ell_p \text{\textrm{d}} x
+ \int_{\Omega_e } \bB (\uud_h) \ell_p \partial_x \uud_h \text{\textrm{d}} x
\\
& \qquad + \big [\fnum_\eph + \bB (\uu_\eph^-)\unum_\eph - \pf_\eph^- - \bB (\uu_\eph^-)\uu_\eph^-\big ] \ell_p (x_\eph) \\
& \qquad - \big[ \fnum_\emh + \bB (\uu_\emh^+)\unum_\emh - \pf_\emh^+ - \bB (\uu_\emh^+)\uu_\emh^+ \big ] \ell_p (x_\emh),
\end{aligned}
\qquad \forall p=0,\ldots,N.
\label{eq:strong.form.non.cons}
\end{equation}
We choose correction functions $g_{L/R} \in \mathbb{P}_{N + 1}$ to be $g_{\text{Radau}}, g_2$~{\cite{Huynh2007}} if the solution points are Gauss-Legendre (GL), Gauss-Lobatto-Legendre (GLL) respectively.
For both of the cases, we have the following identities from the proof of equivalence of FR and DG in~\cite{Huynh2007}\footnote{For Radau correction functions, the identities follow from taking $\phi = \ell_p$ in (7.8) of~{\cite{Huynh2007}}. For $g_2$, see Section 2.4 of~{\cite{Grazia2014}}. A detailed proof is also available in Appendix B of~\cite{babbar2024thesis}.}
\begin{equation}
g_R'(\xi_p) = \ell_p (1) / w_p, \qquad g_L' (\xi_p) = - \ell_p (0) / w_p.
\label{eq:correction.identities}
\end{equation}
Thus, by using the identity $\partial_x = \frac{1}{\Delta x_e} \partial_\xi$ obtained from the reference map~\eqref{eq:ref.map} and performing quadrature on the solution points, the strong form of the DG method~\eqref{eq:strong.form.non.cons} can be equivalently written as the FR method
\begin{equation}
\begin{aligned}
& \partial_t \uep
+ (\partial_x \discf)_{e,p}
+ \bB (\uep) (\partial_x \uud_h)_{e,p}
\\
& \qquad + \big[\fnum_\eph + \bB (\uu_\eph^-) \unum_\eph - \pf_\eph^- - \bB (\uu_\eph^-)\uu_\eph^-\big ] \partial_x g_R(x_\eph) \\
& \qquad + \big[\fnum_\emh + \bB (\uu_\emh^+) \unum_\emh - \pf_\emh^+ - \bB (\uu_\emh^+)\uu_\emh^+ \big ] \partial_x g_L(x_\emh)
% \qquad + \bB (\uu_\eph^-)(\unum_\eph -  \uu_\eph^-) \partial_x g_R(x_\eph) + \bB (\uu_\emh^+)  (\unum_\emh - \uu_\emh^+) \partial_x g_L(x_\emh)
= \bzero,
\end{aligned} \qquad \forall p=0,\ldots,N.
\label{eq:fr.form.non.cons}
\end{equation}
By using the definitions of non-conservative numerical fluxes~\eqref{eq:combined.num.flux.fr} and extrapolations~\eqref{eq:total.flux.non.cons}, it can directly be seen that the FR method~\eqref{eq:fr.form.non.cons} is the same as the FR method described in~\eqref{eq:general.non.conservative.fr.one}.
If we use the polynomial degree $N=0$, the volume terms in~\eqref{eq:fr.form.non.cons} vanish, the correction functions become linear functions so that $\partial_x g_R(x_\eph) = 1/\Delta x_e, \partial_x g_L(x_\emh) = -1/\Delta x_e$, leading to cancellation of the extrapolation terms $\pf_\emh^\pm - \bB (\uu_\emh^\pm)\uu_\emh^\pm$, and we obtain
\begin{equation}
\begin{aligned}
\partial_t \uez + \frac{1}{\Delta x_e} \big[(\fnum_\eph + \bB (\uu_\eph^-) \unum_\eph) - (\fnum_\emh + \bB (\uu_\emh^+) \unum_\emh)\big] = \bzero,
\end{aligned}
\label{eq:fr.fv}
\end{equation}
which is the first order finite volume scheme using the numerical fluxes~(\ref{eq:combined.num.flux},~\ref{eq:non.cons.num.flux.basic}).
The next specialization of the scheme~\eqref{eq:fr.form.non.cons} we now discuss is the one with GLL points and $g_2$ correction functions.
In this case, by~\eqref{eq:correction.identities}, the correction functions satisfy $g_R'(\xi_p) = \delta_{p,N} / w_p, g_L' (\xi_p) = - \delta_{p,0} / w_p$, where $\delta_{i,j}$ is the Kronecker delta function that is equal to $1$ if $i=j$ and $0$ otherwise.
This property implies that the correction functions only influence the solution for $p=0,N$.
In addition, noting that for GLL points, we have $\bB(\uu_\eph^-) = \bB(\ueN), \bB(\uu_\emh^+) = \bB(\uez)$, the FR scheme~\eqref{eq:fr.form.non.cons} with GLL points and $g_2$ correction functions can be written as
\begin{equation}
\begin{aligned}
& \partial_t \uep
+ \underbrace{(\partial_x \discf)_{e,p}
+ (\fnum_\eph - \pf_\eph^-) \partial_x g_R(x_\eph)
+ (\fnum_\emh - \pf_\emh^+) \partial_x g_L(x_\emh)}_{=(\dfrx \pf(\uud_h))_{e,p}~\eqref{eq:dfrx.defn}}
\\
& + \bB (\uep) \big[\underbrace{(\partial_x \uud_h)_{e,p} + (\unum_\eph - \uu_\eph^-)\partial_x g_R(x_\eph)  + (\unum_\emh - \uu_\emh^+)\partial_x g_L(x_\emh)}_{=(\dfrx \uud_h)_{e,p}~\eqref{eq:fr.non.cons.gll} }\big]
= \bzero,
\end{aligned} \quad \forall p=0,\ldots,N.
\label{eq:fr.gll}
\end{equation}
As can be seen from the identities under the braces, the scheme~\eqref{eq:fr.gll} is equivalent to~\eqref{eq:fr.non.cons.gll}.

\section{Examples of IMEX schemes} \label{app:imex.schemes}
In this section, we provide the specific Butcher tableaux of the IMEX schemes~\eqref{eq:imex.butcher} used in this work, while also giving the fully discrete IMEX compact Runge-Kutta flux reconstruction scheme for the cheapest second order IMEX method HT (1,1,2)~\cite{Hairer1996}.
IMEX methods follow the notation $\text{NAME} (s_E, s_I, p)$~\cite{boscarino2024} where $s_E, s_I,$ and $p$ represent, respectively, the number of stages of the explicit part, the number of stages of the implicit one, and the order of the IMEX-RK method.
The numbers $s_E, s_I$ need not be the same as $s$ because, as we will see in examples below, the initial rows of $A, \tilde{A}$ are often zero.

In terms of the structure of the implicit coefficient matrix $A$, IMEX-RK methods can be classified into two types known as type I (or type A) and type II (or type CK).
Since this classification did not influence the results of this work, we refer the reader to~\cite{boscarino2024} for the definitions.
% \begin{definition}
% An IMEX-RK method is said to be of type I (called also type A) if the matrix
% $A \in \mathbb{R}^{s \times s}$ is invertible, or equivalently $a_{ii} \neq 0$, $i = 1, \ldots, s$.
% An IMEX-RK method is said to be of type II
% (called also type CK) if the matrix $A$ can be written as
% \begin{equation}
% A = \left(\begin{array}{cc}
% 0 & 0\\
% \ba & \hat{A}
% \end{array}\right),
% \end{equation}
% with $\ba = (a_{21}, \ldots, a_{s 1})^T \in \mathbb{R}^{s - 1}$ and the sub-matrix $\hat{A} \in \mathbb{R}^{(s - 1) \times (s - 1)}$ with $\hat{a}_{ij} = a_{ij}$, $i, j = 2, \ldots, s$, is invertible.
% If the vector $\ba = \bzero$, then the IMEX-RK method is said to be of type ARS~\cite{ascher1997,ascher1995}.
% \end{definition}

The following terminologies are taken from~\cite{Hairer1996}.
The stability polynomial of the implicit part of the IMEX scheme is defined as
\begin{equation}
R(z) = \frac{\det(I - z A + z I b^T)}{\det(I - z A)}, \label{eq:stability.poly}
\end{equation}
where $I$ is the identity matrix of size $s \times s$.
The stability region is defined as $S = \{ z \in \mathbb{C} : |R(z)| \leq 1 \}$.
The method is said to be A-stable if the stability region contains the left half of the complex plane, i.e., $\{ z \in \mathbb{C} : \text{Re} (z) \leq 0 \} \subseteq S$.
If an A-stable method has the property that $\lim_{z \to \infty} R(z) = 0$, then it is said to be L-stable.
An implicit method is said to be stiffly accurate if the last row of the matrix $A$ is equal to the weights vector $b^T$, i.e., $a_{s j} = b_j$, $j = 1, \ldots, s$.
The IMEX method is said to be globally stiffly accurate (GSA)~\cite{boscarino2024} if the last row of both matrices $A$ and $\tilde{A}$ are equal to the respective weights vectors, i.e., $a_{s j} = b_j$, $\tilde{a}_{s j} = \tilde{b}_j$, $j = 1, \ldots, s$.
As proven in~\cite{Hairer1996}, if an A-stable implicit method is stiffly accurate, then it is also L-stable.

\paragraph{HT (1,1,2).}
The Butcher tableaux for the HT (1,1,2)~\cite{Hairer1991} scheme are given by
\begin{equation}
\begin{array}{c|cc}
0 & 0 & 0\\
1 & 1 & 0\\
\hline
& 1 / 2 & 1 / 2
\end{array} \qquad \begin{array}{c|cc}
0 & 0 & 0\\
1 & 1 / 2 & 1 / 2\\
\hline
& 1 / 2 & 1 / 2
\end{array}. \label{eq:ht112.butcher}
\end{equation}
Thus, the fully discrete IMEX cRKFR scheme~\eqref{eq:crkfr.inner}--\eqref{eq:crkfr} for the HT (1,1,2) scheme is given by
\begin{align*}
\uu^{(1)} & = \uu^n,\\
\uu^{(2)} & = \uu^n - \Delta t \dlocx  \pf(\uu^{(1)}) + \frac{\Delta t}{2}  (\bss (\uu^{(1)}) + \bss (\uu^{(2)})),\\
\uu^{n + 1} & = \uu^n - \Delta t \dfrx \F + \Delta t \bS_h^{\delta}, \\
&\F_h^\delta = \frac{1}{2}  (\pf_h^\delta (\uu^{(1)}) +\pf_h^\delta (\uu^{(2)}) ), \qquad \bS_h^\delta = \frac{1}{2}  (\bss(\uu^{(1)}) + \bss(\uu^{(2)})).
\end{align*}
The implicit part of the HT (1,1,2) scheme is A-stable but not L-stable, which we suspect to be the explanation for why this scheme did not work for Jin-Xin relaxation
(Section~\ref{sec:jin.xin}) with very small $\varepsilon$.
\paragraph{SSP3-IMEX(4,3,3).} This is one of the third order schemes from~\cite{pareschi2005} that we use
\begin{equation}
\begin{array}{c}
\begin{array}{c|cccc}
0 & 0 & 0 & 0 & 0\\
0 & 0 & 0 & 0 & 0\\
1 & 0 & 1 & 0 & 0\\
1 / 2 & 0 & 1 / 4 & 1 / 4 & 0\\
\hline
& 0 & 1 / 6 & 1 / 6 & 2 / 3
\end{array} \qquad \begin{array}{c|cccc}
\alpha & \alpha & 0 & 0 & 0\\
0 & 0 & \alpha & 0 & 0\\
1 & 1 / 3 & 1 / 3 & \alpha & 0\\
1 & \beta & \eta & 1 / 2 - \beta - \eta - \alpha & \alpha\\
\hline
& 0 & 1 / 6 & 1 / 6 & 2 / 3
\end{array} ,\\
\alpha = 0.241694260788, \quad \beta = 0.0604235651970, \quad \eta = 0.12915286960590.
\end{array}
\label{eq:ssp33.butcher}
\end{equation}
This scheme is L-stable, and we could use it for all the test cases including the test with arbitrary small $\varepsilon$ in the Jin-Xin relaxation system (Section~\ref{sec:jin.xin}).

% \paragraph{ARS-443.} This is the third order scheme from~\cite{ascher1997}, and has the following Butcher tableaus
% \begin{equation}
% \begin{array}{c|ccccc}
% 0 & 0 & 0 & 0 & 0 & 0\\
% 1 / 2 & 1 / 2 & 0 & 0 & 0 & 0\\
% 2 / 3 & 11 / 18 & 1 / 18 & 0 & 0 & 0\\
% 1 / 2 & 5 / 6 & - 5 / 6 & 1 / 2 & 0 & 0\\
% 1 & 1 / 4 & 7 / 4 & 3 / 4 & - 7 / 4 & 0\\
% \hline
% & 1 / 4 & 7 / 4 & 3 / 4 & - 7 / 4 & 0
% \end{array}, \qquad \begin{array}{c|ccccc}
% 0 & 0 & 0 & 0 & 0 & 0\\
% 1 / 2 & 0 & 1 / 2 & 0 & 0 & 0\\
% 2 / 3 & 0 & 1 / 6 & 1 / 2 & 0 & 0\\
% 1 / 2 & 0 & - 1 / 2 & 1 / 2 & 1 / 2 & 0\\
% 1 & 0 & 3 / 2 & - 3 / 2 & 1 / 2 & 1 / 2\\
% \hline
% & 0 & 3 / 2 & - 3 / 2 & 1 / 2 & 1 / 2
% \end{array} .
% \end{equation}
% The scheme is L-stable and globally stiffly accurate (GSA).
% The GSA property is possibly the reason why the scheme gave the least errors in the stiff collision source term test (Table~\ref{tab:imex.convergence}).
For Butcher tableaux of other IMEX schemes, we refer the reader to~\cite{boscarino2024}.
\printbibliography

\end{document}